\documentclass[12pt]{article}
\usepackage{amsfonts}
\usepackage{full page,
amssymb, amscd, amsmath,graphicx, 
}

 \title{{\bf The first cohomology, derivations and the reductivity of a (meromorphic open-string) vertex algebra}}
 \author{Yi-Zhi Huang and Fei Qi}
    \date{}
    \begin{document}
    \bibliographystyle{alpha}

\newtheorem{thm}{Theorem}[section]
\newtheorem{defn}[thm]{Definition}
\newtheorem{prop}[thm]{Proposition}
\newtheorem{cor}[thm]{Corollary}
\newtheorem{lemma}[thm]{Lemma}
\newtheorem{rema}[thm]{Remark}
\newtheorem{app}[thm]{Application}
\newtheorem{prob}[thm]{Problem}
\newtheorem{conv}[thm]{Convention}
\newtheorem{conj}[thm]{Conjecture}
\newtheorem{cond}[thm]{Condition}
    \newtheorem{exam}[thm]{Example}
\newtheorem{assum}[thm]{Assumption}
     \newtheorem{nota}[thm]{Notation}
\newcommand{\halmos}{\rule{1ex}{1.4ex}}
\newcommand{\pfbox}{\hspace*{\fill}\mbox{$\halmos$}}

\newcommand{\nn}{\nonumber \\}

 \newcommand{\res}{\mbox{\rm Res}}
 \newcommand{\ord}{\mbox{\rm ord}}
\renewcommand{\hom}{\mbox{\rm Hom}}
\newcommand{\edo}{\mbox{\rm End}\ }
 \newcommand{\pf}{{\it Proof.}\hspace{2ex}}
 \newcommand{\epf}{\hspace*{\fill}\mbox{$\halmos$}}
 \newcommand{\epfv}{\hspace*{\fill}\mbox{$\halmos$}\vspace{1em}}
 \newcommand{\epfe}{\hspace{2em}\halmos}
\newcommand{\nord}{\mbox{\scriptsize ${\circ\atop\circ}$}}
\newcommand{\wt}{\mbox{\rm wt}\ }
\newcommand{\swt}{\mbox{\rm {\scriptsize wt}}\ }
\newcommand{\lwt}{\mbox{\rm wt}^{L}\;}
\newcommand{\rwt}{\mbox{\rm wt}^{R}\;}
\newcommand{\slwt}{\mbox{\rm {\scriptsize wt}}^{L}\,}
\newcommand{\srwt}{\mbox{\rm {\scriptsize wt}}^{R}\,}
\newcommand{\clr}{\mbox{\rm clr}\ }
\newcommand{\tr}{\mbox{\rm Tr}}
\newcommand{\C}{\mathbb{C}}
\newcommand{\Z}{\mathbb{Z}}
\newcommand{\R}{\mathbb{R}}
\newcommand{\Q}{\mathbb{Q}}
\newcommand{\N}{\mathbb{N}}
\newcommand{\CN}{\mathcal{N}}
\newcommand{\F}{\mathcal{F}}
\newcommand{\I}{\mathcal{I}}
\newcommand{\V}{\mathcal{V}}
\newcommand{\one}{\mathbf{1}}
\newcommand{\BY}{\mathbb{Y}}
\newcommand{\ds}{\displaystyle}

        \newcommand{\ba}{\begin{array}}
        \newcommand{\ea}{\end{array}}
        \newcommand{\be}{\begin{equation}}
        \newcommand{\ee}{\end{equation}}
        \newcommand{\bea}{\begin{eqnarray}}
        \newcommand{\eea}{\end{eqnarray}}
         \newcommand{\lbar}{\bigg\vert}
        \newcommand{\p}{\partial}
        \newcommand{\dps}{\displaystyle}
        \newcommand{\bra}{\langle}
        \newcommand{\ket}{\rangle}

        \newcommand{\ob}{{\rm ob}\,}
        \renewcommand{\hom}{{\rm Hom}}

\newcommand{\A}{\mathcal{A}}
\newcommand{\Y}{\mathcal{Y}}

\newcommand{\dlt}[3]{#1 ^{-1}\delta \bigg( \frac{#2 #3 }{#1 }\bigg) }

\newcommand{\dlti}[3]{#1 \delta \bigg( \frac{#2 #3 }{#1 ^{-1}}\bigg) }

 \makeatletter
\newlength{\@pxlwd} \newlength{\@rulewd} \newlength{\@pxlht}
\catcode`.=\active \catcode`B=\active \catcode`:=\active
\catcode`|=\active
\def\sprite#1(#2,#3)[#4,#5]{
   \edef\@sprbox{\expandafter\@cdr\string#1\@nil @box}
   \expandafter\newsavebox\csname\@sprbox\endcsname
   \edef#1{\expandafter\usebox\csname\@sprbox\endcsname}
   \expandafter\setbox\csname\@sprbox\endcsname =\hbox\bgroup
   \vbox\bgroup
  \catcode`.=\active\catcode`B=\active\catcode`:=\active\catcode`|=\active
      \@pxlwd=#4 \divide\@pxlwd by #3 \@rulewd=\@pxlwd
      \@pxlht=#5 \divide\@pxlht by #2
      \def .{\hskip \@pxlwd \ignorespaces}
      \def B{\@ifnextchar B{\advance\@rulewd by \@pxlwd}{\vrule
         height \@pxlht width \@rulewd depth 0 pt \@rulewd=\@pxlwd}}
      \def :{\hbox\bgroup\vrule height \@pxlht width 0pt depth
0pt\ignorespaces}
      \def |{\vrule height \@pxlht width 0pt depth 0pt\egroup
         \prevdepth= -1000 pt}
   }
\def\endsprite{\egroup\egroup}
\catcode`.=12 \catcode`B=11 \catcode`:=12 \catcode`|=12\relax
\makeatother

\def\hboxtr{\FormOfHboxtr} 
\sprite{\FormOfHboxtr}(25,25)[0.5 em, 1.2 ex] 

:BBBBBBBBBBBBBBBBBBBBBBBBB | :BB......................B |
:B.B.....................B | :B..B....................B |
:B...B...................B | :B....B..................B |
:B.....B.................B | :B......B................B |
:B.......B...............B | :B........B..............B |
:B.........B.............B | :B..........B............B |
:B...........B...........B | :B............B..........B |
:B.............B.........B | :B..............B........B |
:B...............B.......B | :B................B......B |
:B.................B.....B | :B..................B....B |
:B...................B...B | :B....................B..B |
:B.....................B.B | :B......................BB |
:BBBBBBBBBBBBBBBBBBBBBBBBB |

\endsprite
\def\shboxtr{\FormOfShboxtr} 
\sprite{\FormOfShboxtr}(25,25)[0.3 em, 0.72 ex] 

:BBBBBBBBBBBBBBBBBBBBBBBBB | :BB......................B |
:B.B.....................B | :B..B....................B |
:B...B...................B | :B....B..................B |
:B.....B.................B | :B......B................B |
:B.......B...............B | :B........B..............B |
:B.........B.............B | :B..........B............B |
:B...........B...........B | :B............B..........B |
:B.............B.........B | :B..............B........B |
:B...............B.......B | :B................B......B |
:B.................B.....B | :B..................B....B |
:B...................B...B | :B....................B..B |
:B.....................B.B | :B......................BB |
:BBBBBBBBBBBBBBBBBBBBBBBBB |

\endsprite

\vspace{2em}



\date{}
\bibliographystyle{alpha}
\maketitle

\begin{abstract}
We give a criterion for the complete reducibility of modules  
satisfying a composability condition for a 
meromorphic open-string vertex algebra $V$ using the first cohomology of the algebra.  
For a $V$-bimodule $M$, let $\hat{H}^{1}_{\infty}(V, M)$ 
be the first cohomology of $V$ with the coefficients in $M$, which is canonically 
isomorphic to the 
quotient space of the space of derivations from $V$ to $W$ by the subspace of inner derivations. 
Let $\hat{Z}^{1}_{\infty}(V, M)$ be the subspace of $\hat{H}^{1}_{\infty}(V, M)$ canonically isomorphic 
to the space of  derivations obtained from the zero mode of the right vertex operators of weight $1$ 
elements such that the  difference between the skew-symmetric opposite action of the left action and  the
right action on these elements are Laurent polynomials in the variable.  If
$\hat{H}^{1}_{\infty}(V, M)= \hat{Z}^{1}_{\infty}(V, M)$ for every 
$\Z$-graded $V$-bimodule $M$, then every left $V$-module satisfying a composability condition 
is completely reducible.  
In particular, since a lower-bounded $\Z$-graded vertex algebra $V$ is a special meromorphic open-string
vertex algebra and left $V$-modules are in fact what has been 
called generalized $V$-modules with lower-bounded weights (or lower-bounded generalized $V$-modules), 
this result provides a cohomological criterion for the 
complete reducibility of lower-bounded generalized modules for such a vertex algebra. 
We conjecture that the converse of the main theorem above is also true. 
We also prove that when a grading-restricted vertex algebra $V$ contains a subalgebra
satisfying some familiar conditions, the composability condition for grading-restricted generalized $V$-modules 
always holds and we need $\hat{H}^{1}_{\infty}(V, M)= \hat{Z}^{1}_{\infty}(V, M)$ only for every $\Z$-graded
$V$-bimodule $M$ generated by a grading-restricted subspace
in our complete reducibility theorem. 
\end{abstract}

\renewcommand{\theequation}{\thesection.\arabic{equation}}
\renewcommand{\thethm}{\thesection.\arabic{thm}}
\setcounter{equation}{0}
\setcounter{thm}{0}
\section{Introduction}

In the representation theory of various algebras, one of the main tools is the cohomological
method. The powerful tool of homological 
algebra often provides a unified treatment of many results in 
representation theory. Such a unified treatment not only gives solutions 
to open problems, but also provides a 
conceptual understanding of the results. 

In the representation theory of vertex (operator) algebras, though the cohomology 
for a grading-restricted vertex algebra has been introduced by the first author in \cite{Hcoh},
the cohomological method 
in the representation theory of vertex (operator) algebras still needs 
to be fully developed. In this paper, we study the relation between the first 
cohomology and the complete reducibility of suitable modules for 
a grading-restricted vertex algebra or 
more generally for a meromorphic open-string vertex algebra introduced by the first author
in \cite{HMOSVA}.

In \cite{Hcoh} and \cite{H1st-sec-coh}, the first author 
introduced the cohomology of a grading-restricted vertex algebra and proved 
that the first cohomology and second cohomology indeed have the properties that 
they should have. The cohomology introduced in \cite{Hcoh} can be viewed
as an analogue of the Harrison cohomology of a commutative associative algebra.
In particular,  the first author also introduced in  \cite{Hcoh} a cohomology that should be
viewed as an analogue of the Hochschild cohomology of a commutative 
associative algebra viewed only as an associative algebra. The generalization 
of this cohomology to a meromorphic open-string vertex algebra (a noncommutative 
generalization of a grading-restricted vertex algebra) has been given by the second 
author in \cite{Q2} and \cite{Q3}. In this paper, we give a criterion for the complete 
reducibility of left modules for a
meromorphic open-string vertex algebra 
satisfying a composability condition  using the first cohomology of 
this analogue of the Hochschild cohomology.  In particular, we obtain a
criterion for such a module (in fact such a generalized module in the terminology 
of  \cite{HLZ}) for a lower-bounded or grading-restricted vertex algebra to be completely reducible.

We describe the main results of this paper more precisely here:
For a  meromorphic open-string vertex algebra
$V$ and a  $V$-bimodule $M$ (note that we do not require that the algebra $V$ and 
$V$-modules be grading restricted, though they are lower bounded), let $\hat{H}^{1}_{\infty}(V, M)$ 
be the first cohomology of $V$ with the coefficients in $M$ introduced in 
\cite{Hcoh}, \cite{Q2} and \cite{Q3}.  The first cohomology $\hat{H}^{1}_{\infty}(V, M)$  
is in fact canonically isomorphic to the quotient space of the space of derivations from $V$ to $M$ 
by the subspace of inner derivations.  Let $\hat{Z}^{1}_{\infty}(V, M)$ be the subspace of 
$\hat{H}^{1}_{\infty}(V, M)$ canonically isomorphic 
to the space of what we call zero-mode derivations, that is,
derivations obtained from the zero mode of the right vertex operators of weight $1$ 
elements of $M$ such that the  difference between the skew-symmetric opposite action of the left action and  the
right action on these elements are Laurent polynomials in the variable.
For a left $V$-module $W$, a left $V$-submodule
$W_{2}$ of $W$ and a graded subspace $W_{1}$ of $W$ such that as a graded vector space,
$W=W_{1}\oplus W_{2}$, let $\pi_{W_{1}}$ and $\pi_{W_{2}}$ be the projections
from $W$ to $W_{1}$ and $W_{2}$, respectively. For a left $V$-module $W$ and a left $V$-submodule
$W_{2}$,  we say that the pair
$(W, W_{2})$ satisfies the composability condition
if there exists a graded subspace $W_{1}$ of $W$ such that 
$W=W_{1}\oplus W_{2}$ and such that for $k, l\in \N$, $w'_{2}\in W_{2}$, $w_{1}\in W_{1}$, 
$v_{1}, \dots, v_{k+l}, v\in V$, the series
$$\langle w'_{2}, Y_{W_{2}}(v_{1}, z_{1})\cdots Y_{W_{2}}(v_{k}, z_{k})
\pi_{W_{2}}Y_{W}(v, z) \pi_{W_{1}}Y_{W}(v_{k+1}, z_{k+1})\cdots \pi_{W_{1}}Y_{W}(v_{k+l}, z_{k+l})w_{1}
\rangle$$ 
 is absolutely convergent the region $|z_{1}|>\cdots >|z_{k}|>|z|>\cdots >|z_{k+l}|>0$
 to a suitable rational function. We say that a left $V$-module
$W$ satisfies the composability condition if for every proper nonzero left 
$V$-submodule $W_{2}$ of $W$, the pair $(W, W_{2})$ satisfies the composability condition.
We prove in this paper that
if $\hat{H}^{1}_{\infty}(V, M)= \hat{Z}^{1}_{\infty}(V, M)$ for every 
 $\Z$-graded $V$-bimodule $M$,
then every left $V$-module satisfying the  composability condition is
completely reducible. Since the first cohomology of $V$ with coefficients in $M$
is the quotient of the space of derivations from $V$ to $M$ by the space of inner derivations,
the condition $\hat{H}^{1}_{\infty}(V, M)= \hat{Z}^{1}_{\infty}(V, M)$ 
in our main theorem above can also be formulated 
as the condition that every derivation from 
$V$ to $M$ is the sum of an inner derivation and a zero-mode derivation.

In particular, since a lower-bounded $\Z$-graded vertex algebra $V$ is a special meromorphic open-string
vertex algebra and left $V$-modules are in fact what has been 
called generalized $V$-modules with lower-bounded weights (or lower-bounded generalized $V$-modules), 
this result provides a cohomological criterion for the 
complete reducibility of lower-bounded generalized modules  for such a vertex algebra $V$. 
We also prove that when the grading-restricted vertex algebra $V$ contains a subalgebra $V_{0}$
such that products of intertwining operators are convergent absolutely in the usual region and the sums 
can be analytically extended, 
the associativity of intertwining operators holds and grading-restricted $V_{0}$-modules (grading-restricted 
generalized $V_{0}$-modules in the terminology of \cite{HLZ}) are completely reducible, 
the composability condition holds for every $V$-modules. 
We also prove that in this case, we need $\hat{H}^{1}_{\infty}(V, M)= \hat{Z}^{1}_{\infty}(V, M)$ 
only for every $\Z$-graded
$V$-bimodule $M$ generated by a grading-restricted subspace
in our complete reducibility theorem.

The composability condition on left $V$-modules needed in our results is in fact a condition on
convergence and analytic extension.
In the representation theory of vertex operator algebras, the convergence and analytic extension of suitable
series always play a crucial role. The main difficult parts of the proofs of a number of 
major results are in fact on suitable convergence and analytic extension. For example,  in the proof of the 
associativity of intertwining operators, the main difficult part of the proof in the paper  \cite{Hdiff-eqn}
is to use the $C_{1}$-cofiniteness of the modules to prove the convergence and extension property
of the products and iterates of intertwining operators. Another important example 
is the modular invariance of intertwining operators proved in \cite{Hint-mod-inv}. 
The most difficult part of the proof of this modular invariance is the proof of the 
convergence and analytic extension of the $q$-traces of products of at least two intertwining operators. 
In fact, the method of proving such convergence and analytic extension by reducing it using 
formal series recurrent relations to the convergence 
of $q$-traces of one vertex operator does not work for general intertwining operators. 
One needs to derive modular invariant differential equations satisfied by these $q$-traces 
directly from the $C_{2}$-cofiniteness of the modules. Note that the associativity of 
intertwining operators and the modular invariance of intertwining operators are the 
two main conjectures stated in the important paper of Moore and Seiberg \cite{MS} 
that led to the Verlinde formula and 
the modular tensor category structures (see
\cite{Hv-conj} and \cite{Hmod-cat} for the mathematical proof and mathematical construction, 
respectively, 
based on the associativity and modular invariance of intertwining operators
proved in \cite{Hdiff-eqn} and \cite{Hint-mod-inv}, respectively). 

In fact, in Section 6, we prove that every grading-restricted left $V$-module satisfies
the composability condition using the theory of intertwining operators 
for a subalgebra of a grading-restricted vertex algebra $V$ satisfying suitable conditions. 
This reveals a deep connection 
between the composability condition and the theory of intertwining operators. 
It is also proved in Section 6 that when $V$ is a grading-restricted vertex algebra containing a subalgebra 
satisfying suitable conditions,  the relevant $\Z$-graded $V$-bimodule constructed from 
a grading-restricted $V$-module and a $V$-submodule  appearing in the proofs of our results
is in fact generated by a grading-restricted subspace and thus we need consider only $\Z$-graded
$V$-bimodules generated by such subspaces. This proof uses the $Q(z)$-tensor product first introduced in 
\cite{HL1} and \cite{HL2} and studied further in the more general nonsemisimple case
in \cite{HLZ1.5}. This use also reveals a deep connection between the cohomology theory 
and the tensor category theory for module categories for a suitable vertex (operator) algebra. 
We expect that 
the study of the composability condition introduced in this paper 
and the $V$-bimodules constructed in Section 3 from 
two left $V$-modules  will lead to deep 
results in the cohomological method to the representation theory of (meromorphic open-string) vertex 
(operator) algebras.

The present paper is organized as follows: In Section 2, we recall the notions of meromorphic 
open-string vertex algebra, left module, right module
and bimodule for such an algebra. We also recall briefly the cohomology theory for such an algebra.
We prove that the first cohomology $\hat{H}^{1}_{\infty}(V, W)$ of such an algebra $V$ 
with coefficients in a bimodule $W$ is 
isomorphic to the quotient space of the space of derivations from $V$ to $W$ by the subspace of 
inner derivations. We introduce what we call zero-mode derivations and the subspace 
$\hat{Z}^{1}_{\infty}(V, W)$ of $\hat{H}^{1}_{\infty}(V, W)$. We construct $\Z$-graded $V$-bimodules from two  
left $V$-modules in Section 3. In Section 4, we construct a $1$-cocycle from a  
left $V$-module and a left $V$-submodule. This cocyle  is in fact the obstruction 
for the left $V$-module to be decomposed as the direct sum of the left $V$-submodule and 
another left $V$-submodule. We prove our main theorem on the complete 
reducibility in Section 5. In Section 6, 
we prove that  the composability condition holds for suitable modules for a grading-restricted vertex algebra
containing a subalgebra satisfying suitable conditions. We also prove in this section
that in this case, we need $\hat{H}^{1}_{\infty}(V, M)= \hat{Z}^{1}_{\infty}(V, M)$ 
only for every 
$V$-bimodule $M$ generated by a grading-restricted subspace
in the complete reducibility theorem.

\renewcommand{\theequation}{\thesection.\arabic{equation}}
\renewcommand{\thethm}{\thesection.\arabic{thm}}
\setcounter{equation}{0}
\setcounter{thm}{0}
\section{Meromorphic open-string vertex algebras, modules and cohomology}

We first recall the notion of meromorphic open-string vertex algebra. 

\begin{defn}
{\rm A {\it meromorphic open-string vertex algebra} is a $\Z$-graded vector space 
$V=\coprod_{n\in\Z} V_{(n)}$ (graded by {\it weights}) equipped with a {\it vertex operator map}
\begin{eqnarray*}
   Y_V:  V\otimes V &\to & V[[x,x^{-1}]]\\
	u\otimes v &\mapsto& Y_V(u,x)v
  \end{eqnarray*}
and a {\it vacuum} $\one\in V$, satisfying the following axioms:
\begin{enumerate}
\item Axioms for the grading:
(a) {\it Lower bound condition}: When $n$ is sufficiently negative,
$V_{(n)}=0$.
(b) {\it $\mathbf{d}$-commutator formula}: Let $\mathbf{d}_{V}: V\to V$
be defined by $\mathbf{d}_{V}v=nv$ for $v\in V_{(n)}$. Then
$$[\mathbf{d}_{V}, Y_{V}(v, x)]=\frac{d}{dx}Y_{V}(v, x)+Y_{V}(\mathbf{d}_{V}v, x)$$
for $v\in V$.


\item Axioms for the vacuum: (a) {\it Identity property}:
Let $1_{V}$ be the identity operator on $V$. Then
$Y_{V}(\mathbf{1}, x)=1_{V}$. (b)
{\it Creation property}: For $u\in V$, $Y_{V}(u, x)\mathbf{1}\in V[[x]]$ and 
$\lim_{x\to 0}Y_{V}(u, x)\mathbf{1}=u$.

\item {\it $D$-derivative property} and {\it $D$-commutator formula}:
Let $D_{V}: V\to V$ be the operator
given by
$$D_{V}v=\lim_{x\to 0}\frac{d}{dx}Y_{V}(v, x)\one$$
for $v\in V$. Then for $v\in V$,
$$\frac{d}{dx}Y_{V}(u, x)=Y_{V}(D_{V}u, x)=[D_{V}, Y_{V}(u, x)].$$

 \item {\it Rationality}: Let $V'=\coprod_{n\in \Z}V_{(n)}^*$ be the graded dual of $V$. 
For $u_1, \cdots, u_n, v\in V, v'\in V'$, the series
$$\langle v', Y_V(u_1, z_1)\cdots Y_V(u_n, z_n)v\rangle$$
converges absolutely when $|z_1|>\cdots>|z_n|>0$ to a rational function in $z_1,\cdots, z_n$, 
with the only possible poles at $z_i=0, i=1,\cdots, n$ and $z_i=z_j, 1\leq i\neq j\leq n$. 
For $u_1, u_2, v\in V$ and $v' \in V'$, the series
$$\langle v', Y_V(Y_V(u_1,z_1-z_2)u_2,z_2)v\rangle$$
converges absolutely when $|z_2|>|z_1-z_2|>0$ to a rational function with the only possible poles at 
$z_1=0, z_2=0$ and $z_1=z_2$.

\item {\it Associativity}: For $u_{1}, u_{2}, v\in V$ and
$v'\in V'$, we have 
$$\langle v', Y_{V}(u_{1}, z_{1})Y_{V}(u_{2}, z_{2})v\rangle=\langle v', Y_{V}(Y_{V}(u_{1}, z_{1}-z_{2})u_{2}, z_{2})v\rangle$$
when $|z_{1}|>|z_{2}|>|z_{1}-z_{2}|>0$.

\end{enumerate} A meromorphic open-string vertex algebra is said to be {\it grading restricted} if 
$\dim V_{(n)}<\infty$ for $n\in \Z$. }
\end{defn}

The meromorphic open-string vertex algebra just defined is denoted by $(V, Y_V, \one)$ or simply by $V$. 

\begin{rema}\label{va-defn}
{\rm Let $V$ be a meromorphic open-string vertex algebra. We say that $V$ satisfies the {\it skew-symmetry}
if $Y_{V}(u, x)v=e^{xD_{V}}Y(v, -x)u$ for $u, v\in V$. A meromorphic open-string vertex algebra 
satisfies the skew-symmetry is a lower-bounded vertex algebra, that is, a vertex algebra with a lower-bounded 
$\Z$-grading. A grading-restricted meromorphic open-string vertex algebra satisfying the skew-symmetry 
is a grading-restring vertex algebra, that is, a vertex algebra with a $\Z$-grading satisfying the two 
grading restriction conditions.  }
\end{rema}

The following notion of left $V$-module  was introduced 
in \cite{HMOSVA}:

\begin{defn}
{\rm A {\it left module for $V$} or a {\it left $V$-module} is a $\C$-graded vector space 
$W=\coprod_{n\in \C}W_{[n]}$ (graded by {\it weights}), equipped with 
a {\it vertex operator map}
\begin{eqnarray*}
Y_{W}: V\otimes W&\to& W[[x, x^{-1}]]\nn
u\otimes w&\mapsto& Y_{W}(u, x)w,
\end{eqnarray*}
and an operator $D_{W}$ of weight $1$, satisfying the 
following axioms:
\begin{enumerate}

\item Axioms for the grading: (a) {\it Lower bound condition}:  When the real part $\Re{(n)}$ 
of $n$ is sufficiently negative,
$W_{[n]}=0$. (b) {\it $\mathbf{d}$-commutator formula}: Let $\mathbf{d}_{W}: W\to W$
be defined by $\mathbf{d}_{W}w=nw$ for $w\in W_{[n]}$. Then for $u\in V$, 
$$[\mathbf{d}_{W}, Y_{W}(u,x)]= Y_{W}(\mathbf{d}_{V}u,x)+x\frac{d}{dx}Y_{W}(u,x).$$

\item The {\it identity property}:
$Y_{W}(\one,x)=1_{W}$.

\item The {\it $D$-derivative property} and the  {\it $D$-commutator formula}: 
For $u\in V$,
\begin{eqnarray*}
\frac{d}{dx}Y_{W}(u, x)
&=&Y_{W}(D_{V}u, x) \nn
&=&[D_{W}, Y_{W}(u, x)].
\end{eqnarray*}

\item {\it Rationality}: For $u_{1}, \dots, u_{n}, w\in W$
and $w'\in W'$, the series 
$$
\langle w', Y_{W}(u_{1}, z_1)\cdots Y_{W}(u_{n}, z_n)w\rangle
$$
converges absolutely 
when $|z_1|>\cdots >|z_n|>0$ to a rational function in $z_{1}, \dots, z_{n}$
with the only possible poles at $z_{i}=0$ for $i=1, \dots, n$ and $z_{i}=z_{j}$ 
for $i\ne j$. For $u_{1}, u_{2}, w\in W$
and $w'\in W'$, the series 
$$
\langle w', Y_{W}(Y_{V}(u_{1}, z_1-z_{2})u_{2}, z_2)w\rangle
$$
converges absolutely when $|z_{2}|>|z_{1}-z_{2}|>0$ to a rational function
with the only possible poles at $z_{1}=0$, $z_{2}=0$ and $z_{1}=z_{2}$. 

\item {\it Associativity}: For $u_{1}, u_{2}, w\in W$, 
$w'\in W'$, 
$$
\langle w', 
Y_{W}(u_{1},z_1)Y_{W}(u_{2},z_2)w\rangle
=
\langle w', 
Y_{W}(Y_{V}(u_{1},z_{1}-z_{2})u_{2},z_2)w\rangle
$$
when  $|z_{1}|>|z_{2}|>|z_{1}-z_{2}|>0$.

\end{enumerate} 

A left $V$-module is said to be {\it grading restricted} if 
$\dim W_{[n]}<\infty$ for $n\in \C$. }
\end{defn}

We denote the left $V$-module just defined by $(W, Y_{W},  D_{W})$ or simply $W$. 

The following notions of right $V$-module and $V$-bimodule were introduced by the first author and were
explicitly written down in \cite{Q1} and \cite{Q3} by the second author:

\begin{defn}
{\rm A {\it right module for $V$} or a {\it right $V$-module} is a $\C$-graded vector space 
$W=\coprod_{n\in \C}W_{[n]}$ (graded by {\it weights}), equipped with 
a {\it vertex operator map}
\begin{eqnarray*}
Y_{W}: W\otimes V&\to& W[[x, x^{-1}]]\nn
w\otimes u&\mapsto& Y_{W}(w, x)u,
\end{eqnarray*}
and an operator $D_{W}$ of weight $1$, satisfying the 
following axioms:
\begin{enumerate}

\item Axioms for the grading: (a) {\it Lower bound condition}:  When $\Re{(n)}$ is sufficiently negative,
$W_{[n]}=0$. (b) {\it $\mathbf{d}$-commutator formula}: Let $\mathbf{d}_{W}: W\to W$
be defined by $\mathbf{d}_{W}w=nw$ for $w\in W_{[n]}$. Then for $w\in W$, 
$$\mathbf{d}_{W}Y_{W}(w,x)-Y_{W}(w,x)\mathbf{d}_{V}= Y_{W}(\mathbf{d}_{W}w,x)+x\frac{d}{dx}Y_{W}(w,x).$$

\item The {\it creation property}: For $w\in W$, $Y_{W}(w,x)\one\in W[[x]]$ and 
$\lim_{x\to 0}Y_{W}(w,x)\one=w$.

\item The {\it $D$-derivative property} and the  {\it $D$-commutator formula}: 
For $u\in V$,
\begin{eqnarray*}
\frac{d}{dx}Y_{W}(w, x)
&=&Y_{W}(D_{W}w, x) \nn
&=&D_{W}Y_{W}(w, x)-Y_{W}(w, x)D_{V}.
\end{eqnarray*}

\item {\it Rationality}: For $u_{1}, \dots, u_{n}, w\in W$
and $w'\in W'$, the series 
$$
\langle w', Y_{W}(w, z_1)Y_{V}(u_{1}, z_2)\cdots Y_{V}(u_{n-1}, z_n)u_{n}\rangle
$$
converges absolutely 
when $|z_1|>\cdots >|z_n|>0$ to a rational function in $z_{1}, \dots, z_{n}$
with the only possible poles at $z_{i}=0$ for $i=1, \dots, n$ and $z_{i}=z_{j}$ 
for $i\ne j$. For $u_{1}, u_{2}, w\in W$
and $w'\in W'$, the series 
$$
\langle w', Y_{W}(Y_{W}(w, z_1-z_{2})u_{1}, z_2)u_{2}\rangle
$$
converges absolutely when $|z_{2}|>|z_{1}-z_{2}|>0$ to a rational function
with the only possible poles at $z_{1}=0$, $z_{2}=0$ and $z_{1}=z_{2}$. 

\item {\it Associativity}: For $u_{1}, u_{2}, w\in W$, 
$w'\in W'$, 
$$
\langle w', 
Y_{W}(w, z_1)Y_{V}(u_{1},z_2)u_{2}\rangle
=
\langle w', 
Y_{W}(Y_{W}(w,z_{1}-z_{2})u_{1},z_2)u_{2}\rangle
$$
when  $|z_{1}|>|z_{2}|>|z_{1}-z_{2}|>0$.

\end{enumerate} 

A right $V$-module is said to be {\it grading restricted} if 
$\dim W_{[n]}<\infty$ for $n\in \C$. }
\end{defn}

We also denote the generalized right $V$-module just defined by $(W, Y_{W}, D_{W})$ or simply $W$. 

We now give the definition of $V$-bimodule. Roughly speaking, a 
$V$-bimodule is a $\C$-graded vector space equipped with a left $V$-module structure and a 
right-module structure such that these two structures are compatible. More precisely, we have:

\begin{defn}
{\rm A {\it 
$V$-bimodule} is a $\C$-graded vector space 
$W=\coprod_{n\in \C}W_{[n]}$
(graded by {\it weights})
equipped with a {\it left vertex operator map}
\begin{eqnarray*}
Y_{W}^{L}: V\otimes W&\to& W[[x, x^{-1}]]\nn
u\otimes w&\mapsto& Y_{W}^{L}(u, x)w,
\end{eqnarray*}
a {\it right vertex operator map}
\begin{eqnarray*}
Y_{W}^{R}: W\otimes V&\to& W[[x, x^{-1}]]\nn
w\otimes u&\mapsto& Y_{W}^{R}(w, x)u,
\end{eqnarray*}
and $D_W$ on $W$ satisfying the following conditions.
\begin{enumerate}

 \item $(W, Y_W^L, D_W)$ is a left $V$-module.

 \item $(W, Y_W^R, D_W)$ is a  right $V$-module.

 \item  {\it Compatibility}: 

\begin{enumerate}

\item {\it Rationality for left and right vertex operator maps}: For $u_{1}, \dots, u_{m+n}\in V$, $w\in W$ and 
$w'\in W'$, the series
$$\langle w', Y_{W}^{L}(u_{1}, z_{1})\cdots Y_{W}^{L}(u_{k}, z_{m})Y_{W}^{R}(w, z_{m+1})
Y_{V}(v_{m+1}, z_{m+2})\cdots Y_{V}(v_{m+n-1}, z_{m+n})v_{m+n}\rangle$$
is  absolutely convergent in the region $|z_{1}| > \cdots |z_{m+n}|>0$ to a
rational function in $z_{1}, \dots, z_{m+n}$ with the only possible poles $z_{i}=0$
for $i=1, \dots, m+n$ and $z_{i}=z_{j}$ for $i, j=1, \dots, m+n$, $i\ne j$. For 
$u, v\in V$, $w\in W$ and $w'\in W'$, the series
$$\langle w', Y_W^R(Y_W^L(u,z_1-z_2)w, z_2)v\rangle$$
converges absolutely when $|z_{2}|>|z_{1}-z_{2}|>0$ to a rational function
with the only possible poles at $z_{1}=0$, $z_{2}=0$ and $z_{1}=z_{2}$.

\item {\it Associativity for left and right vertex operator maps}:
For $u, v\in V$, $w\in W$ and $w'\in W'$, 
$$\langle w', Y_W^L(u, z_1)Y_W^R(w,z_2)v\rangle=\langle w', Y_W^R(Y_W^L(u,z_1-z_2)w, z_2)v\rangle$$
when  $|z_{1}|>|z_{2}|>|z_{1}-z_{2}|>0$. 
\end{enumerate}
 
\end{enumerate}}
\end{defn}

The $V$-bimodule just defined is denoted 
by  $(W, Y_W^L, Y_W^R, D_W)$ or simply by $W$. There is also a notion
of $V$-bimodule with different $D_{W}$ for the left and right
$V$-module structure. Since we do not need such $V$-bimodules in this paper, we shall not 
discuss this more general notion.

We shall use the notations and terminology in \cite{Hcoh} and \cite{H-2-cons}. For example, 
we shall denote 
the  rational function that $\langle w', Y_W^L(u, z_1)Y_W^R(w,z_2)v\rangle$ converges to 
by 
$$R(\langle w', Y_W^L(u, z_1)Y_W^R(w,z_2)v\rangle).$$
In particular, we have, for example, 
$$R(\langle w', Y_W^L(u, z_1)Y_W^R(w,z_2)v\rangle)=R(\langle w', Y_W^R(Y_W^L(u,z_1-z_2)w, z_2)v\rangle).$$
We shall also use, for example, $Y_W^L(u, z_1)Y_W^R(w,z_2)v$ to denote the element 
of $\overline{W}\subset (W')^{*}$ given by
$w' \mapsto \langle w', Y_W^L(u, z_1)Y_W^R(w,z_2)v\rangle$ for $z_{1}$ and $z_{2}$
satisfying $|z_1|>|z_2|>0$ (see Remark \ref{alg-compl}). 
In particular, we shall say, for example,  that $Y_W^L(u, z_1)Y_W^R(w,z_2)v$ is equal to 
$Y_W^R(Y_W^L(u,z_1-z_2)w, z_2)v$   in the region $|z_1|>|z_2|>|z_{1}-z_{2}|>0$. 
In addition, we shall use, for example,
$E(Y_W^L(u, z_1)Y_W^R(w,z_2)v)$ to denote the the element 
of $\overline{W}$ given by
$w' \mapsto R(\langle w', Y_W^L(u, z_1)Y_W^R(w,z_2)v\rangle)$ for $z_{1}, z_{2}\in \C^{\times}$
satisfying $z_{1}\ne z_{2}$. In particular, we have, for example,
$$E(Y_W^L(u, z_1)Y_W^R(w,z_2)v)=E(Y_W^R(Y_W^L(u,z_1-z_2)w, z_2)v).$$


The following pole-order conditions on a meromorphic open-string vertex algebra $V$
and various $V$-modules are introduced in \cite{Q2} and \cite{Q3}:

\begin{defn}
{\rm A meromorphic open-string vertex algebra $V$ is said to satisfy the {\it pole-order condition}
if for $v_{1}, \dots, v_{n}, u\in V$, there 
exist $r_{i}\in \N$ depending only on the pair $(v_{i}, u)$ 
for $i=1, \dots, n$, $p_{ij} \in \N$ depending only on the pair  $(v_{i}, v_{j})$ 
for $i, j=1, \dots, n$, $i\ne j$ 
and $g(z_{1}, \dots, z_{n})\in V[[z_{1}, \dots, z_{n}]]$ such that for $u'\in V'$, 
$$\prod_{i=1}^{n}z_{i}^{r_{i}}\prod_{1\le i<j\le n}(z_{i}-z_{j})^{p_{ij}}
R(\langle u', Y_{V}(v_{1}, z_{1})
\cdots Y_{V}(v_{n}, z_{n})u\rangle)$$
is a polynomial and is equal to 
$\langle u', g(z_{1}, \dots, z_{n})\rangle$. A left $V$-module $W$
is said to satisfy the {\it pole-order condition}
if for $v_{1}, \dots, v_{n}\in V$ and $w\in W$, there 
exist $r_{i}\in \N$ depending only on the pair $(v_{i}, w)$
for $i=1, \dots, n$, $p_{ij} \in \N$ depending only on the pair  $(v_{i}, v_{j})$
for $i, j=1, \dots, n$, $i\ne j$ 
and $g(z_{1}, \dots, z_{n})\in W[[z_{1}, \dots, z_{n}]]$ such that for $w'\in W'$, 
$$\prod_{i=1}^{n}z_{i}^{r_{i}}\prod_{1\le i<j\le n}(z_{i}-z_{j})^{p_{ij}}
R(\langle w', Y_{W}(v_{1}, z_{1})
\cdots Y_{W}(v_{n}, z_{n})w\rangle)$$
is a polynomial and is equal to 
$\langle w', g(z_{1}, \dots, z_{n})\rangle$. A right $V$-module $W$
is said to satisfy the {\it pole-order condition}
if for $v_{1}, \dots, v_{n}\in V$ and $w\in W$, there 
exist $r_{1}\in \N$ depending only on the pair $(w, v_{n})$, 
$r_{i}\in \N$ depending only on the pair $(v_{i-1}, v_{n})$
for $i=2, \dots, n$, $p_{1j} \in \N$ depending only on the pair  $(w, v_{j-1})$
for $j=2, \dots, n$, $p_{ij}\in \N$ depending only on the pair  $(v_{i-1}, v_{j-1})$
for $i, j=2, \dots, n$, $i\ne j$ 
and $g(z_{1}, \dots, z_{n})\in W[[z_{1}, \dots, z_{n}]]$ such that for $w'\in W'$, 
$$\prod_{i=1}^{n}z_{i}^{r_{i}}\prod_{1\le i<j\le n}(z_{i}-z_{j})^{p_{ij}}
R(\langle w', Y_{W}(w, z_{1})Y_{V}(v_{1}, z_{2})
\cdots Y_{V}(v_{n-1}, z_{n})v_{n}\rangle)$$
is a polynomial and is equal to 
$\langle w', g(z_{1}, \dots, z_{n})\rangle$. A $V$-bimodule $W$
is said to satisfy the {\it pole-order condition}
if for $v_{1}, \dots, v_{k+l}\in V$ and $w\in W$, there exist
$r_{i}\in \N$ depending only on the pair $(v_{i}, u)$
for $i=1, \dots, k+l$, $m\in \N$ depending on the pair $(w, v_{n})$, $p_{ij} \in \N$ depending only 
on the pair $(v_{i}, v_{j})$
for $i, j=1, \dots, k+l$, $i\ne j$, $s_{i}\in \N$ depending only on the pair $(v_{i}, w)$ for $i=1, \dots, k+l$
and $g(z_{1}, \dots, z_{k+l}, z)\in W[[z_{1}, \dots, z_{k+l}]]$ such that for $w'\in W'$, 
\begin{align*}
z^{m}\prod_{i=1}^{k+l}&z_{i}^{r_{i}}\prod_{1\le i<j\le k+l}(z_{i}-z_{j})^{p_{ij}}
\prod_{i=1}^{k+l}(z_{i}-z)^{s_{i}}\cdot\nn
&\cdot R(\langle w', Y_{W}^{L}(v_{1}, z_{1})
\cdots Y_{W}^{L}(v_{k}, z_{k})Y_{W}^{R}(w, z)
Y_{V}(v_{k+1}, z_{k+1})\cdots Y_{V}(v_{k+l}, z_{k+l})u\rangle)
\end{align*}
is a polynomial and is equal to 
$\langle w', g(z_{1}, \dots, z_{k+l}, z)\rangle$. }
\end{defn}

\begin{rema}\label{alg-compl}
{\rm In the definition above,  $p_{ij}$ are required to be dependent only on the pair  $(v_{i}, v_{j})$.
This condition is automatically satisfied by lower-bounded vertex algebras and lower-bounded 
modules for such vertex algebras because of the commutativity satisfied by them. 
This condition allows us to formulate equivalent definitions of meromorphic open-string vertex algebras
and modules using formal variables. On the other hand, the condition that there exists 
$g(z_{1}, \dots, z_{n})$ guarantees that the product of vertex operators 
acting on an element of the algebra or module is in its algebraic completion, which is 
smaller than the full dual space of the graded dual of the algebra or module when 
its homogeneous subspaces are not finite dimensional. For example, 
from this existence,  it can  be shown easily that the element of $(V')^{*}$ given by 
$v'\mapsto \langle v', Y_{V}(v_{1}, z_{1})
\cdots Y_{V}(v_{n}, z_{n})u\rangle$ for $v'\in V'$ is in fact in $\overline{V}$. Similarly for modules. See \cite{Q2} and \cite{Q3}
for detailed discussions.}
\end{rema}


In this paper, we consider only those meromorphic open-string vertex algebras, 
left modules, right modules and bimodules satisfying these pole-order conditions. 
From now on, when we mention such algebras and modules, we shall omit the phrase 
``satisfying the pole-order condition."

\begin{rema}
{\rm If $V$ is a lower-bounded vertex algebra, then a left $V$-module (or a right-$V$-module) $W$
has a $V$-bimodule structure whose right vertex operator map (or left vertex operator map) 
is defined by $Y_{W}^{R}(w, x)v=e^{xD_{W}}Y_{W}^{L}(v, -x)w$ for $v\in V$ and $w\in W$
(or $Y_{W}^{L}(v, x)w=e^{xD_{W}}Y_{W}^{R}(w, -x)v$ for $v\in V$ and $w\in W$. The proof was
in fact given in \cite{FHL}.}
\end{rema}

We now briefly review the cohomology $\hat{H}^{\bullet}_{\infty}(V, W)$ for 
a meromorphic open-string vertex algebra $V$ and a $V$-bimodule $W$. 
We refer the reader to \cite{Hcoh}, \cite{Q2} and \cite{Q3} for details.

Let $V$ be a meromorphic open-string vertex algebra and $W$ a $V$-bimodule. 
Note that since we do not assume that $W$ is grading restricted, 
$\overline{W}=\prod_{n\in \\C}W_{[n]}$ might not be isomorphic to 
$(W')^{*}$. 
For $n>0$, a map $f$ from the configuration space $F_{n}\C=\{(z_{1}, \dots, z_{n})\in \C^{n}
\;|\; z_{i}\ne z_{j}, i, j=1, \dots, n, \;i\ne j\}$ to $\overline{W}=\prod_{n\in \C}W_{[n]}$
is called a $\overline{W}$-valued 
rational functions in $z_{1}, \dots, z_{n}$ with the only possible poles $z_{i}=z_{j}$ for 
$i, j=1, \dots, n$, $i\ne j$ if for $w'\in W'$, 
$\langle w', f(z_{1}, \dots, z_{n})\rangle$ is a rational function of such a form. 
For example, for $v_{1}, \dots, v_{n}\in V$, 
$E(Y_{V}(v_{1}, z_{1})\cdots Y_{V}(v_{n}, z_{n})\one)$
is a $\overline{V}$-valued 
rational functions in $z_{1}, \dots, z_{n}$ with the only possible poles $z_{i}=z_{j}$ for 
$i, j=1, \dots, n$, $i\ne j$.
Let $\widetilde{W}_{z_{1}, \dots, z_{n}}$ be the space of such $\overline{W}$-valued 
rational functions.

For $n>0$,  let $\hat{C}^{n}_{\infty}(V, W)$ be the subspace 
of $\hom(V^{\otimes n}, \widetilde{W}_{z_{1}, \dots, z_{n}})$
satisfying the $D$-derivative property, the $\mathbf{d}$-conjugation property and 
being composable with arbitrary numbers of vertex operators. In this paper we are mainly 
concerned with the first cohomology. So we shall  give the definitions of the $D$-derivative property, 
the $\mathbf{d}$-conjugation property and 
being composable with an arbitrary number of vertex operators only  in the case of $n=1$. In the general case,
see \cite{Hcoh}, \cite{Q2} and \cite{Q3} for the precise meaning of these properties. 
An element $\Psi$ of $\hom(V, \widetilde{W}_{z_{1}})$ is said to satisfy
the $D$-derivative property if for $v_{1}\in V$ and $z_{1}\in \C$, 
$$\frac{d}{dz_{1}} (\Psi(v_{1}))(z_{1})=(\Psi(D_{V}v_{1}))(z_{1})
=D_{W}(\Psi(v_{1}))(z_{1})$$
and is said to satisfy the $\mathbf{d}$-conjugation property if for $a\in \C^{\times}$, 
$v_{1}\in V$ and $z_{1}\in \C$, 
$$a^{\mathbf{d}_{W}}(\Psi(v_{1}))(z_{1})=(\Psi(a^{\mathbf{d}_{V}}v_{1}))(az_{1}).$$
Such an element is said to be composable with an arbitrary number of  vertex operators if 
for $k, l, m\in \N$, $v_{1}, \dots, v_{k+l+m}\in V$, $w'\in W'$, 
\begin{align*}
\langle w', Y^{L}_{W}(v_{1}, z_{1})\cdots Y^{L}_{W}(v_{k}, &z_{k})
 Y^{R}_{W}(\Psi(Y_{V}(v_{k+1}, z_{k+1}-\xi)\cdots 
Y_{V}(v_{k+l}, z_{k+l}-\xi)\one)(\xi-\zeta), \zeta)\cdot\nn
&\cdot Y_{V}(v_{k+l+1}, z_{k+l+1})
\cdots Y_{V}(v_{k+l+m}, z_{k+l+m})\one\rangle
\end{align*}
is absolutely convergent in the region given by
$|z_{1}|>\cdots >|z_{k}|$,  $|z_{k+1}-\xi|>\cdots >|z_{k+l}-\xi|$,
$|z_{k+l+1}|>\cdots > |z_{k+l+m}|$,
$|z_{i}|>|z_{k+j}-\xi|+|\xi-\zeta|+|\zeta|$
and $|\zeta|>|z_{k+j}-\xi|+|z_{k+l+p}|+|\xi-\zeta|$ for  $i=1, \dots, k$,
 $j=1, \dots, l$ and $p=1, \dots, m$
to a rational function of the form $\langle w', f(z_{1}, \dots, z_{k+l+m})\rangle$,
where $f$ is a $\overline{W}$-valued 
rational function  in $z_{1}, \dots, z_{k+l+m}$,
with the only possible poles $z_{i}=z_{j}$ for $i, j=1, \dots, k+l+m$, $i\ne j$. 
Moreover, we also require that there 
exist $p_{ij}\in \N$ depending only on $v_{i}$ and $v_{j}$ for $i, j=1, \dots, k+l+m$, $i\ne j$,
and $g(z_{1}, \dots, z_{k+l+m})\in W[[z_{1}, \dots, z_{k+l+m}]]$ such that for $w'\in W'$, 
$$\prod_{1\le i<j\le k+l+m}(z_{i}-z_{j})^{p_{ij}}\langle w', f(z_{1}, \dots, z_{k+l+m})\rangle$$
is a polynomial and is equal to 
$\langle w', g(z_{1}, \dots, z_{k+l+m})\rangle$. Note that the last part of the condition 
of being composable with an arbitrary number of  vertex operators implies that 
the order of the pole $z_{i}=z_{j}$ is bounded by a number depending only on
$v_{i}$ and $v_{j}$. Using the $\mathbf{d}$-conjugation property and this last 
part of the condition of being composable with an arbitrary number of  vertex operators, 
we can also show that the expansions of $f(z_{1}, \dots, z_{k+l+m})$ in various regions 
are in fact in $\overline{W}$ instead of just $(W')^{*}$. Note that when $W$ is not grading restricted,
$\overline{W}\ne (W')^{*}$. 

The formulation of the property above that an element of  $\hom(V, \widetilde{W}_{z_{1}})$
is composable with an arbitrary number of  vertex operators
uses directly the left and right
vertex operators. In \cite{Q2} and \cite{Q3}, the second author uses the skew-symmetry 
opposite vertex operators instead of the right vertex operators. It is easy to see 
from the definitions that the formulation here and the formulation in \cite{Q2} and \cite{Q3} are the same. 

For $n=0$, let 
$\hat{C}^{0}_{\infty}(V, W)$ be the subspace of $W$ (which is in fact 
isomorphic to $\hom(\C, W)$) consisting of elements, say $w$, such that 
$a^{\mathbf{d}_{W}}w=w$ (that is, $w\in W_{[0]}$ and can be interpreted as
a version of the $\mathbf{d}$-conjugation property),
$D_{W}w=0$ (can be interpreted as
a version of the $D$-derivative property).
Note that because $D_{W}w=0$, the $D$-derivative property 
of $Y_{W}^{R}$ implies that for $v\in V$, $Y_{W}^{R}(w, -z)v$ is independent of $z$ and
hence $e^{zD_{W}}Y_{W}^{R}(w, -z)v$ is a power series in $z$. 
Also, using the $D$-derivative property, the $D$-commutator formula and 
$D_{W}w=0$, we see that for $v\in V$,  the derivative of $e^{-zD_{W}}Y_{W}^{L}(v, z)w$ with respect to $z$
is in fact $0$. Thus for $v\in V$, $Y_{W}^{L}(v, z)w=e^{zD_{W}}(e^{-zD_{W}}Y_{W}^{L}(v, z)w)$ is 
a power series in $z$. An element $w$ of $\hat{C}^{0}_{\infty}(V, W)$ has almost all the properties 
that the vacuum $\one\in V$ has and can be called a vacuum-like element.

For $n>0$ and  $\Psi\in \hat{C}^{n}_{\infty}(V, W)$,
we define $\hat{\delta}_{\infty}^{n} \Psi\in  \hat{C}^{n+1}_{\infty}(V, W)$ by
\begin{eqnarray}\label{def-delta}
\lefteqn{\langle w', ((\hat{\delta}^{n}_{\infty}(\Psi))(v_{1}\otimes \cdots\otimes v_{n+1}))
(z_{1}, \dots, z_{n+1})\rangle}\nn
&&=R(\langle w', Y_{W}^{L}(v_{1}, z_{1})(\Psi(v_{2}\otimes \cdots\otimes v_{n+1}))
(z_{2}, \dots, z_{n+1})\rangle)\nn
&&\quad +\sum_{i=1}^{n}(-1)^{i}R(\langle w', 
(\Psi(v_{1}\otimes \cdots \otimes v_{i-1} \otimes 
Y_{V}(v_{i}, z_{i}-z_{i+1})v_{i+1}\nn
&&\quad\quad\quad\quad\quad\quad\quad\quad\quad 
\quad\quad\otimes \cdots \otimes v_{n+1}))
(z_{1}, \dots, z_{i-1}, z_{i+1}, \dots, z_{n+1})\rangle)\nn
&&\quad + (-1)^{n+1}R(\langle w', e^{z_{n+1}D_{W}}
Y_{W}^{R}((\Psi(v_{1}\otimes \cdots \otimes v_{n}))(z_{1}, \dots, z_{n}), -z_{n+1})v_{n+1}
\rangle).
\end{eqnarray}
For $w\in  \hat{C}^{0}(V, W)$, we define 
$\hat{\delta}^{0} w \in \hat{C}^{1}(V, W)$ by
$$((\hat{\delta}^{0}(w))(v))(z)=Y_{W}^{L}(v, z)w-e^{zD_{W}}Y_{W}^{R}(w, -z)v$$
for $v\in V$. 

It was proved in \cite{Hcoh}, \cite{Q2} and \cite{Q3} that 
$\hat{\delta}_{\infty}^{n}\circ \hat{\delta}_{\infty}^{n-1}=0$ for $n\in \Z_{+}$. Thus we have
the $n$-th cohomology $\hat{H}^{n}_{\infty}(V, W)=\ker \hat{\delta}_{\infty}^{n}/
\hat{\delta}_{\infty}^{n-1}(\hat{C}^{n-1}(V, W))$ 
for $n\in \Z_{+}$. 

We now discuss the relation between $\hat{H}^{1}_{\infty}(V, W)$ and derivations from $V$ to $W$.

\begin{defn}
{\rm Let $W$ be a $V$-bimodule. A {\it derivation} from $V$ to $W$ is a grading-preserving 
linear map $f: V\to W$ satisfying
$$f(Y_{V}(u, z)v)=Y_{W}^{R}(f(u), z)v+Y_{W}^{L}(u, z)f(v)$$
for $u, v\in V$. }
\end{defn}

Let $w\in \hat{C}^{0}_{\infty}(V, W)$. By definition, for $v\in V$, $Y_{W}^{L}(v, z)w$ and 
$e^{zD_{W}}Y_{W}^{R}(w, -z)v$ are power series in $z$. In particular,  
$\lim_{z\to 0}(Y_{W}^{L}(v, z)w-e^{zD_{W}}Y_{W}^{R}(w, -z)v)$ exists. 
Let $f_{w}: V\to W$ be defined by these limits. that is,
$$f_{w}(v)=\lim_{z\to 0}(Y_{W}^{L}(v, z)w-e^{zD_{W}}Y_{W}^{R}(w, -z)v)$$
for $v\in V$. By the $\mathbf{d}$-conjugation property for $Y_{W}^{L}$ and $Y_{W}^{R}$ and the 
fact that the weight of $w$ is $0$, 
\begin{eqnarray*}
a^{\mathbf{d}_{W}}f_{w}(v)&=&a^{\mathbf{d}_{W}}\lim_{z\to 0}
(Y_{W}^{L}(v, z)w-e^{zD_{W}}Y_{W}^{R}(w, -z)v)\nn
&=&\lim_{z\to 0}(Y_{W}^{L}(a^{\mathbf{d}_{V}}v, az)a^{\mathbf{d}_{W}}w
-e^{azD_{W}}Y_{W}^{R}(a^{\mathbf{d}_{W}}w, -az)a^{\mathbf{d}_{V}}v)\nn
&=&\lim_{z'\to 0}Y_{W}^{L}(a^{\mathbf{d}_{V}}v, z')w-e^{z'D_{W}}Y_{W}^{R}(w, -z')a^{\mathbf{d}_{V}}v\nn
&=&f_{w}(a^{\mathbf{d}_{V}}v)
\end{eqnarray*}
for $a\in \C^{\times}$. So $f_{w}$ preserves weights. 

From the $D$-derivative property of $Y_{W}^{L}$ and $Y_{W}^{R}$, $D_{W}w=0$
and Taylor's theorem, we obtain
\begin{eqnarray*}
\lefteqn{e^{z'D_{W}}(Y_{W}^{L}(v, z)w-e^{zD_{W}}Y_{W}^{R}(w, -z)v)}\nn
&&=e^{z'D_{W}}Y_{W}^{L}(v, z)w-e^{(z+z')D_{W}}Y_{W}^{R}(w, -z)v)\nn
&&=Y_{W}^{L}(e^{z'D_{W}}v, z)w-e^{(z+z')D_{W}}Y_{W}^{R}(e^{-z'D_{W}}w, -z)v)\nn
&&=Y_{W}^{L}(v, z+z')w-e^{(z+z')D_{W}}Y_{W}^{R}(w, -(z+z'))v
\end{eqnarray*}
for $z'\in \C$.
 Taking the limit $z\to 0$ on both sides and then replacing
$z'$ by $z$, we obtain
$$Y_{W}^{L}(v, z)w-e^{zD_{W}}Y_{W}^{R}(w, -z)v
=e^{zD_{W}}f_{w}(v).$$
Let $\Psi_{w}(v, z)=Y_{W}^{L}(v, z)w-e^{zD_{W}}Y_{W}^{R}(w, -z)v$ for $v\in V$.
Then $\Psi_{w}(v, z)=e^{zD_{W}}f_{w}(v)$ and
$f_{w}(v)=\lim_{z\to 0}\Psi_{w}(v, z)$. 

\begin{prop}
For $w\in \hat{C}^{0}_{\infty}(V, W)$, $f_{w}$ is a derivation from 
$V$ to $W$.
\end{prop}
\pf
We already know that $f_{w}$ preserves weights. 

For $u, v\in V$ and $w'\in W'$, 
\begin{align*}
&R(\langle w', \Psi_{w}(Y_{V}(u, z_{1}-z_{2})v, z_{2})\rangle)\nn
&\;=R(\langle w', Y_{W}^{L}(Y_{V}(u, z_{1}-z_{2})v, z_{2})w\rangle)
-R(\langle w', e^{z_{2}D_{W}}Y_{W}^{R}(w, -z_{2})Y_{V}(u, z_{1}-z_{2})v\rangle)\quad\quad\quad\quad\quad\nn
&\;=R(\langle w', Y_{W}^{L}(u, z_{1})Y_{W}^{L}(v, z_{2})w\rangle)
-R(\langle w', e^{z_{2}D_{W}}Y_{W}^{R}(Y_{W}^{R}(w, -z_{1})u, z_{1}-z_{2})v\rangle)\nn
&\;=R(\langle w', Y_{W}^{L}(u, z_{1})Y_{W}^{L}(v, z_{2})w\rangle)
-R(\langle w', Y_{W}^{L}(u, z_{1})e^{z_{2}D_{W}}Y_{W}^{R}(w, -z_{2})v\rangle)\nn
&\quad +R(\langle w', Y_{W}^{L}(u, z_{1})e^{z_{2}D_{W}}Y_{W}^{R}(w, -z_{2})v\rangle)
-R(\langle w', e^{z_{2}D_{W}}Y_{W}^{R}(Y_{W}^{R}(w, -z_{1})u, z_{1}-z_{2})v\rangle)
\end{align*}
\begin{align*}
&=R(\langle w', Y_{W}^{L}(u, z_{1})Y_{W}^{L}(v, z_{2})w\rangle)
-R(\langle w', Y_{W}^{L}(u, z_{1})e^{z_{2}D_{W}}Y_{W}^{R}(w, -z_{2})v\rangle)\nn
&\quad 
+R(\langle w', e^{z_{2}D_{W}}Y_{W}^{L}(u, z_{1}-z_{2})Y_{W}^{R}(w, -z_{2})v\rangle)
-R(\langle w', e^{z_{2}D_{W}}Y_{W}^{R}(Y_{W}^{R}(w, -z_{1})u, z_{1}-z_{2})v\rangle)\nn
&=R(\langle w', Y_{W}^{L}(u, z_{1})Y_{W}^{L}(v, z_{2})w\rangle)
-R(\langle w', Y_{W}^{L}(u, z_{1})e^{z_{2}D_{W}}Y_{W}^{R}(w, -z_{2})v\rangle)\nn
&\quad 
+R(\langle w', e^{z_{2}D_{W}}Y_{W}^{R}(Y_{W}^{L}(u, z_{1})w, -z_{2})v\rangle)
-R(\langle w', e^{z_{2}D_{W}}Y_{W}^{R}(e^{z_{1}D_{W}}Y_{W}^{R}(w, -z_{1})u, -z_{2})v\rangle)\nn
&=R(\langle w', (Y_{W}^{L}(u, z_{1})\Psi_{w}(v, z_{2})+e^{z_{2}D_{W}}
Y_{W}^{R}(\Psi_{w}(u, z_{1}), -z_{2})v)\rangle)\nn
&=R(\langle w', (Y_{W}^{L}(u, z_{1})\Psi_{w}(v, z_{2})+e^{z_{2}D_{W}}
Y_{W}^{R}(e^{z_{1}D_{W}}f_{w}(u), -z_{2})v)\rangle)\nn
&=R(\langle w', (Y_{W}^{L}(u, z_{1})\Psi_{w}(v, z_{2})+e^{z_{2}D_{W}}
Y_{W}^{R}(f_{w}(u), z_{1}-z_{2})v)\rangle)
\end{align*}
Taking the limit $z_{2}\to 0$ on both sides, we obtain
$$\langle w',  f_{w}(Y_{V}(u, z_{1})v)\rangle=
\langle w', (Y_{W}^{L}(u, z_{1})f_{w}(v)+Y_{W}^{R}(f_{w}(u), z_{1})v)\rangle,$$
proving that $f_{w}$ is a derivation from $V$ to $W$. 
\epfv

The derivation of the form $f_{w}$ is called an {\it inner derivation}. 
By definition, $\hat{\delta}^{0} \hat{C}^{0}_{\infty}(V, W)$ is exactly the space of 
$\Psi_{w}(\cdot, z)$ for $w\in \hat{C}^{0}_{\infty}(V, W)$. We have a linear isomorphism 
from $\hat{\delta}^{0}_{\infty}(\hat{C}^{0}_{\infty}(V, W))$ to the space of inner derivations
given by $\Psi_{w}(\cdot, z)\mapsto f_{w}$. 

Our next result and some other results in this paper
need the following lemma on iterated series convergent to suitable
rational functions:

\begin{lemma}\label{it-series-conv}
Let $f(\zeta_{1}, \dots, \zeta_{n})$ be a rational function in $\zeta_{1}, \dots, \zeta_{n}$
whose denominator is a product of homogeneous polynomials of degree one. 
Let 
\begin{align}\label{iterated-series}
\sum_{k_{n}\in \Z}\left(\cdots \left(\sum_{k_{1}\in \Z}
a_{k_{1}\cdots k_{n}}\zeta_{1}^{k_{1}}\right)\cdots \zeta_{n}^{k_{n}} \right)
\end{align}
be an iterated series satisfying the following conditions:
\begin{enumerate}

\item The iterated series  (\ref{iterated-series}) is convergent 
in a region to $f(\zeta_{1}, \dots, \zeta_{n})$.

\item The truncation condition: 
The coefficients of (\ref{iterated-series}) in fixed powers of $\zeta_{i+1}, \dots, \zeta_{n}$
have only finitely many negative power terms in $\zeta_{i}$.

\item A Laurent expansion of  $f(\zeta_{1}, \dots, \zeta_{n})$ by expanding 
powers of homogeneous polynomials of degree one using the binomial expansion 
satisfies the same truncation property and its region of convergence contains 
the region of convergence of  (\ref{iterated-series}).

\end{enumerate}
Then the multisum of (\ref{iterated-series})  is equal 
 to the Laurent series of $f(\zeta_{1}, \dots, \zeta_{n})$ in Condition 3. In particular, 
the multisum of (\ref{iterated-series}) is convergent absolutely  
to $f(\zeta_{1}, \dots, \zeta_{n})$ in the region of convergence of 
the Laurent series of $f(\zeta_{1}, \dots, \zeta_{n})$  in Condition 3. 
\end{lemma}
\pf
Taking $i=n$ in Condition 2, we know that  (\ref{iterated-series})
has only finitely many negative powers of $\zeta_{n}$. 
Since  (\ref{iterated-series}) is convergent absolutely in a region to 
$f(\zeta_{1}, \dots, \zeta_{n})$, there must be a region containing 
a disk with $0$ deleted in the variable $\zeta_{n}$ such that 
(\ref{iterated-series}) is convergent absolutely in this possibly larger region to 
$f(\zeta_{1}, \dots, \zeta_{n})$. Then the coefficients of  (\ref{iterated-series})
as a series in $\zeta_{n}$ must be convergent absolutely 
to the coefficients of the Laurent series in $\zeta_{n}$ of $f(\zeta_{1}, \dots, \zeta_{n})$
expanded in the same region for $\zeta_{n}$. 
Since the Laurent series of $f(\zeta_{1}, \dots, \zeta_{n})$
in Condition 3 satisfies the same truncation condition as the one for 
(\ref{iterated-series}), we see that the coefficients of  (\ref{iterated-series})
and the corresponding coefficients of the Laurent series of $f(\zeta_{1}, \dots, \zeta_{n})$
as series in $\zeta_{n}$ must be convergent in a region to a common 
rational function in $\zeta_{1}, 
\dots, \zeta_{n-1}$ whose denominator is a product of homogeneous polynomials of degree one.  
For these coefficient series in 
$\zeta_{1}, \dots, \zeta_{n-1}$, 
using the above argument repeatedly, we see that the coefficients of 
the series are all the same as the corresponding coefficients 
of the Laurent series of $f(\zeta_{1}, \dots, \zeta_{n})$ in Condition 3. Thus the lemma 
is proved. 
\epf

\begin{rema}
{\rm Lemma \ref{it-series-conv} can also be derived from Lemma 4.5 or Lemma 4.7 in \cite{Q1}.
The rational function in Lemma \ref{it-series-conv} is more special but is exactly 
what we need in this paper. Also the conditions in Lemma \ref{it-series-conv} are what 
we can see easily in our proofs of the results in this paper.}
\end{rema}

\begin{thm}\label{1-coh-der}
Let $V$ be a meromorphic open-string vertex algebra and $W$ a $V$-bimodule. 
For
$\Psi\in \hat{C}^{1}_{\infty}(V, W)$, let $f_{\Psi}: V\to W$ be defined by 
$f_{\Psi}(v)=(\Psi(v))(0)$ for $v\in V$. Then if $\Psi$ is closed,
$f_{\Psi}$ is a derivation from $V$ to $W$
and the map given by 
$\Psi\mapsto f_{\Psi}$ is a linear isomorphism from the space of closed  $1$-cochains
to the space of derivations from $V$ to $W$. 
In particular, 
the first cohomology $\hat{H}^{1}_{\infty}(V, W)$ of $V$ with coefficients in $W$ is isomorphic to 
the quotient of the space of derivations from $V$ to $W$ by the space 
of inner derivations. 
\end{thm}
\pf 
Given an element 
$\Psi\in \hat{C}^{1}_{\infty}(V, W)$, by definition,
\begin{eqnarray*}
\lefteqn{\langle w', ((\hat{\delta}^{1}\Psi)(v_{1}\otimes v_{2}))(z_{1}, z_{2})\rangle}\nn
&&=R(\langle w', Y_{W}^{L}(v_{1}, z_{1})(\Psi(v_{2}))(z_{2})\rangle)
-R(\langle w', (\Psi(Y_{V}(v_{1}, z_{1}-z_{2})v_{2}))(z_{2})\rangle)\nn
&&\quad 
+R(\langle w', e^{z_{2}D_{W}}Y_{W}^{R}((\Psi(v_{1}))(z_{1}), -z_{2})v_{2}\rangle).
\end{eqnarray*}
Then the $\mathbf{d}$-conjugation property satisfied by $\Psi$ 
implies that $f_{\Psi}$ preserves weights. The $D$-derivative property satisfied by
$\Psi$ gives $(\Psi(v))(z)=e^{zD}f_{\Psi}(v)$ for $v\in V$.
If $\Psi$ is closed, then we obtain
\begin{eqnarray}\label{vertex-derivation}
\lefteqn{R(\langle w', (\Psi(Y_{V}(v_{1}, z_{1}-z_{2})v_{2}))(z_{2})\rangle)}\nn
&&=R(\langle w',  Y_{W}^{L}(v_{1}, z_{1})(\Psi(v_{2}))(z_{2})\rangle)
+R(\langle w', e^{z_{2}D_{W}}Y_{W}^{R}((\Psi(v_{1}))(z_{1}), -z_{2})v_{2}\rangle)\nn
&&=R(\langle w', Y_{W}^{L}(v_{1}, z_{1})(\Psi(v_{2})(z_{2})\rangle)
+R(\langle w', e^{z_{2}D_{W}}Y_{W}^{R}(e^{z_{1}D_{W}}f_{\Psi}(v_{1}), -z_{2})v_{2}\rangle)\nn
&&= R(\langle w', Y_{W}^{L}(v_{1}, z_{1})(\Psi(v_{2})(z_{2})\rangle)
+R(\langle w', e^{z_{2}D_{W}}Y_{W}^{R}(f_{\Psi}(v_{1}), z_{1}-z_{2})v_{2}\rangle).
\end{eqnarray}
Letting $z_{2}= 0$ on both sides of (\ref{vertex-derivation}),
we obtain 
\begin{equation}\label{derivation}
\langle w', f_{\Psi}(Y_{V}(v_{1}, z_{1})v_{2})\rangle =
 \langle w', Y_{W}^{L}(v_{1}, z_{1})f_{\Psi}(v_{2})\rangle
+\langle w', Y_{W}^{R}(f_{\Psi}(v_{1}), z_{1})v_{2}\rangle.
\end{equation}
Since (\ref{derivation}) holds for all $w'\in W'$, we have proved that 
$f_{\Psi}$ is a derivation from $V$ to $W$. We obtain a linear map defined by
$\Psi\mapsto f_{\Psi}$ from 
the space of closed $1$-cochains to the space of derivations from $V$ to $W$.

It is clear that given any derivation $f$ from $V$ to $W$, 
for $v\in V$,  $e^{zD_{W}}f(v)$ is a $\overline{W}$-valued rational function 
in $z$. 
To see that  the linear map defined by
$\Psi\mapsto f_{\Psi}$ is invertible, we first prove that 
$\Psi_{f}\in \hom(V, \widetilde{W}_{z})$
defined by $(\Psi_{f}(v))(z)=e^{zD_{W}}f(v)$ for $v\in V$ is a closed $1$-cochain. 
We have $\frac{d}{dz}e^{zD_{w}}f(v)
=D_{W}e^{zD_{W}}f(v)$. On the other hand, the same proof as the proof of Lemma 2.1 in 
\cite{H1st-sec-coh} shows that $f(\one)=0$. Thus we obtain 
\begin{align*}
f(D_{V}v)&=\lim_{z\to 0}\frac{d}{dz}f(e^{zD_{V}}v)\nn
&=\lim_{z\to 0}\frac{d}{dz}f(Y_{V}(v, z)\one)\nn
&=\lim_{z\to 0}\frac{d}{dz}Y_{W}^{R}(f(v), z)\one
+\lim_{z\to 0}\frac{d}{dz}Y_{W}^{L}(v, z)f(\one)\nn
&=\lim_{z\to 0}Y_{W}^{R}(D_{W}f(v), z)\one\nn
&=D_{W}f(v)
\end{align*}
for $v\in V$. So $\Psi_{f}$ satisfies the $D$-derivative property. 
Since $D_{W}$ is an operator of weight $1$
and $f$ preserves weights, we have $a^{\mathbf{d}_{W}}e^{zD_{W}}f(v)
=e^{za^{\mathbf{d}_{W}}D_{W}a^{-\mathbf{d}_{W}}}a^{\mathbf{d}_{W}}f(v)
=e^{azD_{w}}f(a^{\mathbf{d}_{W}}v)$ for $a\in \C^{\times}$, 
proving the $\mathbf{d}$-conjugation property for $\Psi_{f}$. 

For $k, l, m\in \N$, $v_{1}, \dots, v_{k+l+m}\in V$, $w'\in W'$, 
\begin{align}\label{1-coh-der-1}
&\langle w', Y^{L}_{W}(v_{1}, z_{1})\cdots Y^{L}_{W}(v_{k}, z_{k})
Y^{R}_{W}(e^{(\xi-\zeta)D_{W}}f(Y_{V}(v_{k+1}, z_{k+1}-\xi)\cdots Y_{V}(v_{k+l}, z_{k+l}-\xi)\one),
\zeta)\cdot\nn
&\quad\quad\quad\quad\cdot Y_{V}(v_{k+l+1}, z_{k+l+1})
\cdots Y_{V}(v_{k+l+m}, z_{k+l+m})\one\rangle\nn
&\quad=\sum_{i=1}^{l}\langle w', Y^{L}_{W}(v_{1}, z_{1})\cdots Y^{L}_{W}(v_{k}, z_{k})\cdot\nn
&\quad\quad\quad\quad\cdot
Y^{R}_{W}(e^{(\xi-\zeta)D_{W}}Y_{W}^{L}(v_{k+1}, z_{k+1}-\xi)\cdots
Y_{W}^{L}(v_{k+i-1}, z_{k+i-1}-\xi)\cdot\nn
&\quad\quad\quad\quad\cdot Y_{W}^{R}(f(v_{k+i}), z_{k+i}-\xi) Y_{V}(v_{k+i+1}, z_{k+i+l}-\xi)
\cdots Y_{V}(v_{k+l}, z_{k+l}-\xi)\one,
\zeta)\cdot\nn
&\quad\quad\quad\quad\cdot Y_{V}(v_{k+l+1}, z_{k+l+1})
\cdots Y_{V}(v_{k+l+m}, z_{k+l+m})\one\rangle\nn
&\quad=\sum_{i=1}^{l}\langle w', Y^{L}_{W}(v_{1}, z_{1})\cdots Y^{L}_{W}(v_{k}, z_{k})
Y^{R}_{W}(Y_{W}^{L}(v_{k+1}, (z_{k+1}-\xi)+(\xi-\zeta))\cdot\nn
&\quad\quad\quad\quad\cdots
Y_{W}^{L}(v_{k+i-1}, (z_{k+i-1}-\xi)+(\xi-\zeta))
Y_{W}^{R}(f(v_{k+i}), (z_{k+i}-\xi)+(\xi-\zeta)) \cdot\nn
&\quad\quad\quad\quad\cdot 
Y_{V}(v_{k+i+1}, (z_{k+i+l}-\xi)+(\xi-\zeta))
\cdots Y_{V}(v_{k+l}, (z_{k+l}-\xi)+(\xi-\zeta))\one,
\zeta)\cdot\nn
&\quad\quad\quad\quad\cdot Y_{V}(v_{k+l+1}, z_{k+l+1})
\cdots Y_{V}(v_{k+l+m}, z_{k+l+m})\one\rangle,
\end{align}
where in the right-hand side of (\ref{1-coh-der-1}), the negative powers of 
$(z_{k+j}-\xi)+(\xi-\zeta)=z_{k+j}-\zeta$ for $j=1, \dots, l$ are expanded
as power series in $\xi-\zeta$. 
Using (\ref{1-coh-der-1}) and the properties of the $V$-bimodule $W$, we now prove
that $\Psi_{f}$ can be composed with an arbitrary number of vertex operators.

We first prove that every term in the right-hand side of (\ref{1-coh-der-1})
is  convergent absolutely in the region given by
$|z_{1}|>\cdots >|z_{k}|$,  $|z_{k+1}-\xi|>\cdots >|z_{k+l}-\xi|$,
$|z_{k+l+1}|>\cdots > |z_{k+l+m}|$,
$|z_{a}|>|z_{k+j}-\xi|+|\xi-\zeta|+|\zeta|$
and $|\zeta|>|z_{k+j}-\xi|+|z_{k+l+p}|+|\xi-\zeta|$ for  $a=1, \dots, k$,
 $j=1, \dots, l$ and $p=1, \dots, m$. 
We use induction on $l$ to prove that the $i$-th term in the right-hand side of (\ref{1-coh-der-1})
is convergent absolutely. When $l=1$, 
the only term in the right-hand side of (\ref{1-coh-der-1})
is 
\begin{align}\label{1-coh-der-1.2}
\langle w', Y^{L}_{W}(v_{1}, z_{1})&\cdots Y^{L}_{W}(v_{k}, z_{k})
Y^{R}_{W}(Y_{W}^{R}(f(v_{k+1}), (z_{k+1}-\xi)+(\xi-\zeta)) \one,
\zeta)\cdot\nn
&\cdot Y_{V}(v_{k+2}, z_{k+2})
\cdots Y_{V}(v_{k+1+m}, z_{k+1+m})\one\rangle.
\end{align}
By the associativity of $Y_{W}^{R}$, we know that 
$Y^{R}_{W}(Y_{W}^{R}(f(v_{k+1}), (z_{k+1}-\xi)+(\xi-\zeta)) \one,
\zeta)v$  for $v\in V$ is convergent absolutely in the region 
$|\zeta|>|z_{k+1}-\xi|+|\xi-\zeta|$ and is equal 
to $Y^{R}_{W}(f(v_{k+1}), z_{k+1})v$. So in this region, the sum 
involving $Y^{R}_{W}(Y_{W}^{R}(f(v_{k+1}), (z_{k+1}-\xi)+(\xi-\zeta)) \one,
\zeta)v$ in (\ref{1-coh-der-1.2}) is  convergent absolutely
to the series 
\begin{align}\label{1-coh-der-1.4}
\langle w', Y^{L}_{W}(v_{1}, z_{1})\cdots Y^{L}_{W}(v_{k}, z_{k})
Y^{R}_{W}(f(v_{k+1}), z_{k+1})
Y_{V}(v_{k+2}, z_{k+2})
\cdots Y_{V}(v_{k+1+m}, z_{k+1+m})\one\rangle.
\end{align}
By the rationality for products of 
left and right vertex operators 
for the $V$-bimodule $W$, (\ref{1-coh-der-1.4}) is convergent absolutely 
in the region $|z_{1}|>\cdots> |z_{k+1+m}|$ to a rational function 
in $z_{1}, \dots, z_{k+1+m}$ with the only possible poles at 
$z_{a}=z_{b}$ for $a\ne b$.
Expand this rational function as a Laurent series in the variables 
$z_{1}, \dots, z_{k}$, $z_{k+1}-\xi$, $z_{k+2}, \dots, z_{k+1+m}$,
$\zeta$ and $\xi-\zeta$ in the region $|z_{1}|>\cdots>|z_{k}|>0$,
$|z_{k+2}|>|z_{k+1+m}|$, $|z_{i}|>
|z_{k+1}-\xi|+|\xi-\zeta|+|\zeta|$ for $i=1, \dots, k$, 
$|\zeta|>|z_{k+1}-\xi|+|z_{k+1+p}|+|\xi-\zeta|$ for $p=1, \dots, m$.
This expansion has the same form as the 
series  (\ref{1-coh-der-1.2}) and  one of the sums in the series 
(\ref{1-coh-der-1.2}) is convergent absolutely in the region $|\zeta|>|z_{k+1}-\xi|+|\xi-\zeta|$
to the series (\ref{1-coh-der-1.4}). By Lemma \ref{it-series-conv}, 
we see that this expansion must be equal to (\ref{1-coh-der-1.2}).
In particular, we see that (\ref{1-coh-der-1.2}) must be  
convergent absolutely to this rational function in the region $|z_{1}|>\cdots>|z_{k}|>0$
$|z_{k+2}|>|z_{k+1+m}|$, $|z_{a}|>
|z_{k+1}-\xi|+|\xi-\zeta|+|\zeta|$ for $a=1, \dots, k$, 
$|\zeta|>|z_{k+1}-\xi|+|z_{k+1+p}|+|\xi-\zeta|$ for $p=1, \dots, m$.
This proves that in the case $l=1$,  the right-hand side of (\ref{1-coh-der-1})
is  convergent absolutely in the region given by
$|z_{1}|>\cdots >|z_{k}|$, 
$|z_{k+2}|>\cdots > |z_{k+1+m}|$,
$|z_{a}|>|z_{k+1}-\xi|+|\xi-\zeta|+|\zeta|$
and $|\zeta|>|z_{k+1}-\xi|+|z_{k+1+p}|+|\xi-\zeta|$ for  $a=1, \dots, k$
and $p=1, \dots, m$. 

Now we assume that for $l=q-1$, every term in the right-hand side of (\ref{1-coh-der-1})
is convergent absolutely  in the region given by
$|z_{1}|>\cdots >|z_{k}|$,  $|z_{k+1}-\xi|>\cdots >|z_{k+l}-\xi|$,
$|z_{k+l+1}|>\cdots > |z_{k+l+m}|$,
$|z_{a}|>|z_{k+j}-\xi|+|\xi-\zeta|+|\zeta|$
and $|\zeta|>|z_{k+j}-\xi|+|z_{k+l+p}|+|\xi-\zeta|$ for  $a=1, \dots, k$,
 $j=1, \dots, l$ and $p=1, \dots, m$.  In the case $l=q$ and $i=1$,
$Y^{R}_{W}(Y_{W}^{R}(f(v_{k+1}), (z_{k+1}-\xi)+(\xi-\zeta)) u,
\zeta)v$  for $u, v\in V$ is convergent absolutely in the region 
$|\zeta|>|z_{k+1}-\xi|+|\xi-\zeta|$ and is equal 
to $Y^{R}_{W}(f(v_{k+1}), z_{k+1})Y_{V}(u, \zeta)v$. In the case $l=q$
 and $i>1$, 
$Y^{R}_{W}(Y_{W}^{L}(v_{k+1}, (z_{k+1}-\xi)+(\xi-\zeta)w,
\zeta)v$ for $v\in V$ and $w\in W$ is  convergent absolutely in the region 
$|\zeta|>|z_{k+1}-\xi|+|\xi-\zeta|$ and is equal 
to $Y^{L}_{W}(v_{k+1}, z_{k+1})Y^{R}_{W}(w, \zeta)v$.
So in the region $|\zeta|>|z_{k+1}-\xi|+|\xi-\zeta|$, the sum 
involving $Y^{R}_{W}(Y_{W}^{R}(f(v_{k+1}), (z_{k+1}-\xi)+(\xi-\zeta)) u,
\zeta)v$ in the case $i=1$ or the sum inolving 
$Y^{R}_{W}(Y_{W}^{L}(v_{k+1}, (z_{k+1}-\xi)+(\xi-\zeta)w,
\zeta)v$ in the case $i>1$ in the $i$-th term in 
(\ref{1-coh-der-1.2}) is  convergent absolutely
to the series 
\begin{align}\label{1-coh-der-1.5}
&\langle w', Y^{L}_{W}(v_{1}, z_{1})\cdots Y^{L}_{W}(v_{k}, z_{k})
Y^{R}_{W}(f(v_{k+1}), z_{k+1})\cdot\nn
&\quad\quad\cdot
Y_{V}(Y_{V}(v_{k+2}, (z_{k+2}-\xi)+(\xi-\zeta))
\cdots Y_{V}(v_{k+l}, (z_{k+l}-\xi)+(\xi-\zeta))\one,
\zeta)\cdot\nn
&\quad\quad\cdot Y_{V}(v_{k+l+1}, z_{k+l+1})
\cdots Y_{V}(v_{k+l+m}, z_{k+l+m})\one\rangle
\end{align}
or
\begin{align}\label{1-coh-der-1.6}
&\langle w', Y^{L}_{W}(v_{1}, z_{1})\cdots Y^{L}_{W}(v_{k}, z_{k})Y^{L}_{W}(v_{k+1}, z_{k+1})
Y^{R}_{W}(Y_{W}^{L}(v_{k+2}, (z_{k+2}-\xi)+(\xi-\zeta))\cdot\nn
&\quad\quad\cdots
Y_{W}^{L}(v_{k+i-1}, (z_{k+i-1}-\xi)+(\xi-\zeta))
Y_{W}^{R}(f(v_{k+i}), (z_{k+i}-\xi)+(\xi-\zeta)) \cdot\nn
&\quad\quad\cdot 
Y_{V}(v_{k+i+1}, (z_{k+i+l}-\xi)+(\xi-\zeta))
\cdots Y_{V}(v_{k+l}, (z_{k+l}-\xi)+(\xi-\zeta))\one,
\zeta)\cdot\nn
&\quad\quad\cdot Y_{V}(v_{k+l+1}, z_{k+l+1})
\cdots Y_{V}(v_{k+l+m}, z_{k+l+m})\one\rangle,
\end{align}
respectively. 

We first discuss the case $i=1$. By the induction assumption with $W=V$, 
for $v'\in V'$, 
\begin{align*}
\langle v', Y_{V}(Y_{V}&(v_{k+2}, (z_{k+2}-\xi)+(\xi-\zeta))
\cdots Y_{V}(v_{k+l}, (z_{k+l}-\xi)+(\xi-\zeta))\one,
\zeta)\cdot\nn
&\quad\quad\cdot Y_{V}(v_{k+l+1}, z_{k+l+1})
\cdots Y_{V}(v_{k+l+m}, z_{k+l+m})\one\rangle
\end{align*}
is  convergent absolutely in the region $|z_{k+1}-\xi|>\cdots >|z_{k+l}-\xi|$,
$|z_{k+l+1}|>\cdots > |z_{k+l+m}|$
and $|\zeta|>|z_{k+j}-\xi|+|z_{k+l+p}|+|\xi-\zeta|$ for 
 $j=1, \dots, l$ and $p=1, \dots, m$ to the rational function 
 $$R(\langle v',
Y_{V}(v_{k+2}, z_{k+2})
\cdots Y_{V}(v_{k+1+m}, z_{k+1+m})\one\rangle),$$
or equivalently,
\begin{align*}
Y_{V}(Y_{V}&(v_{k+2}, (z_{k+2}-\xi)+(\xi-\zeta))
\cdots Y_{V}(v_{k+l}, (z_{k+l}-\xi)+(\xi-\zeta))\one,
\zeta)\cdot\nn
&\quad\quad\cdot Y_{V}(v_{k+l+1}, z_{k+l+1})
\cdots Y_{V}(v_{k+l+m}, z_{k+l+m})\one
\end{align*}
is convergent absolutely in this region to 
$$E(Y_{V}(v_{k+2}, z_{k+2})
\cdots Y_{V}(v_{k+l+m}, z_{k+l+m})\one).$$
But we know that in the region $|z_{1}|>\cdots> |z_{k+l+m}|$,
\begin{align*}
&\langle w', Y^{L}_{W}(v_{1}, z_{1})\cdots Y^{L}_{W}(v_{k}, z_{k})
Y^{R}_{W}(f(v_{k+1}), z_{k+1})
E(Y_{V}(v_{k+2}, z_{k+2})
\cdots Y_{V}(v_{k+l+m}, z_{k+l+m})\one)\rangle\nn
&\quad =\langle w', Y^{L}_{W}(v_{1}, z_{1})\cdots Y^{L}_{W}(v_{k}, z_{k})
Y^{R}_{W}(f(v_{k+1}), z_{k+1})
Y_{V}(v_{k+2}, z_{k+2})
\cdots Y_{V}(v_{k+l+m}, z_{k+l+m})\one\rangle
\end{align*}
is  convergent absolutely to  a rational function 
in $z_{1}, \dots, z_{k+1+m}$ with the only possible poles at 
$z_{a}=z_{b}$ for $a\ne b$. Thus in the region $|z_{1}|>\cdots> |z_{k+l+m}|$, 
$|z_{k+2}-\xi|>\cdots >|z_{k+l}-\xi|$, $|z_{a}|>|z_{k+j}-\xi|+|\xi-\zeta|+|\zeta|$
and $|\zeta|>|z_{k+j}-\xi|+|z_{k+l+p}|+|\xi-\zeta|$ for $a=1, \dots, k+1$,
 $j=2, \dots, l$ and $p=1, \dots, m$, 
\begin{align*}
&\langle w', Y^{L}_{W}(v_{1}, z_{1})\cdots Y^{L}_{W}(v_{k}, z_{k})
Y^{R}_{W}(f(v_{k+1}), z_{k+1})\cdot\nn
&\quad\quad\quad\quad\cdot
Y_{V}(Y_{V}(v_{k+2}, (z_{k+2}-\xi)+(\xi-\zeta))
\cdots Y_{V}(v_{k+l}, (z_{k+l}-\xi)+(\xi-\zeta))\one,
\zeta)\cdot\nn
&\quad\quad\quad\quad\cdot Y_{V}(v_{k+l+1}, z_{k+l+1})
\cdots Y_{V}(v_{k+l+m}, z_{k+l+m})\one\rangle\nn
&\quad =\langle w', Y^{L}_{W}(v_{1}, z_{1})\cdots Y^{L}_{W}(v_{k}, z_{k})
Y^{R}_{W}(f(v_{k+1}), z_{k+1})\cdot\nn
&\quad\quad\quad \quad\cdot 
E(Y_{V}(v_{k+2}, z_{k+2})
\cdots Y_{V}(v_{k+1+m}, z_{k+1+m})\one)\rangle
\end{align*}
is convergent absolutely to this rational function. 
Expand this rational function as a Laurent series in the variables 
$z_{1}, \dots, z_{k}$, $z_{k+1}-\xi, \dots, z_{k+l}-\xi$, $z_{k+l+1}, \dots, z_{k+l+m}$,
$\zeta$ and $\xi-\zeta$ in the region $|z_{1}|>\cdots >|z_{k}|$,
$|z_{k+1}-\xi|>\cdots >|z_{k+l}-\xi|$,
$|z_{k+l+1}|>\cdots > |z_{k+l+m}|$,
$|z_{a}|>|z_{k+j}-\xi|+|\xi-\zeta|+|\zeta|$
and $|\zeta|>|z_{k+j}-\xi|+|z_{k+l+p}|+|\xi-\zeta|$ for  $a=1, \dots, k$,
 $j=1, \dots, l$ and $p=1, \dots, m$. 
This expansion has the same form as the first term in the
series  (\ref{1-coh-der-1}) and  one of the sums in the series 
(\ref{1-coh-der-1}) is convergent absolutely in the region $|\zeta|>|z_{k+1}-\xi|+|\xi-\zeta|$
 to the series (\ref{1-coh-der-1.5}) which is in turn convergent absolutely  
to this rational function. By Lemma \ref{it-series-conv}, 
we see that this expansion must be equal to the first term in the
series  (\ref{1-coh-der-1}).
In particular, we see that the first term in the
series  (\ref{1-coh-der-1}) must be  
convergent absolutely to this rational 
function in 
$z_{1}, \dots, z_{k}$, $z_{k+1}-\xi, \dots, z_{k+l}-\xi$, $z_{k+l+1}, \dots, z_{k+l+m}$,
$\zeta$ and $\xi-\zeta$ in the region $|z_{1}|>\cdots >|z_{k}|$,
$|z_{k+1}-\xi|>\cdots >|z_{k+l}-\xi|$,
$|z_{k+l+1}|>\cdots > |z_{k+l+m}|$,
$|z_{a}|>|z_{k+j}-\xi|+|\xi-\zeta|+|\zeta|$
and $|\zeta|>|z_{k+j}-\xi|+|z_{k+l+p}|+|\xi-\zeta|$ for  $a=1, \dots, k$,
 $j=1, \dots, l$ and $p=1, \dots, m$. 
This proves that in the case $l=q$, the first term in the right-hand side of (\ref{1-coh-der-1})
is convergent absolutely to the rational function above in this region.

We now discuss the case $i>1$ which is more straightforward 
than the case $i=1$. Using the induction assumption, 
(\ref{1-coh-der-1.6}) is  convergent absolutely in the region $|z_{1}|>\cdots >|z_{k+1}|$,
$|z_{k+2}-\xi|>\cdots >|z_{k+l}-\xi|$,
$|z_{k+l+1}|>\cdots > |z_{k+l+m}|$,
$|z_{a}|>|z_{k+j}-\xi|+|\xi-\zeta|+|\zeta|$
and $|\zeta|>|z_{k+j}-\xi|+|z_{k+l+p}|+|\xi-\zeta|$ for  $a=1, \dots, k+1$,
 $j=2, \dots, l$ and $p=1, \dots, m$ to a rational function 
in $z_{1}, \dots, z_{k+1+m}$ with the only possible poles at 
$z_{a}=z_{b}$ for $a\ne b$. Expanding this rational function 
as a series in 
 $z_{1}, \dots, z_{k}$, $z_{k+1}-\xi, \dots, z_{k+l}-\xi$, $z_{k+l+1}, \dots, z_{k+l+m}$,
$\zeta$ and $\xi-\zeta$ in the region $|z_{1}|>\cdots >|z_{k}|$,
$|z_{k+1}-\xi|>\cdots >|z_{k+l}-\xi|$,
$|z_{k+l+1}|>\cdots > |z_{k+l+m}|$,
$|z_{a}|>|z_{k+j}-\xi|+|\xi-\zeta|+|\zeta|$
and $|\zeta|>|z_{k+j}-\xi|+|z_{k+l+p}|+|\xi-\zeta|$ for  $a=1, \dots, k$,
 $j=1, \dots, l$ and $p=1, \dots, m$. The remaining steps are the same as 
in the $i=1$ case except that we replace the first term in (\ref{1-coh-der-1}) 
and (\ref{1-coh-der-1.5}) by the $i$-th term and (\ref{1-coh-der-1.6}), respectively. 
This proves that in the case $l=q$, the $i$-th term in the right-hand side of (\ref{1-coh-der-1})
is convergent absolutely  to the rational function above in the region above.

Since the poles of the rational function  that 
 (\ref{1-coh-der-1}) converges to are the poles of 
\begin{align}\label{1-coh-der-2}
\sum_{i=1}^{l}R(\langle w', &Y^{L}_{W}(v_{1}, z_{1})\cdots 
Y_{W}^{L}(v_{k+i-1}, z_{k+i-1})\cdot\nn
&\quad\quad\cdot 
Y_{W}^{R}(f(v_{k+i}), z_{k+i}) 
Y_{V}(v_{k+i+1}, z_{k+i+l})
\cdots  Y_{V}(v_{k+l+m}, z_{k+l+m})\one\rangle),
\end{align}
we see from the proof of the convergence above
that the existence of $p_{ij}$ for $i, j=1, \dots, k+l+m$ 
follows from the pole-order condition satisfied by $W$. 
Thus $\Psi_{f}$ can be composed with an arbitrary number of vertex operators.

For $v_{1}, v_{2}\in V$
and $w'\in W'$, by the definition of $\hat{\delta}^{1}_{\infty}$ and the fact that $f$ is 
a derivation, 
\begin{align*}
\langle w',& ((\hat{\delta}^{1}_{\infty}(\Psi_{f}))(v_{1}\otimes  v_{2}))
(z_{1}, z_{2})\rangle\nn
&=R(\langle w', Y_{W}^{L}(v_{1}, z_{1})(e^{z_{2}D_{W}}f(v_{2}))\rangle)\nn
&\quad -R(\langle w', 
(e^{z_{2}D_{W}}f(Y_{V}(v_{1}, z_{1}-z_{2})v_{2})\rangle)\nn
&\quad + R(\langle w', e^{z_{2}D_{W}}
Y_{W}^{R}(e^{z_{1}D_{W}}f(v_{1}), -z_{2})v_{2}
\rangle)
\end{align*}
\begin{align*}
&\;=R(\langle w', e^{z_{2}D_{W}}Y_{W}^{L}(v_{1}, z_{1}-z_{2})f(v_{2})\rangle)\nn
&\;\quad -R(\langle w', 
e^{z_{2}D_{W}}f(Y_{V}(v_{1}, z_{1}-z_{2})v_{2}\rangle)\nn
&\;\quad + R(\langle w', e^{z_{2}D_{W}}
Y_{W}^{R}(f(v_{1}), z_{1}-z_{2})v_{2}
\rangle)\nn
&\;=0.
\end{align*}
So  $\Psi_{f}$ is  closed.

Now from the formula $(\Psi(v))(z)=e^{zD_{W}}f_{\Psi}(v)$, 
we see that the map $\Psi\mapsto f_{\Psi}$ above 
is invertible with the inverse given by $f\mapsto \Psi_{f}$. 
We obtain an isomorphism from the space of closed $1$-cochains
to the space of derivations from $V$ to $W$. Moreover, we have proved that  this isomorphism 
maps exact $1$-cochains to inner derivations. Thus $\hat{H}^{1}_{\infty}(V, W)$ is isomorphic to the 
quotient space of the space of derivations from $V$ to $W$ by the space of inner derivations. 
\epfv

We now discuss another type of derivations from $V$ to $W$ which play an important role in our main theorem.

\begin{prop}\label{wt-1-der}
Let $w$ be an element of $W_{[1]}$ such that 
\begin{equation}\label{wt-1-der-1}
e^{xD_{W}}Y_{W}^{L}(v, -x)-Y_{W}^{R}(w, x)v\in W[x, x^{-1}]
\end{equation}
for $v\in V$. Then $g_{w}: V\to W$ defined by 
$g_{w}(v)=\res_{x}Y_{W}^{R}(w, x)v=(Y_{W}^{R})_{0}(w)v$ for $v\in V$ is a derivation from $V$ to $W$.
\end{prop}
\pf
Since the weight of $w$ is $1$, the weight of the map $g_{w}=(Y_{W}^{R})_{0}(w): V\to W$ is 
$\wt w-0-1=0$. So $g_{w}$ preserves weights. 

For $u, v\in V$ and $w'\in W'$, by the definition of bimodule, we have the associativity
\begin{equation}\label{wt-1-der-2}
R(\langle w', Y_{W}^{R}(w, z_{1})Y_{V}(u, z_{2})v\rangle)
=R(\langle w', Y_{W}^{R}(Y_{W}^{R}(w, z_{1}-z_{2})u, z_{2})v\rangle).
\end{equation}
On the other hand, by (\ref{wt-1-der-1}) and the compatibility of the left and right actions, 
there exists a rational function  $h(z_{1}, z_{2})$
 in $z_{1}$ and $z_{2}$ with the only possible poles at $z_{2}=0$ and $z_{1}=z_{2}$ such that
\begin{align}\label{wt-1-der-3}
R&(\langle w', Y_{W}^{R}(Y_{W}^{R}(w, z_{1}-z_{2})u, z_{2})v\rangle)\nn
&=R(\langle w', Y_{W}^{R}(e^{(z_{1}-z_{2})D_{W}}Y_{W}^{L}(u, -(z_{1}-z_{2}))w, z_{2})v\rangle)+h(z_{1}, z_{2})\nn
&=R(\langle w', Y_{W}^{R}(Y_{W}^{L}(u, z_{2}-z_{1})w, z_{1})v\rangle)+h(z_{1}, z_{2})\nn
&=R(\langle w', Y_{W}^{L}(u, z_{2})Y_{W}^{R}(w, z_{1})v\rangle)+h(z_{1}, z_{2}).
\end{align}
Integrate the same rational function  in (\ref{wt-1-der-2}) and (\ref{wt-1-der-3}) with respect to $z_{1}$
along a circle $C_{\infty}$  centered at $0$ and containing $z_{2}$ in its interior.  Since $h(z_{1}, z_{2})$ is analytic 
in the disk $|z_{1}|<|z_{2}|$ as a function of $z_{1}$, its integral along a circle of radius less than $|z_{2}|$ is $0$. 
By this fact,  (\ref{wt-1-der-2}) and (\ref{wt-1-der-3}), Cauchy's integral theorem  and the convergence region of 
$\langle w', Y_{W}^{R}(w, z_{1})Y_{V}(u, z_{2})v\rangle$, 
$\langle w', Y_{W}^{R}(Y_{W}^{R}(w, z_{1}-z_{2})u, z_{2})v\rangle$ and 
$\langle w', Y_{W}^{L}(u, z_{2})Y_{W}^{R}(w, z_{1})v\rangle$, we obtain 
\begin{align}\label{wt-1-der-4}
\oint_{C_{\infty}}&\langle w', Y_{W}^{R}(w, z_{1})Y_{V}(u, z_{2})v\rangle dz_{1}\nn
&=\oint_{C_{z_{2}}}\langle w', Y_{W}^{R}(Y_{W}^{R}(w, z_{1}-z_{2})u, z_{2})v\rangle dz_{1}
+\oint_{C_{0}}(\langle w', Y_{W}^{L}(u, z_{2})Y_{W}^{R}(w, z_{1})v\rangle+h(z_{1}, z_{2})) dz_{1}\nn
&=\oint_{C_{z_{2}}}\langle w', Y_{W}^{R}(Y_{W}^{R}(w, z_{1}-z_{2})u, z_{2})v\rangle dz_{1}
+\oint_{C_{0}}\langle w', Y_{W}^{L}(u, z_{2})Y_{W}^{R}(w, z_{1})v\rangle dz_{1},
\end{align}
where $C_{z_{2}}$ and $C_{0}$ are circles centered at $z_{2}$ and $0$, respectively,
with radii less than $|z_{2}|$. By the definition of $g_{w}$, (\ref{wt-1-der-4}) gives 
$$g_{w}(Y_{V}(u, z_{2})v)=Y_{W}^{R}(g_{w}(u), z_{2})v+Y_{W}^{L}(u, z_{2})g_{w}(v),$$
proving that $g_{w}$ is a derivation from $V$ to $W$. 
\epfv

We shall call the derivation in Proposition \ref{wt-1-der} a {\it zero-mode derivation}
since it is in fact the zero-mode of the right vertex operator of the element $w$. We shall 
denote the subspace of $\hat{H}_{\infty}^{1}(V, W)$ consisting of the cosets 
containing zero-mode derivations by $\hat{Z}_{\infty}^{1}(V, W)$.

\renewcommand{\theequation}{\thesection.\arabic{equation}}
\renewcommand{\thethm}{\thesection.\arabic{thm}}
\setcounter{equation}{0}
\setcounter{thm}{0}
\section{A $V$-bimodule constructed from two left $V$-modules}

In this section, we construct a $V$-bimodule $H^{N}$ for $N\in \Z$ from two left $V$-modules. 
The bimodule $H^{N}$ is analogous to the bimodule $\hom(M_{1},
M_{2})$ for an associative algebra $A$ constructed from left $A$-modules 
$M_{1}$ and $M_{2}$. In the category of modules for a
vertex operator algebra, the first analogue of this bimodule was in fact given 
in the construction of the $Q(z)$-tensor product by Lepowsky and the first author
in \cite{HL0} (although the term $Q(z)$-tensor product was introduced later in 
\cite{HL1}). This is the reason why a $Q(z_{0})$-tensor product appears in 
Theorem \ref{F(v)-g-r} below. But even when $V$ is a vertex operator algebra, 
the analogue of the bimodule  $\hom(M_{1},
M_{2})$ needed in this paper is different. 
What we want is an analogue when $V$ is viewed as
a meromorphic open-string vertex algebra. This is the main reason why the construction 
in this section is difficult. 

In the rest of this paper, we fix a meromorphic open-string vertex algebra $V$. 
Let  $W_{1}$ and $W_{2}$ be two 
left $V$-modules. Recall that by our convention, $V$,  
$W_{1}$ and $W_{2}$,  in particular, satisfy the 
pole-order conditions. 
Let $\widehat{(W_{2})}_{z}$ be the space of $\overline{W_{2}}$-valued rational 
functions with the only possible pole at $z=0$. Recall that $\widetilde{(W_{2})}_{z}$ is the space of 
$\overline{W_{2}}$-valued holomorphic functions. Thus $\widehat{(W_{2})}_{z}\supset \widetilde{(W_{2})}_{z}$.

Let  $H$ be the subspace  of $\hom(W_{1}, \widehat{(W_{2})}_{z})$ spanned by
elements, denoted by $\phi$, satisfying the following conditions:

\begin{enumerate}


\item The {\it $\mathbf{d}$-conjugation property}: There exists $n\in \Z$ (called the {\it weight of $\phi$}
and denoted by $\wt \phi$) such that 
for $a\in \C^{\times}$ and $w_{1}\in W_{1}$,
$$a^{\mathbf{d}_{W_{2}}}(\phi(w_{1}))(z)=a^{n}(\phi(a^{\mathbf{d}_{W_{1}}}w_{1}))(az).$$

\item The {\it composability}: For $k, l\in \N$ and $v_{1}, \dots, v_{k+l}\in V$, $w_{1}\in W_{1}$
and $w_{2}'\in W_{2}'$,
the series 
$$\langle w_{2}', Y_{W_{2}}(v_{1}, z_{1})\cdots Y_{W_{2}}(v_{k}, z_{k})
(\phi(Y_{W_{1}}(v_{k+1}, z_{k+1})\cdots Y_{W_{1}}(v_{k+l}, z_{k+l})w_{1}))(z)\rangle$$
is absolutely convergent in the region $|z_{1}|>\cdots >|z_{k}|>|z|>|z_{k+1}|>\cdots |z_{k+l}|>0$
to a rational function 
\begin{equation}\label{r-fun}
R(\langle w_{2}', Y_{W_{2}}(v_{1}, z_{1})\cdots Y_{W_{2}}(v_{k}, z_{k})
(\phi(Y_{W_{1}}(v_{k+1}, z_{k+1})\cdots Y_{W_{1}}(v_{k+l}, z_{k+l})w_{1}))(z)\rangle)
\end{equation}
in $z_{1}, \dots, z_{k+l}$ and $z$ with the only possible 
poles $z_{i}=0$ for $i=1, \dots, k+l$, $z=0$, $z_{i}=z_{j}$ for $i, j=1, \dots, k+l$, $i\ne j$ and 
$z_{i}=z$ for $i=1, \dots, k+l$. Moreover, there 
exist $r_{i}\in \N$ depending only on the pair $(v_{i}, w_{1})$ for $i=1, \dots, k+l$, $m\in \N$
depending only on the pair $(\phi, w_{1})$, $p_{ij} \in \N$ depending only on the pair $(v_{i}, v_{j})$
for $i, j=1, \dots, k+l$, $i\ne j$, $s_{i}\in \N$ depending only on the pair $(v_{i}, \phi)$  for $i=1, \dots, k+l$
and $g(z_{1}, \dots, z_{k+l}, z)\in W_{2}[[z_{1}, \dots, z_{k+l}]]$ such that for $w_{2}'\in W_{2}'$, 
\begin{align*}
z^{m}\prod_{i=1}^{k+l}&z_{i}^{r_{i}}\prod_{1\le i<j\le k+l}(z_{i}-z_{j})^{p_{ij}}
\prod_{i=1}^{k+l}(z_{i}-z)^{s_{i}}\cdot\nn
&\cdot R(\langle w_{2}', Y_{W_{2}}(v_{1}, z_{1})
\cdots Y_{W_{2}}(v_{k}, z_{k})
(\phi(Y_{W_{1}}(v_{k+1}, z_{k+1})\cdots Y_{W_{1}}(v_{k+l}, z_{k+l})w_{1}))(z)\rangle)
\end{align*}
is a polynomial and is equal to 
$\langle w_{2}', g(z_{1}, \dots, z_{k+l}, z)\rangle$.

\end{enumerate}


Let $H_{[n]}$ be the 
subspace of $H$ consisting of elements of weight $n$. Then $H=\coprod_{n\in \Z}H_{[n]}$. 

Next we define the left and right vertex operator maps:
\begin{eqnarray*}
Y_{H}^{L}: V\otimes H&\to& H[[x, x^{-1}]]\nn
v\otimes \phi&\mapsto & Y_{H}^{L}(v, x)\phi,\\
Y_{H}^{R}: H\otimes V&\to& H[[x, x^{-1}]]\nn
\phi\otimes v&\mapsto & Y_{H}^{R}(\phi, x)v.
\end{eqnarray*}
Heuristically we would like to define them using the formulas
\begin{eqnarray}
\langle w_{2}', ((Y_{H}^{L}(v, z_{1})\phi)(w_{1}))(z_{2})\rangle&=&\langle w_{2}', 
Y_{W_{2}}(v, z_{1}+z_{2})(\phi(w_{1}))(z_{2})\rangle,\label{y-l-rigor}\\
\langle w_{2}',  ((Y_{H}^{R}(\phi, z_{1})v)(w_{1}))(z_{2})\rangle&=& 
\langle w_{2}',  \phi(Y_{W_{1}}(v, z_{2})w_{1})(z_{1}+z_{2})\rangle\label{y-r-rigor}
\end{eqnarray}
for $v\in V$, $\phi\in H$, $w_{1}\in W_{1}$ and $w_{2}'\in W_{2}'$.  But we need to make these
heuristic definitions precise. 

We first give the precise definition of $Y_{H}^{L}$. Let $\phi\in H$. 
Since $\phi$ satisfies the composability, for $v\in V$, $w_{1}\in W_{1}$ and $w_{2}'\in W_{2}'$,
$$\langle w_{2}', Y_{W_{2}}(v_{1}, z_{1}+z_{2})(\phi(w_{1}))(z_{2})\rangle$$
is absolutely 
convergent in the region given by $|z_{1}+z_{2}|>|z_{2}|>0$ to a rational 
function
$$R(\langle w_{2}', Y_{W_{2}}(v, z_{1}+z_{2})(\phi(w_{1}))(z_{2})\rangle)$$
in $z_{1}$ and $z_{2}$ with the only possible poles at $z_{1}=0$, $z_{2}=0$, $z_{1}+z_{2}=0$. 
Expanding this rational function in the region $|z_{2}|>|z_{1}|>0$, we obtain a lower truncated Laurent series 
$$\sum_{p\in \Z}a^{L}_{p}(w_{2}'\otimes v\otimes \phi \otimes w_{1}; z_{2})z_{1}^{-p-1}$$
in $z_{1}$. The coefficients $a^{L}_{p}(w_{2}'\otimes v\otimes \phi\otimes w_{1}; z_{2})$ for $p\in \Z$ 
of this Laurent series 
are in fact Laurent polynomials in $z_{2}$.  
On the other hand, by the composability, there exist $r, s, m\in \N$ 
and $g(z_{1}+z_{2}, z_{2})\in W_{2}[[z_{1}+z_{2}, z_{2}]]$
such that 
\begin{equation}\label{eta-L-H-0}
z_{2}^{m}(z_{1}+z_{2})^{r}z_{1}^{s}
R(\langle w_{2}', Y_{W_{2}}(v, z_{1}+z_{2})(\phi(w_{1}))(z_{2})\rangle)
=\langle w_{2}', g(z_{1}+z_{2}, z_{2})\rangle.
\end{equation}
For fixed  $z_{2}\ne 0$, the coefficients of the expansion in $z_{1}$ of the left-hand side of (\ref{eta-L-H-0}) 
in the region $|z_{2}|>|z_{1}|>0$ for $w_{2}'\in W_{2}'$ 
define elements of 
$(W_{2}')^{*}$. But since $g(z_{1}+z_{2}, z_{2})\in W_{2}[[z_{1}+z_{2}, z_{2}]]$, 
we see that these elements must be in the subspace $\overline{W_{2}}$. Multiplying 
these elements by $z_{2}^{-m}(z_{1}+z_{2})^{-r}z_{1}^{-s}$ in the region $|z_{2}|>|z_{1}|>0$,
we see that the results are still in $\overline{W_{2}}$, that is, 
for fixed $z_{2}\not=0$, 
the maps given by $w_{2}' \mapsto 
a^{L}_{p}(w_{2}'\otimes v\otimes \phi\otimes w_{1}; z_{2})$  for $p\in \Z$ are in fact elements 
of $\overline{W_{2}}$. When $z_{2}$ 
changes, we obtain elements of $\widehat{(W_{2})}_{z_{2}}$. Then for $v\in V$ and $\phi \in H$,
we have elements $\eta^{L}_{p; v, \phi}\in\hom(W_{1}, \widehat{(W_{2})}_{z})$ for $p\in \Z$ such that
$$\langle w_{2}', (\eta^{L}_{p; v, \phi}(w_{1}))(z)\rangle=a^{L}_{p}(w_{2}'\otimes v
\otimes \phi\otimes w_{1}; z)$$
for $w_{1}\in W_{1}$ and $w_{2}'\in W_{2}'$.

Similarly, since $\phi$ satisfies the composability, for $v\in V$, $w_{1}\in W_{1}$ and $w_{2}'\in W_{2}'$,
$$\langle w_{2}',  \phi(Y_{W_{1}}(v, z_{2})w_{1})(z_{1}+z_{2})\rangle$$
is absolutely 
convergent in the region given by $|z_{1}+z_{2}|>|z_{2}|>0$ to a rational 
function
$$R(\langle w_{2}',  \phi(Y_{W_{1}}(v, z_{2})w_{1})(z_{1}+z_{2})\rangle)$$
in $z_{1}$ and $z_{2}$ with the only possible poles at $z_{1}=0$, $z_{2}=0$, $z_{1}+z_{2}=0$. 
Expanding this rational function in the region $|z_{2}|>|z_{1}|>0$, we obtain a lower truncated Laurent series 
$$\sum_{p\in \Z}a^{R}_{p}(w_{2}'\otimes v\otimes w_{1}; z_{2})z_{1}^{-p-1}$$
in $z_{1}$.
The coefficients $a^{R}_{p}(w_{2}'\otimes v\otimes w_{1}; z_{2})$ for $p\in \Z$ of this Laurent series 
are in fact Laurent polynomials in $z_{2}$. For fixed $z_{2}\not=0$, 
the same argument as above shows that the maps given by $w_{2}' \mapsto 
a_{p}(w_{2}'\otimes v\otimes w_{1}; z_{2})$  for $p\in \Z$ are in fact elements of 
$\overline{W_{2}}$. When $z_{2}$ 
changes, we obtain elements of $\widehat{(W_{2})}_{z_{2}}$. The for $v\in V$ and $\phi \in H$,
we have elements $\eta^{R}_{p; \phi, v}\in\hom(W_{1}, \widehat{(W_{2})}_{z})$ for $p\in \Z$ such that
$$\langle w_{2}', \eta^{R}_{p; \phi, v}(w_{1})\rangle=a^{R}_{p}(w_{2}'\otimes v\otimes w_{1}; z)$$
for $w_{1}\in W_{1}$ and $w_{2}'\in W_{2}'$. 

\begin{prop}\label{eta-L-H}
The maps $\eta^{L}_{p; v, \phi}$ and $\eta^{R}_{p; \phi, v}$ are elements of $H$.  
When both $v$ and $\phi$ are homogeneous, 
$\eta^{L}_{p; v, \phi}$ and  $\eta^{R}_{p; \phi, v}$ are also homogeneous of weight 
$\wt v+\wt \phi-p-1$. In addition, if $\phi$ satisfies the $D$-derivative property
$$\frac{d}{dz}(\phi(w_{1}))(z)=D_{W_{2}}(\phi(w_{1}))(z)-(\phi(D_{W_{1}}w_{1}))(z),$$
or $w_{1}\in W_{1}$, where $D_{W_{2}}$ is the natural extension of $D_{W_{2}}$ on $W_{2}$
to $\overline{W_{2}}$, 
$\eta^{L}_{p; v, \phi}$ and $\eta^{R}_{p; \phi, v}$ also satisfy the $D$-derivative property
\begin{align*}
\frac{d}{dz_{2}}(\eta^{L}_{p; v, \phi}(w_{1}))(z_{2})
&=(D_{W_{2}}(\eta^{L}_{p; v, \phi}(w_{1}))(z_{2})
-(\eta^{L}_{p; v, \phi}(D_{W_{1}}w_{1}))(z_{2}))\nn
\frac{d}{dz_{2}}(\eta^{R}_{p; v, \phi}(w_{1}))(z_{2})
&=(D_{W_{2}}(\eta^{R}_{p; v, \phi}(w_{1}))(z_{2})
-(\eta^{R}_{p; v, \phi}(D_{W_{1}}w_{1}))(z_{2}))
\end{align*}
for $w_{1}\in W_{1}$. 
\end{prop}
\pf
We prove the result only for $\eta^{L}_{p; v, \phi}$. The proof for $\eta^{R}_{p; \phi, v}$ 
is similar and is omitted.

In the case that $v$ and $\phi$ are homogeneous, for $w_{1}\in W_{1}$, $w_{2}'\in W'_{2}$ and $a\in \C^{\times}$, 
in the region $|z_{2}|>|z_{1}|>0$, we have
\begin{align*}
\sum_{p\in \Z}&\langle w_{2}', a^{\mathbf{d}_{W_{2}}}(\eta^{L}_{p; v, \phi}(w_{1}))(z_{2})
\rangle z_{1}^{-p-1}\nn
&=\sum_{p\in \Z}\langle a^{d'_{W_{2}}}w_{2}', (\eta^{L}_{p; v, \phi}(w_{1}))(z_{2})
\rangle z_{1}^{-p-1}\nn
&=\sum_{p\in \Z}a^{L}_{p}(a^{d'_{W_{2}}}w_{2}'\otimes v\otimes \phi\otimes w_{1}; z_{2})z_{1}^{-p-1}\nn
&=R(\langle a^{d'_{W_{2}}}w_{2}', Y_{W_{2}}(v, z_{1}+z_{2})(\phi(w_{1}))(z_{2})\rangle)\nn
&=R(\langle w_{2}', a^{\mathbf{d}_{W_{2}}}Y_{W_{2}}(v, z_{1}+z_{2})(\phi(w_{1}))(z_{2})\rangle)\nn
&=R(\langle w_{2}', Y_{W_{2}}(a^{\mathbf{d}_{V}}v, az_{1}+az_{2})a^{\mathbf{d}_{W_{2}}}(\phi(w_{1}))(z_{2})\rangle)\nn
&=a^{\swt v+\swt \phi}R(\langle w_{2}', Y_{W_{2}}(v, az_{1}+az_{2})(\phi(a^{\mathbf{d}_{W_{1}}}w_{1}))(az_{2})\rangle)\nn
&=a^{\swt v+\swt \phi-p-1}\sum_{p\in \Z}a^{L}_{p}(w_{2}'\otimes v\otimes \phi
\otimes a^{\mathbf{d}_{W_{1}}}w_{1}; az_{2})z_{1}^{-p-1}\nn
&=a^{\swt v+\swt \phi-p-1}\sum_{p\in \Z}\langle w_{2}', (\eta^{L}_{p; v, \phi}( a^{\mathbf{d}_{W_{1}}}w_{1}))(az_{2})
\rangle z_{1}^{-p-1},
\end{align*}
proving the $\mathbf{d}$-conjugation property for $\eta^{L}_{p; v, \phi}$ with weight $\wt v+\wt \phi-p-1$.

Now we prove the composability of $\eta^{L}_{p; v, \phi}$. Since $\phi\in H$, it satisfies the composability. Then
$$\langle w_{2}', Y_{W_{2}}(v_{1}, z_{1})\cdots Y_{W_{2}}(v_{k}, z_{k})Y_{W_{2}}(v, z_{0}+z)
(\phi(Y_{W_{1}}(v_{k+1}, z_{k+1})\cdots Y_{W_{1}}(v_{k+l}, z_{k+l})w_{1}))(z)\rangle$$
is absolutely convergent in the region $|z_{1}|>\cdots |z_{k}|>|z_{0}+z|>|z|>|z_{k+1}|>\cdots >|z_{k+l}|>0$
to a rational function
\begin{equation}\label{eta-L-H-2}
R(\langle w_{2}', Y_{W_{2}}(v_{1}, z_{1})\cdots Y_{W_{2}}(v_{k}, z_{k})Y_{W_{2}}(v, z_{0}+z)
(\phi(Y_{W_{1}}(v_{k+1}, z_{k+1})\cdots Y_{W_{1}}(v_{k+l}, z_{k+l})w_{1}))(z)\rangle)
\end{equation}
 in $z_{0}, z_{1}, \dots, z_{k+l}$ and $z$ with the only possible 
poles $z_{i}=0$ for $i=1, \dots, k+l$, $z_{0}+z=0$, $z=0$,
$z_{i}=z_{j}$ for $i\ne j$, $z_{i}=z_{0}+z$, $z_{i}=z$ for $i=1, \dots, k+l$ and $z_{0}=0$. 


For $q\in \C$, let $\{e^{W_{2}}_{(q; \lambda)}\}_{\lambda\in \Lambda_{q}}$ be a basis of $(W_{2})_{[q]}$ and 
$\{(e^{W_{2}}_{(q; \lambda)})'\}_{\lambda\in \Lambda_{q}}$ be the subset of $W_{2}'$ defined by 
$$\langle (e^{W_{2}}_{(q; \lambda_{1})})', e^{W_{2}}_{(q; \lambda_{2})}\rangle=\delta_{\lambda_{1}\lambda_{2}}.$$
Putting $\{e^{W_{2}}_{(q; \lambda)}\}_{\lambda\in \Lambda_{q}}$ for $q\in \C$ together, we see that
$\{e^{W_{2}}_{(q; \lambda)}\}_{q\in \C,\;\lambda\in \Lambda_{q}}$ is a basis of $W_{2}$. 
For $q\in \C$, let $\pi_{q}$ be the projection from $W_{2}$ 
to $(W_{2})_{[q]}$ and use the same notation to denote its natural extension to 
$\overline{W_{2}}$. Then for $\overline{w}_{2}\in \overline{W_{2}}$, we have 
$$\pi_{q}\overline{w}_{2}=\sum_{\lambda\in \Lambda_{q}}\langle (e^{W_{2}}_{(q; \lambda)})', 
\overline{w}_{2}\rangle e^{W_{2}}_{(q; \lambda)}.$$

For $k, l\in \N$ and $v_{1}, \dots, v_{k+l}\in V$,  we know that 
\begin{align*}
\sum_{q\in \C}&
\langle w_{2}', Y_{W_{2}}(v_{1}, z_{1})\cdots Y_{W_{2}}(v_{k}, z_{k})\cdot \nn
&\quad\quad\quad\quad\quad\quad \cdot\pi_{q}Y_{W_{2}}(v, z_{0}+z)
(\phi(Y_{W_{1}}(v_{k+1}, z_{k+1})\cdots Y_{W_{1}}(v_{k+l}, z_{k+l})w_{1}))(z)\rangle\nn
&=\sum_{q\in \C}\Bigg(\sum_{\lambda\in \Lambda_{q}}
\langle w_{2}', Y_{W_{2}}(v_{1}, z_{1})\cdots Y_{W_{2}}(v_{k}, z_{k})
e^{W_{2}}_{(q; \lambda)}\rangle\cdot\nn
&\quad\quad\quad\quad\quad\quad \cdot \langle (e^{W_{2}}_{(q; \lambda)})', Y_{W_{2}}(v, z_{0}+z)
(\phi(Y_{W_{1}}(v_{k+1}, z_{k+1})\cdots Y_{W_{1}}(v_{k+l}, z_{k+l})w_{1}))(z)\rangle\Bigg),
\end{align*}
is absolutely convergent in the region $|z_{1}|> \cdots> |z_{k}|>|z_{0}+z|> |z|> |z_{k+1}|> \dots> |z_{k+l}|>0$
to the rational function (\ref{eta-L-H-2}) and 
$$\langle w_{2}', Y_{W_{2}}(v_{1}, z_{1})\cdots Y_{W_{2}}(v_{k}, z_{k})e^{W_{2}}_{(q; \lambda)}\rangle$$
and 
$$\langle (e^{W_{2}}_{(q; \lambda)})', Y_{W_{2}}(v, z_{0}+z)
(\phi(Y_{W_{1}}(v_{k+1}, z_{k+1})\cdots Y_{W_{1}}(v_{k+l}, z_{k+l})w_{1}))(z)\rangle$$
are absolutely convergent in the regions $|z_{1}|> \cdots> |z_{k}|$ and 
$|z_{0}+z|> |z|> |z_{k+1}|> \dots> |z_{k+l}|>0$, respectively, to the rational functions
$$R(\langle w_{2}', Y_{W_{2}}(v_{1}, z_{1})\cdots Y_{W_{2}}(v_{k}, z_{k})e^{W_{2}}_{(q; \lambda)}\rangle)$$
and 
\begin{equation}\label{eta-L-H-1.5}
R(\langle (e^{W_{2}}_{(q; \lambda)})', Y_{W_{2}}(v, z_{0}+z)
(\phi(Y_{W_{1}}(v_{k+1}, z_{k+1})\cdots Y_{W_{1}}(v_{k+l}, z_{k+l})w_{1}))(z)\rangle),
\end{equation}
respectively. 
So
\begin{align}\label{eta-L-H-1}
\sum_{q\in \C}&\Biggl(\sum_{\lambda\in \Lambda_{q}}
R(\langle w_{2}', Y_{W_{2}}(v_{1}, z_{1})\cdots Y_{W_{2}}(v_{k}, z_{k})e^{W_{2}}_{(q; \lambda)}\rangle)\cdot\nn
&\quad\quad\quad\quad\quad\quad \cdot R(\langle (e^{W_{2}}_{(q; \lambda)})', Y_{W_{2}}(v, z_{0}+z)
(\phi(Y_{W_{1}}(v_{k+1}, z_{k+1})\cdots Y_{W_{1}}(v_{k+l}, z_{k+l})w_{1}))(z)\rangle)\Biggr)
\end{align}
is absolutely convergent in the region $|z_{1}|> \cdots> |z_{k}|>|z_{0}+z|> |z|> |z_{k+1}|> 
\dots> |z_{k+l}|>0$  
to (\ref{eta-L-H-2}). But the expansion of the rational function (\ref{eta-L-H-2}) in the region $|z_{1}|, \dots, |z_{k}|>|z_{0}+z|, |z|, |z_{k+1}|$, 
$\dots, |z_{k+l}|>0$ is a series of rational functions of the same form as that of (\ref{eta-L-H-1}). 
By Lemma \ref{it-series-conv}, we see that 
this expansion must be equal to (\ref{eta-L-H-1}).
In particular,  (\ref{eta-L-H-1})
 is absolutely convergent in the region $|z_{1}|, \dots$, $|z_{k}|>|z_{0}+z|, |z|, |z_{k+1}|, \dots, |z_{k+l}|>0$
to the rational function (\ref{eta-L-H-2}). 

But each term in the right-hand side of (\ref{eta-L-H-1})
can be further expanded in the region $|z|>|z_{0}|>0$
and $|z_{i}-z|>|z_{0}|>0$ for $i=k+1, \dots, k+l$ as a Laurent series in $z_{0}$. In particular, in the region 
$|z_{1}|>\cdots>|z_{k}|>|z|> |z_{k+1}|>\cdots> |z_{k+l}|>0$, $|z|>|z_{0}|>0$
and $|z_{i}-z|>|z_{0}|>0$ for $i=k+1, \dots, k+l$, the series in the 
right-hand side of (\ref{eta-L-H-1}) is equal to
\begin{align}\label{eta-L-H-3}
\sum_{p\in \Z}&\sum_{q\in \C}\Biggl(\sum_{\lambda\in \Lambda_{q}}
\langle w_{2}', Y_{W_{2}}(v_{1}, z_{1})\cdots Y_{W_{2}}(v_{k}, z_{k})e^{W_{2}}_{(q; \lambda)}\rangle\cdot\nn
&\quad\quad\quad\quad\quad\quad \cdot a_{p}^{L}((e^{W_{2}}_{(q; \lambda)})'\otimes v\otimes \phi
\otimes Y_{W_{1}}(v_{k+1}, z_{k+1})\cdots Y_{W_{1}}(v_{k+l}, z_{k+l})w_{1}; z)\Biggr)
z_{0}^{-p-1}\nn
&=\sum_{p\in \Z}\sum_{q\in \C}\Biggl(\sum_{\lambda\in \Lambda_{q}}
\langle w_{2}', Y_{W_{2}}(v_{1}, z_{1})\cdots Y_{W_{2}}(v_{k}, z_{k})e^{W_{2}}_{(q; \lambda)}\rangle\cdot\nn
&\quad\quad\quad\quad\quad\quad \cdot \langle (e^{W_{2}}_{(q; \lambda)})', 
(\eta^{L}_{p; v, \phi}(Y_{W_{1}}(v_{k+1}, z_{k+1})\cdots 
Y_{W_{1}}(v_{k+l}, z_{k+l})w_{1}))(z)\rangle\Biggr)
z_{0}^{-p-1}\nn
&=\sum_{p\in \Z}
\langle w_{2}', Y_{W_{2}}(v_{1}, z_{1})\cdots Y_{W_{2}}(v_{k}, z_{k})\cdot \nn
&\quad\quad\quad\quad\quad\quad \cdot
(\eta^{L}_{p; v, \phi}(Y_{W_{1}}(v_{k+1}, z_{k+1})\cdots Y_{W_{1}}(v_{k+l}, z_{k+l})w_{1}))(z)\rangle
z_{0}^{-p-1}.
\end{align}
Thus each term in the right-hand side of (\ref{eta-L-H-3}) is absolutely convergent 
in the region given by $|z_{1}|>\cdots>|z_{k}|>|z|> |z_{k+1}|>\cdots> |z_{k+l}|>0$
to a rational function in $z_{1}, \dots, z_{k+l}$ and $z$ with the only possible 
poles $z_{i}=0$ for $i=1, \dots, k+l$, $z=0$, $z_{i}=z_{j}$ for $i, j=1, \dots, k+l$, $i\ne j$ and 
$z_{i}=z$ for $i=1, \dots, k+l$, proving the first part of the composability of $\eta^{L}_{p; v, \phi}$. 
Since $\phi$ satisfies the second part of the composability
and the right-hand side of (\ref{eta-L-H-3}) is absolutely convergent to (\ref{eta-L-H-2}), it is clear that
$\eta^{L}_{p; v, \phi}$ also satisfies the second part of 
the composability. 

Finally when $\phi$ satisfies the $D$-derivative property, 
for $w_{1}\in W_{1}$ and $w_{2}'\in W'_{2}$, in the region $|z_{2}|>|z_{1}|>0$, we have
\begin{align*}
\sum_{p\in \Z}&\left\langle w_{2}', \frac{d}{dz_{2}}(\eta^{L}_{p; v, \phi}(w_{1}))(z_{2})
\right\rangle z_{1}^{-p-1}\nn
&=\frac{d}{dz_{2}}\sum_{p\in \Z}a^{L}_{p}(w_{2}'\otimes v\otimes \phi\otimes w_{1}; z_{2})z_{1}^{-p-1}\nn
&=\frac{d}{dz_{2}}R(\langle w_{2}', Y_{W_{2}}(v, z_{1}+z_{2})(\phi(w_{1}))(z_{2})\rangle)
\quad\quad\quad\quad\quad\quad\quad\quad\quad\quad\quad\quad\quad\quad\quad\quad
\quad\quad\quad\quad
\end{align*}
\begin{align*}
&\;\:=R\left(\left\langle w_{2}', \left(\frac{d}{dz_{2}}Y_{W_{2}}(v, z_{1}+z_{2})\right)(\phi(w_{1}))(z_{2})
\right\rangle\right)\nn
&\;\:\quad +R\left(\left\langle w_{2}', Y_{W_{2}}(v, z_{1}+z_{2})\left(\frac{d}{dz_{2}}(\phi(w_{1}))(z_{2})\right)\right\rangle
\right)\nn
&\;\:=R(\langle w_{2}', D_{W_{2}}Y_{W_{2}}(v, z_{1}+z_{2})(\phi(w_{1}))(z_{2})\rangle)
-R(\langle w_{2}', Y_{W_{2}}(v, z_{1}+z_{2})D_{W_{2}}(\phi(w_{1}))(z_{2})\rangle)\nn
&\;\: \quad +R(\langle w_{2}', Y_{W_{2}}(v, z_{1}+z_{2})D_{W_{2}}(\phi(w_{1}))(z_{2})\rangle)
-R(\langle w_{2}', Y_{W_{2}}(v, z_{1}+z_{2})(\phi(D_{W_{1}}w_{1}))(z_{2})\rangle)\nn
&\;\:=R(\langle D_{W_{2}}'w_{2}', Y_{W_{2}}(v, z_{1}+z_{2})(\phi(w_{1}))(z_{2})\rangle)
-R(\langle w_{2}', Y_{W_{2}}(v, z_{1}+z_{2})(\phi(D_{W_{1}}w_{1}))(z_{2})\rangle)\nn
&\;\:=\sum_{p\in \Z}a_{p}^{L}(D_{W_{2}}'w_{2}'\otimes v\otimes \phi\otimes w_{1}; z_{2})z_{1}^{-p-1}
-\sum_{p\in \Z}a_{p}^{L}(w_{2}'\otimes v\otimes \phi\otimes D_{W_{1}}w_{1}; z_{2})z_{1}^{-p-1}\nn
&\;\:=\sum_{p\in \Z}\langle D_{W_{2}}'w_{2}', (\eta^{L}_{p; v, \phi}(w_{1}))(z_{2})\rangle
-\sum_{p\in \Z}\langle w_{2}', (\eta^{L}_{p; v, \phi}(D_{W_{1}}w_{1}))(z_{2})
\rangle z_{1}^{-p-1}\nn
&\;\:=\sum_{p\in \Z}\langle w_{2}', (D_{W_{2}}(\eta^{L}_{p; v, \phi}(w_{1}))(z_{2})
-(\eta^{L}_{p; v, \phi}(D_{W_{1}}w_{1}))(z_{2}))
\rangle z_{1}^{-p-1},
\end{align*}
proving the $D$-derivative property for $\eta^{L}_{p; v, \phi}$.
\epfv

For $v\in V$ and $\phi\in H$, we define
$$Y_{H}^{L}(v, x)\phi=\sum_{p\in \Z}\eta^{L}_{p; v, \phi}x^{-p-1}$$
and 
$$Y_{H}^{R}(\phi, x)v=\sum_{p\in \Z}\eta^{R}_{p; \phi, v}x^{-p-1}.$$
Using our notations for components of vertex operators, we have
$(Y_{H}^{L})_{p}(v)\phi=\eta^{L}_{p; v, \phi}$ and 
$(Y_{H}^{R})_{p}(\phi)v=\eta^{R}_{p; \phi, v}$.

The proof of Proposition \ref{eta-L-H} in fact has also proved the first half of the following result (the second half 
can be proved similarly and its proof is omitted):

\begin{prop}\label{part-y-h-l-r}
Let $k, l\in \N$, $v_{1}, \dots, v_{k+l}\in V$, $w_{1}\in W_{1}$ and $w_{2}'\in W'_{2}$.
Then the series 
\begin{equation}\label{part-y-h-l-r-1}
\langle w_{2}', Y_{W_{2}}(v_{1}, z_{1})\cdots Y_{W_{2}}(v_{k}, z_{k})
((Y_{H}^{L}(v, z_{0})\phi)(Y_{W_{1}}(v_{k+1}, z_{k+1})\cdots Y_{W_{1}}(v_{k+l}, z_{k+l})w_{1}))(z)\rangle
\end{equation}
is absolutely convergent in the region $|z_{1}|>\cdots |z_{k}|>|z|>|z_{k+1}|>\cdots >|z_{k+l}|$,
$|z|>|z_{0}|>0$, $|z_{i}|>|z+z_{0}|>0$ for $i=1, \dots, k$, $|z_{i}-z|>|z_{0}|>0$ for
$i=k+1, \dots, k+l$ to the rational function 
\begin{equation}\label{part-y-h-l-r-2}
R(\langle w_{2}', Y_{W_{2}}(v_{1}, z_{1})\cdots Y_{W_{2}}(v_{k}, z_{k})Y_{W_{2}}(v, z_{0}+z)
(\phi(Y_{W_{1}}(v_{k+1}, z_{k+1})\cdots Y_{W_{1}}(v_{k+l}, z_{k+l})w_{1}))(z)\rangle). 
\end{equation}
Similarly, 
the series 
$$\langle w_{2}', Y_{W_{2}}(v_{1}, z_{1})\cdots Y_{W_{2}}(v_{k}, z_{k})
((Y_{H}^{R}(\phi, z_{0})v)(Y_{W_{1}}(v_{k+1}, z_{k+1})\cdots Y_{W_{1}}(v_{k+l}, z_{k+l})w_{1}))(z)\rangle$$
is absolutely convergent in the region $|z_{1}|>\cdots |z_{k}|>|z|>|z_{k+1}|>\cdots >|z_{k+l}|$,
$|z|>|z_{0}|>0$, $|z_{i}|>|z+z_{0}|>0$ for $i=1, \dots, k$, $|z_{i}-z|>|z_{0}|>0$ for
$i=k+1, \dots, k+l$ to the rational function 
$$R(\langle w_{2}', Y_{W_{2}}(v_{1}, z_{1})\cdots Y_{W_{2}}(v_{k}, z_{k})
(\phi(Y_{W_{1}}(v, z)Y_{W_{1}}(v_{k+1}, z_{k+1})\cdots Y_{W_{1}}(v_{k+l}, z_{k+l})w_{1}))(z_{0}+z)\rangle).$$ 
In particular, (\ref{y-l-rigor}) and  (\ref{y-r-rigor}) hold
in the region $|z_{1}+z_{2}|>|z_{2}|>z_{1}|>0$
and 
\begin{align}
R(\langle w_{2}', ((Y_{H}^{L}(v, z_{1})\phi)(w_{1}))(z_{2})\rangle)&=R(\langle w_{2}', 
Y_{W_{2}}(v, z_{1}+z_{2})(\phi(w_{1}))(z_{2})\rangle),\label{y-l-rigor-r}\\
R(\langle w_{2}',  ((Y_{H}^{R}(\phi, z_{1})v)(w_{1}))(z_{2})\rangle)&=
R(\langle w_{2}',  \phi(Y_{W_{1}}(v, z_{2})w_{1})(z_{1}+z_{2})\rangle).\label{y-r-rigor-r}
\end{align}
\end{prop}

We also have the following result:

\begin{prop}\label{y-h-l-r}
Let $w_{2}'\in W_{2}'$, $w_{1}\in W_{1}$, $v_{1}, \dots, v_{k+l}\in V$ and 
$\phi\in H$. Then 
\begin{equation}\label{y-h-l-r-1}
\langle w_{2}', ((Y_{H}^{L}(v_{1}, z_{1})\cdots Y_{H}^{L}(v_{k}, z_{k})Y_{H}^{R}(\phi, z)
Y_{V}(v_{k+1}, z_{k+1})\cdots Y_{V}(v_{k+l-1}, z_{k+l-1})
v_{k+l})(w_{1}))(z_{k+l})\rangle
\end{equation}
is absolutely convergent in the region $|z_{k+l}|>|z_{1}|>\cdots>|z_{k}|
>|z|>|z_{k+1}|>\cdots >|z_{k+l-1}|>0$ to the same rational function in $z_{1}, \dots, z_{k+l}$ and $z$
as the rational function that  the series 
\begin{align}\label{y-h-l-r-2}
\langle w_{2}', &Y_{W_{2}}(v_{1}, z_{1}+z_{k+l})\cdots Y_{W_{2}}(v_{k}, z_{k}+z_{k+l})\cdot \nn
&\cdot 
\phi(Y_{W_{1}}(v_{k+1}, z_{k+1}+z_{k+l})\cdots Y_{W_{1}}(v_{k+l-1}, z_{k+l-1}+z_{k+l})
Y_{W_{1}}(v_{k+l}, z_{k+l})w_{1})(z+z_{k+l})\rangle.
\end{align}
is absolutely convergent to in the region $|z_{1}+z_{k+l}|>\cdots >|z_{k}+z_{k+l}|
>|z+z_{k+l}|>|z_{k+1}+z_{k+l}|>\cdots 
|z_{k+l-1}+z_{k+l}|>|z_{k+l}|>0$. In particular, we have the equality
\begin{align}\label{y-h-l-r-3}
R&(\langle w_{2}', ((Y_{H}^{L}(v_{1}, z_{1})\cdots Y_{H}^{L}(v_{k}, z_{k})\cdot\nn
&\quad \cdot Y_{H}^{R}(\phi, z)
Y_{V}(v_{k+1}, z_{k+1})\cdots Y_{V}(v_{k+l-1}, z_{k+l-1})
v_{k+l})(w_{1}))(z_{k+l})\rangle)\nn
&=R(\langle w_{2}', Y_{W_{2}}(v_{1}, z_{1}+z_{k+l})\cdots Y_{W_{2}}(v_{k}, z_{k}+z_{k+l})\cdot \nn
&\quad \cdot 
\phi(Y_{W_{1}}(v_{k+1}, z_{k+1}+z_{k+l})\cdots Y_{W_{1}}(v_{k+l-1}, z_{k+l-1}+z_{k+l})
Y_{W_{1}}(v_{k+l}, z_{k+l})w_{1})(z+z_{k+l})\rangle)
\end{align}
of rational functions 
and the only possible poles of the rational function 
(in both sides of) (\ref{y-h-l-r-3})  are $z_{i}+z_{k+l}=0$ for $i=1, \dots, k+l-1$, 
$z+z_{k+l}=0$,  $z_{i}=z_{j}$ for $ i, j=1, 
\dots, k+l-1$, $i\ne j$,
$z_{i}=0$ for $i=1, \dots,  k+l$, $z_{i}=z$ for $i=1, \dots,  k+l-1$
and $z=0$. More explicitly, (\ref{y-h-l-r-3}) is of the form 
\begin{align}\label{rationality-3}
g(z_{1},& \dots, z_{k+l}, z)z^{-m}(z+z_{k+l})^{-n}\prod_{i=1}^{k+l-1}(z_{i}+z_{k+l})^{-p_{i}}\cdot \nn
&\cdot 
\prod_{1\le i< j\le k+l-1}(z_{i}-z_{j})^{-p_{ij}}
 \prod_{i=1}^{k+l}z_{i}^{-r_{i}}
 \prod_{i=1}^{k+l-1}(z_{i}-z)^{-s_{i}}
\end{align}
where $g(z_{1}, \dots, z_{k+l}, z)$ is a polynomial in $z_{1}, \dots, z_{k+l}$ and $z$
and $m, n, p_{i}, p_{ij}, r_{i}, s_{i}\in \N$ depend only on the pairs $(v_{k+l}, \phi)$, 
$(\phi, w_{1})$, $(v_{i}, w_{1})$, $(v_{i}, v_{j})$, $(v_{i}, v_{k+l})$, $(v_{i}, \phi)$, respectively. 
\end{prop}
\pf
Since $\phi$ satisfies the composability, the series
\begin{align}\label{rationality-1}
\langle w_{2}', &Y_{W_{2}}(v_{1}, z_{1}+z_{k+l})\cdots Y_{W_{2}}(v_{k}, z_{k}+z_{k+l})\cdot \nn
&\cdot 
\phi(Y_{W_{1}}(v_{k+1}, z_{k+1}+z_{k+l})\cdots Y_{W_{1}}(v_{k+l-1}, z_{k+l-1}+z_{k+l})
Y_{W_{1}}(v_{k+l}, z_{k+l})w_{1})(z+z_{k+l})\rangle.
\end{align}
is absolutely convergent in the region $|z_{1}+z_{k+l}|>\cdots >|z_{k}+z_{k+l}|
>|z+z_{k+l}|>|z_{k+1}+z_{k+l}|>\cdots 
|z_{k+l-1}+z_{k+l}|>|z_{k+l}|>0$
to a rational function 
\begin{align}\label{rationality-2}
R(\langle w_{2}', &Y_{W_{2}}(v_{1}, z_{1}+z_{k+l})\cdots Y_{W_{2}}(v_{k}, z_{k}+z_{k+l})\cdot \nn
&\cdot 
\phi(Y_{W_{1}}(v_{k+1}, z_{k+1}+z_{k+l})\cdots Y_{W_{1}}(v_{k+l-1}, z_{k+l-1}+z_{k+l})
Y_{W_{1}}(v_{k+l}, z_{k+l})w_{1})(z+z_{k+l})\rangle).
\end{align}
in $z_{1}, \dots, z_{k+l}$ and $z$ with the only possible 
poles $z_{i}+z_{k+l}=0$ for $i=1, \dots, k+l-1$, 
$z+z_{k+l}=0$,  $z_{i}=z_{j}$ for $ i, j=1, 
\dots, k+l-1$, $i\ne j$,
$z_{i}=0$ for $i=1, \dots,  k+l$, $z_{i}=z$ for $i=1, \dots,  k+l-1$
and $z=0$.
The rational function (\ref{rationality-2}) is of the form (\ref{rationality-3}) such that
$g(z_{1}, \dots, z_{k+l}, z)$ is a polynomial in $z_{1}, \dots, z_{k+l}$ and $z$
and $m, n, p_{i}, p_{ij}, r_{i}, s_{i}\in \N$ depend only on the pairs $(v_{k+l}, \phi)$, 
$(\phi, w_{1})$, $(v_{i}, w_{1})$, $(v_{i}, v_{j})$, $(v_{i}, v_{k+l})$, $(v_{i}, \phi)$, respectively.  .

Let $\{e^{W_{1}}_{(q; \gamma)}\}_{q\in \C,\; \gamma\in \Gamma_{q}}$ be a homogeneous basis 
of $W_{1}$ and $\{(e^{W_{1}}_{(q; \gamma)})'\}_{q\in \C,\; \gamma\in \Gamma_{q}}$ 
the subset of $W_{1}'$ given by
$$\langle (e^{W_{1}}_{(q; \gamma_{1})})', e^{W_{1}}_{(q; \gamma_{2})}\rangle=\delta_{\gamma_{1}\gamma_{2}}.$$
By the associativity of vertex operators, for $q\in \C$, $\gamma\in \Gamma_{q}$, 
\begin{align*}
\langle (e^{W_{1}}_{(q; \gamma)})', &Y_{W_{1}}(v_{k+1}, z_{k+1}+z_{k+l})\cdots 
Y_{W_{1}}(v_{k+l-1}, z_{k+l-1}+z_{k+l})
Y_{W_{1}}(v_{k+l}, z_{k+l})w_{1}\rangle\nn
&=\langle (e^{W_{1}}_{(q; \gamma)})', Y_{W_{1}}(Y_{V}(v_{k+1}, z_{k+1})\cdots 
Y_{V}(v_{k+l-1}, z_{k+l-1})v_{k+l}, z_{k+l})w_{1}\rangle
\end{align*}
in the region $|z_{k+1}+z_{k+l}|>\cdots >|z_{k+l-1}+z_{k+l}|>|z_{k+l}|>0$, 
$|z_{k+l}|>|z_{k+1}|>\cdots >|z_{k+l-1}|>0$. Then by the convergence of (\ref{rationality-1}),
we see that 
\begin{align}\label{rationality-4}
\langle w_{2}', &Y_{W_{2}}(v_{1}, z_{1}+z_{k+l})\cdots Y_{W_{2}}(v_{k}, z_{k}+z_{k+l})\cdot \nn
&\quad\cdot 
\phi(Y_{W_{1}}(Y_{V}(v_{k+1}, z_{k+1})\cdots Y_{V}(v_{k+l-1}, z_{k+l-1})
v_{k+l}, z_{k+l})w_{1})(z+z_{k+l})\rangle
\end{align}
is convergent absolutely to the same rational function (\ref{rationality-2}) or (\ref{rationality-3}) in the region
$|z_{1}+z_{k+l}|>\cdots >|z_{k}+z_{k+l}|>|z+z_{k+l}|>|z_{k+2}+z_{k+l}|>\cdots 
|z_{k+l-1}+z_{k+l}|>0$, $|z_{k+l}|>|z_{k+1}|>\cdots >|z_{k+l-1}|>0$. 

Expand (\ref{rationality-3}) as a Laurent series in $z_{i}+z_{k+l}$ for $i=1, \dots, k$,
$z_{k+j}$ for $j=1, 
\dots, l$, $z+z_{k+l}$ in the region $|z+z_{k+l}|>|z_{k+l}|>|z_{k+1}|>\cdots >|z_{k+l-1}|>0$, 
$|z_{1}+z_{k+l}|>\cdots>|z_{k}+z_{k+l}|>|z_{k+1}+z_{k+l}|, \dots, |z_{k+l-1}+z_{k+l}|,
|z_{k+l}|, |z+z_{k+l}|>0$, $|z+z_{k+l}|>|z_{k+1}+z_{k+l}|, \dots, |z_{k+l-1}+z_{k+l}|>0$.
Since (\ref{rationality-4}) is a Laurent series of the same form, by Lemma \ref{it-series-conv},
it must be convergent absolutely to 
(\ref{rationality-3}) in the same region.

Let $\{e^{V}_{(q; \mu)}\}_{q\in \Z,\; \mu\in M_{q}}$ be a 
homogeneous basis of $V$ and $\{(e^{V}_{(q; \mu)})'\}_{q\in \Z,\;\mu\in M_{q}}$
a subset of $V'$ given by 
$$\langle (e^{V}_{(q; \mu_{1})})', e^{V}_{(q; \mu_{2})}\rangle=\delta_{\mu_{1}\mu_{2}}.$$ 
Let $\{e^{W_{2}}_{(q; \lambda)}\}_{\lambda\in \Lambda_{q}}$ be a basis of $(W_{2})_{(q)}$ and  
$\{(e^{W_{2}}_{(q; \lambda)})'\}_{\lambda\in \Lambda_{q}}$ the subset of $W_{2}'$ 
as above. 
Then the series (\ref{rationality-4}) can be written as
\begin{align}\label{rationality-5}
\sum_{q_{1}\in \Z}\sum_{q_{2}\in \Z}\Biggl(\sum_{\mu\in M_{q_{1}}}\sum_{\lambda\in \Lambda_{q_{2}}}
&\langle w_{2}', Y_{W_{2}}(v_{1}, z_{1}+z_{k+l})\cdots 
Y_{W_{2}}(v_{k}, z_{k}+z_{k+l})e^{W_{2}}_{(q_{2}; \lambda)}\rangle\cdot \nn
&\cdot 
\langle (e^{W_{2}}_{(q_{2}; \lambda)})', \phi(Y_{W_{1}}(e^{V}_{(q_{1}; i)}, 
z_{k+l})w_{1})(z+z_{k+l})\rangle\cdot\nn
&\cdot \langle (e^{V}_{(q_{1}; \mu)})', 
Y_{V}(v_{k+1}, z_{k+1})\cdots Y_{V}(v_{k+l-1}, z_{k+l-1})
v_{k+l}\rangle\Biggr).
\end{align}
Because of (\ref{y-r-rigor}), we can replace 
$$\langle (e^{W_{2}}_{(q_{2}; \lambda)})', 
\phi(Y_{W_{1}}(e^{V}_{(q_{1}; \mu)}, z_{k+l})w_{1})(z+z_{k+l})\rangle$$
in (\ref{rationality-5}) by 
$$\langle (e^{W_{2}}_{(q_{2}; \lambda)})',  ((Y_{H}^{R}(\phi, z)
e^{V}_{(q_{1}; \mu)})(w_{1}))(z_{k+l})\rangle$$
when $|z+z_{k+l}|>|z_{k+l}|>|z|>0$. Thus 
\begin{align}\label{rationality-6}
\sum_{q_{1}\in \Z}&\sum_{q_{2}\in \Z}\Biggl(\sum_{\mu\in M_{q_{1}}}\sum_{\gamma\in \Gamma_{q_{2}}}
\langle w_{2}', Y_{W_{2}}(v_{1}, z_{1}+z_{k+l})\cdots 
Y_{W_{2}}(v_{k}, z_{k}+z_{k+l})e^{W_{2}}_{(q_{2}; \lambda)}\rangle\cdot \nn
&\quad \cdot 
\langle (e^{W_{2}}_{(q_{2}; \lambda)})',  ((Y_{H}^{R}(\phi, z)
e^{V}_{(q_{1}; \mu)})(w_{1}))(z_{k+l})\rangle\cdot\nn
&\quad \cdot \langle (e^{V}_{(q_{1}; \mu)})', 
Y_{V}(v_{k+1}, z_{k+1})\cdots Y_{V}(v_{k+l-1}, z_{k+l-1})
v_{k+l}\rangle\Biggr)\nn
&=\langle w_{2}', Y_{W_{2}}(v_{1}, z_{1}+z_{k+l})\cdots Y_{W_{2}}(v_{k}, z_{k}+z_{k+l})\cdot \nn
&\quad\cdot 
((Y_{H}^{R}(\phi, z)Y_{V}(v_{k+1}, z_{k+1})\cdots Y_{V}(v_{k+l-1}, z_{k+l-1})
v_{k+l})(w_{1}))(z_{k+l})\rangle
\end{align}
is absolutely convergent to (\ref{rationality-3}) in the region $|z+z_{k+l}|>|z_{k+l}|>|z|>
|z_{k+1}|>\cdots >|z_{k+l-1}|>0$, 
$|z_{1}+z_{k+l}|>\cdots>|z_{k}+z_{k+l}|>|z_{k+1}+z_{k+l}|, \dots, |z_{k+l-1}+z_{k+l}|,
|z_{k+l}|, |z+z_{k+l}|>0$, $|z+z_{k+l}|>|z_{k+1}+z_{k+l}|, \dots, |z_{k+l-1}+z_{k+l}|>0$. 
But the form of the Laurent series (\ref{rationality-6}) is the same as the form of 
the expansion of (\ref{rationality-3}) as a Laurent series in $z_{i}+z_{k+l}$ for $i=1, \dots, k$,
$z_{k+j}$ for $j=1, 
\dots, l$, $z$ in the larger region $|z_{k+l}|>|z|>|z_{k+1}|>\cdots >|z_{k+l-1}|>0$, 
$|z_{1}+z_{k+l}|>\cdots>|z_{k}+z_{k+l}|>|z_{k+1}+z_{k+l}|, \dots, |z_{k+l-1}+z_{k+l}|,
|z_{k+l}|, |z+z_{k+l}|>0$, $|z+z_{k+l}|>|z_{k+1}+z_{k+l}|, \dots, |z_{k+l-1}+z_{k+l}|>0$.
Thus by Lemma \ref{it-series-conv}, (\ref{rationality-6})
is in fact absolutely convergent to (\ref{rationality-3}) in this larger region. 

The series in the right-hand side of (\ref{rationality-6}) can be written as
\begin{align*}
\sum_{q_{1}\in \Z}&\sum_{q_{2}\in \Z}\Biggl(\sum_{\mu\in M_{q_{1}}}\sum_{\gamma\in \Gamma_{q_{2}}}
\langle w_{2}', Y_{W_{2}}(v_{1}, z_{1}+z_{k+l})\cdots Y_{W_{2}}(v_{k-1}, z_{k-1}+z_{k+l})
e^{W_{2}}_{(q_{2}; \lambda)}
\rangle
\cdot \nn
&\quad\quad\quad \cdot 
\langle (e^{W_{2}}_{(q_{2}; \lambda)})',  Y_{W_{2}}(v_{k}, z_{k}+z_{k+l})
((Y_{H}^{R}(e^{V}_{(q_{1}; \mu)}, z)\phi)(w_{1}))(z_{k+l})\rangle\cdot\nn
&\quad\quad\quad \cdot \langle (e^{V}_{(q_{1}; \mu)})', 
Y_{V}(v_{k+1}, z_{k+1})\cdots Y_{V}(v_{k+l-1}, z_{k+l-1})
v_{k+l}\rangle\Biggr)\quad\quad\quad\quad\quad\quad\quad\quad\quad\quad\;\;
\end{align*}
\begin{align}\label{rationality-7}
&=\sum_{q_{1}\in \Z}\sum_{q_{2}\in \Z}\sum_{p\in \Z}\Biggl(\sum_{\mu\in M_{q_{1}}}
\sum_{\gamma\in \Gamma_{q_{2}}}
\langle w_{2}', Y_{W_{2}}(v_{1}, z_{1}+z_{k+l})\cdots Y_{W_{2}}(v_{k-1}, z_{k-1}+z_{k+l})
e^{W_{2}}_{(q_{2}; \lambda)}\rangle
\cdot \nn
&\quad\quad\quad  \cdot 
\langle (e^{W_{2}}_{(q_{2}; \lambda)})',  Y_{W_{2}}(v_{k}, z_{k}+z_{k+l})
(\eta^{R}_{p; \phi, e^{V}_{(q_{1}; \mu)}}(w_{1}))(z_{k+l})\rangle\cdot\nn
&\quad\quad\quad \cdot \langle (e^{V}_{(q_{1}; \mu)})', 
Y_{V}(v_{k+1}, z_{k+1})\cdots Y_{V}(v_{k+l-1}, z_{k+l-1})
v_{k+l}\rangle\Biggr) z^{-p-1}
\end{align}
Because of (\ref{y-l-rigor}), we can replace 
$$\langle (e^{W_{2}}_{(q_{2}; \lambda)})',  Y_{W_{2}}(v_{k}, z_{k}+z_{k+l})
(\eta^{R}_{p; e^{V}_{(q_{1}; \mu)}, \phi}(w_{1}))(z_{k+l})\rangle$$
in (\ref{rationality-7}) by 
$$\langle (e^{W_{2}}_{(q_{2}; \lambda)})',  
((Y_{H}^{L}(v_{k}, z_{k})\eta^{R}_{p; \phi}, e^{V}_{(q_{1}; \mu)})(w_{1}))(z_{k+l})\rangle$$
when $|z_{k}+z_{k+l}|>|z_{k+l}|>|z_{k}|>0$. Thus 
\begin{align}\label{rationality-8}
\sum_{q_{1}\in \Z}&\sum_{q_{2}\in \Z}\sum_{p\in \Z}\Biggl(\sum_{\mu\in M_{q_{1}}}
\sum_{\gamma\in \Gamma_{q_{2}}}
\langle w_{2}', Y_{W_{2}}(v_{1}, z_{1}+z_{k+l})\cdots Y_{W_{2}}(v_{k-1}, z_{k-1}+z_{k+l})
e^{W_{2}}_{(q_{2}; \lambda)}\rangle
\cdot \quad\quad\quad\quad\quad\quad\nn
&\quad\quad \cdot 
\langle (e^{W_{2}}_{(q_{2}; \lambda)})',  
((Y_{H}^{L}(v_{k}, z_{k})\eta^{R}_{p; \phi, e^{V}_{(q_{1}; \mu)}})(w_{1}))(z_{k+l})\rangle\cdot\nn
&\quad\quad \cdot \langle (e^{V}_{(q_{1}; \mu)})', 
Y_{V}(v_{k+1}, z_{k+1})\cdots Y_{V}(v_{k+l-1}, z_{k+l-1})
v_{k+l}\rangle\Biggr) z^{-p-1}\nn
&=\sum_{q_{1}\in \Z}\sum_{q_{2}\in \Z}\Biggl(\sum_{\mu\in M_{q_{1}}}
\sum_{\gamma\in \Gamma_{q_{2}}}
\langle w_{2}', Y_{W_{2}}(v_{1}, z_{1}+z_{k+l})\cdots Y_{W_{2}}(v_{k-1}, z_{k-1}+z_{k+l})
e^{W_{2}}_{(q_{2}; \lambda)}\rangle
\cdot \nn
&\quad\quad \cdot 
\langle (e^{W_{2}}_{(q_{2}; \lambda)})',  
((Y_{H}^{L}(v_{k}, z_{k})Y_{H}^{R}(\phi, z)e^{V}_{(q_{1}; \mu)})(w_{1}))(z_{k+l})\rangle\cdot\nn
&\quad\quad \cdot \langle (e^{V}_{(q_{1}; \mu)})', 
Y_{V}(v_{k+1}, z_{k+1})\cdots Y_{V}(v_{k+l-1}, z_{k+l-1})
v_{k+l}\rangle\Biggr)\nn
&=\langle w_{2}', Y_{W_{2}}(v_{1}, z_{1}+z_{k+l})\cdots Y_{W_{2}}(v_{k-1}, z_{k-1}+z_{k+l})\cdot\nn 
&\quad\quad \cdot ((Y_{H}^{L}(v_{k}, z_{k})Y_{H}^{R}(\phi, z)
Y_{V}(v_{k+1}, z_{k+1})\cdots Y_{V}(v_{k+l-1}, z_{k+l-1})
v_{k+l})(w_{1}))(z_{k+l})\rangle
\end{align}
is absolutely convergent to (\ref{rationality-3}) in the larger region 
$|z_{k}+z_{k+l}|>|z_{k+l}|>|z_{k}|>|z|>|z_{k+1}|>\cdots >|z_{k+l-1}|>0$, 
$|z_{1}+z_{k+l}|>\cdots>|z_{k-1}+z_{k+l}|>|z_{k+1}+z_{k+l}|, \dots, |z_{k+l-1}+z_{k+l}|,
|z_{k+l}|, |z+z_{k+l}|>0$, $|z+z_{k+l}|>|z_{k+1}+z_{k+l}|, \dots, |z_{k+l-1}+z_{k+l}|>0$, 
$|z_{k}+z_{k+l}|>|z_{k+l}|>|z_{k}|>|z|>|z_{k+1}|
>\cdots>|z_{k+l-1}|>0$. But the form of the Laurent series (\ref{rationality-8}) is the same as the form of 
the expansion of (\ref{rationality-3}) as a Laurent series in $z_{i}+z_{k+l}$ for $i=1, \dots, k-1$,
$z_{k+j}$ for $j=0, 
\dots, l$, $z$ in the region $|z_{k+l}|>|z_{k}|>|z|>|z_{k+1}|>\cdots >|z_{k+l-1}|>0$, 
$|z_{1}+z_{k+l}|>\cdots>|z_{k-1}+z_{k+l}|>|z_{k+1}+z_{k+l}|, \dots, |z_{k+l-1}+z_{k+l}|,
|z_{k+l}|, |z+z_{k+l}|>0$, $|z+z_{k+l}|>|z_{k+1}+z_{k+l}|, \dots, |z_{k+l-1}+z_{k+l}|>0$.
Thus by Lemma \ref{it-series-conv}, 
(\ref{rationality-8}) is absolutely convergent to (\ref{rationality-3}) in this larger region.

Repeating this last step, we see that the series 
$$
\langle w_{2}', ((Y_{H}^{L}(v_{1}, z_{1})\cdots Y_{H}^{L}(v_{k}, z_{k})Y_{H}^{R}(\phi, z)
Y_{V}(v_{k+1}, z_{k+1})\cdots Y_{V}(v_{k+l-1}, z_{k+l-1})
v_{k+l})(w_{1}))(z_{k+l})\rangle
$$
is absolutely convergent to (\ref{rationality-3}) in the region $|z_{k+l}|>|z_{1}|>\cdots>|z_{k}|
>|z|>|z_{k+1}|>\cdots >|z_{k+l-1}|>0$.
\epfv

Since $H$ is $\Z$-graded, we have an operator $\mathbf{d}_{H}$ on $H$ defined by
$\mathbf{d}_{H}\phi=n\phi$ for $\phi\in H_{[n]}$.  We define $D_{H}$ on $H$ by 
$$((D_{H}\phi)(w_{1}))(z)=\frac{\partial}{\partial z}(\phi(w_{1}))(z)$$
for  $\phi\in H$, $w_{1}\in W_{1}$.

Though $H$ is $\Z$-graded, the grading is in general not lower bounded and hence 
it cannot be a $V$-bimodule. We now consider subspaces of $H$ which indeed have $V$-bimodule 
structures. 

For $N\in \Z$, let $H^{N}$ be the subspace of $H$ spanned by homogeneous elements, say $\phi$, 
satisfying the following condition:

\begin{enumerate}
\setcounter{enumi}{3}
\item The {\it $N$-weight-degree condition}: For $k, l\in \N$, $v_{1}, \dots, v_{k+l}\in V$, $w_{1}\in W_{1}$
and $w_{2}'\in W_{2}'$, expand the rational function
\begin{equation}\label{pole-order-1}
R(\langle w_{2}', Y_{W_{2}}(v_{1}, z_{1})\cdots Y_{W_{2}}(v_{k}, z_{k})
(\phi(Y_{W_{1}}(v_{k+1}, z_{k+1})\cdots Y_{W_{1}}(v_{k+l}, z_{k+l})w_{1}))(z)\rangle)
\end{equation}
in the region $|z_{k+l}|>|z_{1}-z_{k+l}|>\cdots >|z_{k}-z_{k+l}|>|z-z_{k+l}|>|z_{k+1}-z_{k+l}|>\cdots
>|z_{k+l-1}-z_{k+l}|>0$  as a Laurent series in $z_{i}-z_{k+l}$ for $i=1, \dots, k+l-1$ and $z-z_{k+l}$
with Laurent polynomials in $z_{k+l}$ as coefficients. Then the total degree of each monomial 
in $z_{i}-z_{k+l}$ for $i=1, \dots, k+l-1$ and $z-z_{k+l}$ (that is, the sum of the powers of 
$z_{i}-z_{k+l}$ for $i=1, \dots, k+l-1$ and $z-z_{k+l}$) in the expansion 
is larger than or equal to $N-\sum_{i=1}^{k+l}\wt v_{i}-\wt \phi$. 
\end{enumerate}

\begin{thm}\label{bimodule}
The $\Z$-graded space $H^{N}$, equipped with 
 the actions of the restrictions of $Y_{H}^{L}$, $Y_{H}^{R}$,
$\mathbf{d}_{H}$ and $D_{H}$ to $H^{N}$, is a $V$-bimodule.
\end{thm}
\pf
We need to prove that $H^{N}$ is closed under the actions of 
$Y_{H}^{L}$, $Y_{H}^{R}$ to $H^{N}$, $\mathbf{d}_{H}$ and $D_{H}$ and all axioms for $V$-bimodules hold.  
We prove them one by one as follows:

\paragraph{1.} $H^{N}$ is closed under the actions of 
$Y_{H}^{L}$, $Y_{H}^{R}$,  $H^{N}$, $\mathbf{d}_{H}$ and $D_{H}$:
For homogeneous $v\in V$, $\phi\in H^{N}$, we prove $Y_{H}^{L}(v, x)\phi\in H^{N}[[x, x^{-1}]]$.
We need only prove that $(Y_{H}^{L})_{p}(v)\phi$ satisfies the $N$-weight-degree condition. 
For $k, l\in \N$, homogeneous $v_{1}, \dots, v_{k+l}\in V$ , $w_{1}\in W_{1}$
and $w_{2}'\in W_{2}'$, by Proposition \ref{part-y-h-l-r}, 
the series (\ref{part-y-h-l-r-1}) is absolutely convergent 
in the region $|z_{1}|>\cdots |z_{k}|>|z|>|z_{k+1}|>\cdots >|z_{k+l}|$,
$|z|>|z_{0}|>0$, $|z_{i}|>|z+z_{0}|>0$ for $i=1, \dots, k$, $|z_{i}-z|>|z_{0}|>0$ for
$i=k+1, \dots, k+l$ to the rational function (\ref{part-y-h-l-r-2}). 
The rational function (\ref{part-y-h-l-r-2}) can be written as
\begin{align}\label{N-wt-dg-1}
&g(z_{1}, \dots, z_{k+l}, z_{0}, z)\cdot\nn
&\quad \cdot
z^{-\tilde{m}}
\prod_{i=1}^{k+l}z_{i}^{-\tilde{r}_{i}}(z_{0}+z)^{-n}\prod_{1\le i<j\le k+l}(z_{i}-z_{j})^{-p_{ij}}
\prod_{i=1}^{k+l}(z_{i}-z-z_{0})^{-q_{i}}\prod_{i=1}^{k+l}(z_{i}-z)^{-s_{i}}z_{0}^{-\tilde{t}},
\end{align}
where $g(z_{1}, \dots, z_{k+l}, z_{0}, z)$ is a polynomial 
in $z_{i}$ (for $i=1, \dots, k+l$), $z_{0}$, $z$ and 
$\tilde{m}, \tilde{r}_{i}, n, p_{ij}, q_{i}$, $s_{i}, \tilde{t}\in \N$. 
Since $g(z_{1}, \dots, z_{k+l}, z_{0}, z)$ is a polynomial, we can 
write (\ref{N-wt-dg-1}) as a linear combination of rational functions of the form
\begin{equation}\label{N-wt-dg-1.5}
z^{-m}
\prod_{i=1}^{k+l}z_{i}^{-r_{i}}(z_{0}+z)^{-n}\prod_{1\le i<j\le k+l}(z_{i}-z_{j})^{-p_{ij}}
\prod_{i=1}^{k+l}(z_{i}-z-z_{0})^{-q_{i}}\prod_{i=1}^{k+l}(z_{i}-z)^{-s_{i}}z_{0}^{-t},
\end{equation}
where $m, r_{i}, t\in \Z$. 
We can write (\ref{N-wt-dg-1.5}) as
\begin{align}\label{N-wt-dg-2}
&z_{k+l}^{-m}\left(1+\frac{z-z_{k+l}}{z_{k+l}}\right)^{-m}
\left(\prod_{i=1}^{k+l-1}z_{k+l}^{-r_{i}}\left(1+\frac{z_{i}-z_{k+l}}{z_{k+l}}\right)^{-r_{i}}\right)
z_{k+l}^{-r_{k+l}}\cdot\nn
&\quad\cdot z_{k+l}^{-n}\left(1+\frac{z_{0}+z-z_{k+l}}{z_{k+l}}\right)^{-n}
\prod_{1\le i<j\le k+l-1}(z_{i}-z_{k+l})^{-p_{ij}}
\left(1-\frac{z_{j}-z_{k+l}}{z_{i}-z_{k+l}}\right)^{-p_{ij}}\cdot\nn
&\quad \cdot \prod_{i=1}^{k+l-1}(z_{i}-z_{k+l})^{-p_{i, k+l}}
\left(\prod_{i=1}^{k}(z_{i}-z_{k+l})^{-q_{i}}
\left(1-\frac{z+z_{0}-z_{k+l}}{z_{i}-z_{k+l}}\right)^{-q_{i}}\right)\cdot\nn
&\quad \cdot
\left(\prod_{i=k+1}^{k+l-1}(-(z+z_{0}-z_{k+l}))^{-q_{i}}
\left(1-\frac{z_{i}-z_{k+l}}{z+z_{0}-z_{k+l}}\right)^{-q_{i}}\right)\cdot\nn
&\quad \cdot (-(z+z_{0}-z_{k+l}))^{-q_{k+l}}\cdot\nn
&\quad\cdot  \prod_{i=1}^{k}(z_{i}-z_{k+l})^{-s_{i}}\left(1-\frac{z-z_{k+l}}{z_{i}-z_{k+l}}\right)^{-s_{i}}
\prod_{i=k+1}^{k+l-1}(-(z-z_{k+l}))^{-s_{i}}\left(1-\frac{z_{i}-z_{k+l}}{z-z_{k+l}}\right)^{-s_{i}}\cdot\nn
&\quad\cdot 
(-(z-z_{k+l}))^{-s_{k+l}}(z_{0}+z-z_{k+l})^{-t}\left(1-\frac{z-z_{k+l}}{z_{0}+z-z_{k+l}}\right)^{-t}.
\end{align}
When (\ref{N-wt-dg-2})) is expanded
in the region $|z_{k+l}|>|z_{1}-z_{k+l}|>\cdots >|z_{k}-z_{k+l}|>|z_{0}+z-z_{k+l}|>|z-z_{k+l}|
>|z_{k+1}-z_{k+l}|>\cdots
>|z_{k+l-1}-z_{k+l}|>0$  as a Laurent series 
in $z_{i}-z_{k+l}$ for $i=1, \dots, k+l-1$, $z-z_{k+l}$ and $z_{0}-z-z_{k+l}$ with 
Laurent polynomials in $z_{k+l}$ as coefficients, the lowest total degree is 
$$-\sum_{1\le i< j\le k+l}p_{ij}-\sum_{i=1}^{k+l}q_{i}-\sum_{i=1}^{k+l}s_{i}-t.$$
Since $\phi$ satisfies 
the $N$-weight-degree condition, 
we have
\begin{align}\label{N-wt-dg-3}
-\sum_{1\le i< j\le k+l}p_{ij}-\sum_{i=1}^{k+l}q_{i}-\sum_{i=1}^{k+l}s_{i}-t
\ge N-\sum_{i=1}^{k+l}\wt v_{i}-\wt v-\wt \phi.
\end{align}

We can also write  (\ref{N-wt-dg-1.5}) as
\begin{align}\label{N-wt-dg-4}
&z_{k+l}^{-m}\left(1+\frac{z-z_{k+l}}{z_{k+l}}\right)^{-m}
\left(\prod_{i=1}^{k+l-1}z_{k+l}^{-r_{i}}\left(1+\frac{z_{i}-z_{k+l}}{z_{k+l}}\right)^{-r_{i}}\right)
z_{k+l}^{-r_{k+l}}\cdot\nn
&\quad\cdot z_{k+l}^{-n}\left(1+\frac{z-z_{k+l}}{z_{k+l}}\right)^{-n}
\left(1+\frac{z_{0}}{z_{k+l}\left(1+\frac{z-z_{k+l}}{z_{k+l}}\right)}\right)^{-n}\cdot\nn
&\quad \cdot
\prod_{1\le i<j\le k+l-1}(z_{i}-z_{k+l})^{-p_{ij}}
\left(1-\frac{z_{j}-z_{k+l}}{z_{i}-z_{k+l}}\right)^{-p_{ij}}\prod_{i=1}^{k+l-1}(z_{i}-z_{k+l})^{-p_{i, k+l}}\cdot\nn
&\quad \cdot 
\prod_{i=1}^{k+l-1}(z_{i}-z_{k+l})^{-q_{i}}\left(1-\frac{z-z_{k+l}}{z_{i}-z_{k+l}}\right)^{-q_{i}}
\left(1-\frac{z_{0}}{(z_{i}-z_{k+l})\left(1-\frac{z-z_{k+l}}{z_{i}-z_{k+l}}\right)}\right)^{-q_{i}}
\cdot\nn
&\quad \cdot (-(z-z_{k+l}))^{-q_{k+l}}\left(1+\frac{z_{0}}{z-z_{k+l}}\right)^{-q_{k+l}}\cdot\nn
&\quad \cdot  \prod_{i=1}^{k}(z_{i}-z_{k+l})^{-s_{i}}\left(1-\frac{z-z_{k+l}}{z_{i}-z_{k+l}}\right)^{-s_{i}}
\prod_{i=k+1}^{k+l-1}(-(z-z_{k+l}))^{-s_{i}}\left(1-\frac{z_{i}-z_{k+l}}{z-z_{k+l}}\right)^{-s_{i}}\cdot\nn
&\quad \cdot 
(-(z-z_{k+l}))^{-s_{k+l}}z_{0}^{-t},
\end{align}
Now expand (\ref{N-wt-dg-4}) first as a 
Laurent series in $z_{0}$ in the region $|z|>|z_{0}|>0$, $|z_{i}-z|>|z_{0}|>0$ and then expand 
the coefficient of the $-p-1$ power of $z_{0}$ 
in the region $|z_{k+l}|>|z_{1}-z_{k+l}|>\cdots >|z_{k}-z_{k+l}|>|z-z_{k+l}|
>|z_{k+1}-z_{k+l}|>\cdots
>|z_{k+l-1}-z_{k+l}|>0$  as a Laurent series in $z_{i}-z_{k+l}$ for $i=1, \dots, k+l-1$ and $z-z_{k+l}$
with Laurent polynomials in $z_{k+l}$ as coefficients. The $-p-1$ powers of $z_{0}$ in the expansion of
(\ref{N-wt-dg-4}) comes from the expansion of 
\begin{align}\label{N-wt-dg-5}
& \left(1+\frac{z_{0}}{z_{k+l}\left(1+\frac{z-z_{k+l}}{z_{k+l}}\right)}\right)^{-n}
\prod_{i=1}^{k+l-1}\left(1-\frac{z_{0}}{(z_{i}-z_{k+l})\left(1-\frac{z-z_{k+l}}{z_{i}-z_{k+l}}\right)}\right)^{-q_{i}}
\cdot \nn
&\quad \cdot 
\left(1+\frac{z_{0}}{z-z_{k+l}}\right)^{-q_{k+l}}z_{0}^{-t}.
\end{align}
The expansion of (\ref{N-wt-dg-5})  is an infinite linear combination 
of the monomial of $z_{0}$ of the form
\begin{equation}\label{N-wt-dg-6}
\left(\frac{z_{0}}{z_{k+l}\left(1+\frac{z-z_{k+l}}{z_{k+l}}\right)}\right)^{a}
\prod_{i=1}^{k+l-1}
\left(\frac{z_{0}}{(z_{i}-z_{k+l})\left(1-\frac{z-z_{k+l}}{z_{i}-z_{k+l}}\right)}\right)^{b_{i}}
\left(\frac{z_{0}}{z-z_{k+l}}\right)^{c}z_{0}^{-t}
\end{equation}
for $a, b_{i}, c\in \N$. In the case that the power of $z_{0}$ is $-p-1$, we have 
\begin{equation}\label{N-wt-dg-7}
a+\sum_{i=1}^{k+l-1}b_{i}+c-t=-p-1.
\end{equation}
In this case, the 
total degree of (\ref{N-wt-dg-6}) in $z_{i}-z_{k+l}$ for $i=1, \dots, k+l$ 
and $z-z_{k+l}$ is $-\sum_{i=1}^{k+l-1}b_{i}-c$.
From  (\ref{N-wt-dg-4}) , (\ref{N-wt-dg-5})  and (\ref{N-wt-dg-6}) and this fact, 
we see that the total degrees of each monomial 
in the expansion of the coefficient of the $-p-1$ power of $z_{0}$ is larger than or equal to 
\begin{equation}\label{N-wt-dg-8}
-\sum_{1\le i< j\le k+l}p_{ij}-\sum_{i=1}^{k+l}q_{i}-\sum_{i=1}^{k+l-1}b_{i}-c-\sum_{i=1}^{k+l}s_{i}.
\end{equation}
Using  (\ref{N-wt-dg-7}), we see that (\ref{N-wt-dg-8}) is equal to
\begin{equation}\label{N-wt-dg-9}
-\sum_{1\le i< j\le k+l}p_{ij}-\sum_{i=1}^{k+l}q_{i}+a-t+p+1-\sum_{i=1}^{k+l}s_{i}.
\end{equation}
Using $a\in \N$ and (\ref{N-wt-dg-3}), we see that (\ref{N-wt-dg-9}) is larger than or equal to 
\begin{align*}
&-\sum_{1\le i< j\le k+l}p_{ij}-\sum_{i=1}^{k+l}q_{i}-t+p+1-\sum_{i=1}^{k+l}s_{i}\nn
&\quad \ge N-\sum_{i=1}^{k+l}\wt v_{i}-\wt v-\wt \phi+p+1\nn
&\quad =N-\sum_{i=1}^{k+l}\wt v_{i}-\wt (Y_{H}^{L})_{p}(v)\phi,
\end{align*}
where in the last step we have used the formula
$$\wt (Y_{H}^{L})_{p}(v)\phi=\wt v+\wt \phi-p-1,$$
which in turn follows from the $\mathbf{d}$-conjugation property of $Y_{H}^{L}$
that we shall prove below. Since the proof the $\mathbf{d}$-conjugation property of $Y_{H}^{L}$
does not need the result that $H^{N}$ is closed under the actions of 
$Y_{H}^{L}$ and its proof, we can indeed use this formula here. This 
proves that $(Y_{H}^{L})_{p}(v)\phi$ satisfies the $N$-weight-degree condition.

We omit the proof that $Y_{H}^{R}(\phi, x)v\in H^{N}[[x, x^{-1}]]$ for
$v\in V$ and $\phi\in H^{N}$. It is similar to the proof above for $Y_{H}^{L}$.

Since for homogeneous $\phi\in H^{N}$, 
$\mathbf{d}_{H}\phi=(\wt \phi)\phi$, $\mathbf{d}_{H}\phi$ 
also satisfies the $N$-weight-degree condition. So 
$\mathbf{d}_{H}$ maps $H^{N}$ to $H^{N}$. 

For homogeneous 
$\phi\in H^{N}$, $k, l\in \N$, homogeneous $v_{1}, \dots, v_{k+l}\in V$ , $w_{1}\in W_{1}$
and $w_{2}'\in W_{2}'$,
by definition, we have 
\begin{align}\label{closed-1}
R&(\langle w_{2}', Y_{W_{2}}(v_{1}, z_{1})\cdots Y_{W_{2}}(v_{k}, z_{k})
((D_{H}\phi)(Y_{W_{1}}(v_{k+1}, z_{k+1})\cdots Y_{W_{1}}(v_{k+l}, z_{k+l})w_{1}))(z)\rangle)\nn
&=\frac{\partial}{\partial z}R(\langle w_{2}', Y_{W_{2}}(v_{1}, z_{1})\cdots Y_{W_{2}}(v_{k}, z_{k})
(\phi(Y_{W_{1}}(v_{k+1}, z_{k+1})\cdots Y_{W_{1}}(v_{k+l}, z_{k+l})w_{1}))(z)\rangle).
\end{align}
Then the expansion of (\ref{closed-1}) in the region $|z_{k+l}|>
|z_{i}-z_{k+l}>\cdots >|z_{k}-z_{k+l}|>|z-z_{k+l}|>|z_{k+1}|>\cdots
>|z_{k+l-1}-z_{k+l}|>0$  as a Laurent series in $z_{i}-z_{k+l}$ for $i=1, \dots, k+l-1$ and $z-z_{k+l}$
with Laurent polynomials in $z_{k+l}$ as coefficients is equal to the derivative with respect to $z$ of the
same expansion of (\ref{pole-order-1}). In particular, the total degree of each monomial 
in $z_{i}-z_{k+l}$ for $i=1, \dots, k+l-1$ and $z-z_{k+l}$ in the expansion of 
(\ref{closed-1})  is equal to the total degree of a monomial in the expansion 
of (\ref{pole-order-1}) minus $1$. Thus the total degree of each monomial 
in $z_{i}-z_{k+l}$ for $i=1, \dots, k+l-1$ and $z-z_{k+l}$ in the expansion of 
(\ref{closed-1}) is larger than or equal to 
$N-\sum_{i=1}^{k+l}\wt v_{i}-\wt \phi-1$.
But from the definition of $D_{H}$ and $\wt \phi$, we have
\begin{align*}
a^{\mathbf{d}_{W_{2}}}((D_{H}\phi)(w_{1}))(z)
&=a^{\mathbf{d}_{W_{2}}}\frac{\partial}{\partial z}(\phi(w_{1}))(z)\nn
&=a^{\swt \phi}\frac{\partial}{\partial z}(\phi(a^{\mathbf{d}_{W_{1}}}w_{1}))(az)\nn
&=a^{\swt \phi+1}\frac{\partial}{\partial z'}(\phi(a^{\mathbf{d}_{W_{1}}}w_{1}))(z')\lbar_{z'=az}\nn
&=a^{\swt \phi+1}((D_{H}\phi)(a^{\mathbf{d}_{W_{1}}}w_{1}))(az),
\end{align*}
proving $\wt \phi+1=\wt D_{H}\phi$. 
Therefore the total degree of each monomial 
in $z_{i}-z_{k+l}$ for $i=1, \dots, k+l-1$ and $z-z_{k+l}$ in the expansion of 
(\ref{closed-1}) is larger than or equal to 
$N-\sum_{i=1}^{k+l}\wt v_{i}-\wt (D_{H}\phi).$
So $D_{H}\phi$ satisfies the $N$-weight-degree condition and is in $H^{N}$.

\paragraph{2.} Axioms for the grading: 
Take $v_{1}=\cdots=v_{k+l}=\one$ in (\ref{pole-order-1}). Then (\ref{pole-order-1}) becomes
$\langle w_{2}', (\phi(w_{1}))(z)\rangle$. Expand this as a power series in $z-z_{k+l}$ in the region 
$|z_{k+l}|>|z-z_{k+l}|$. Then the lower degree term in this expansion is the constant term 
$\langle w_{2}', (\phi(w_{1}))(z_{k})\rangle$. Thus the lowest of the degrees of the monomials in 
$z_{i}-z_{k+l}$ for $i=1, \dots, k+l-1$ and $z-z_{k+l}$ is $0$. So we obtain 
$0\ge N-\sum_{i=1}^{k+l}v_{i}+\wt \phi=N-\wt \phi$ or $\wt \phi\ge N$. This proves
the lower bound condition of $H^{N}$.

We now prove the $\mathbf{d}$-conjugation property of $Y_{H}^{L}$ which is equivalent to the 
$\mathbf{d}$-commutator formula for $Y_{H}^{L}$. 
Let $v\in V$ and $\phi\in H^{N}$.
Then for $a\in \R$, $w_{2}'\in W_{2}'$ and $w_{1}\in W_{1}$, we have 
\begin{eqnarray*}
\lefteqn{\langle w_{2}', ((a^{\mathbf{d}_{H}}Y_{H}^{L}(v, z_{1})\phi)(w_{1}))(z_{2})\rangle}\nn
&&=\langle w_{2}', a^{\mathbf{d}_{W_{2}}}((Y_{H}^{L}(v, z_{1})\phi)(a^{-\mathbf{d}_{W_{1}}}w_{1}))(a^{-1}z_{2})\rangle\nn
&&=\langle w_{2}', a^{\mathbf{d}_{W_{2}}}Y_{W_{2}}(v, z_{1}+a^{-1}z_{2})(\phi(a^{-\mathbf{d}_{W_{1}}}w_{1}))(a^{-1}z_{2})\rangle\nn
&&=\langle w_{2}', Y_{W_{2}}(a^{\mathbf{d}_{V}}v, az_{1}+z_{2})a^{\mathbf{d}_{W_{2}}}(\phi(a^{-\mathbf{d}_{W_{1}}}w_{1}))(a^{-1}z_{2})\rangle\nn
&&=\langle w_{2}', Y_{W_{2}}(a^{\mathbf{d}_{V}}v, az_{1}+z_{2})((a^{H}\phi)w_{1})(z_{2})\rangle\nn
&&=\langle w_{2}', ((Y_{H}(a^{\mathbf{d}_{V}}v, az_{1})a^{\mathbf{d}_{H}}\phi)(w_{1}))(z_{2})\rangle.
\end{eqnarray*}
Thus we have 
$$a^{\mathbf{d}_{H}}Y_{H}^{L}(v, z_{1})\phi=Y_{H}(a^{\mathbf{d}_{V}}v, az_{1})a^{\mathbf{d}_{H}}\phi,$$
proving the $\mathbf{d}$-conjugation property of 
$Y_{H}^{L}$. Similarly we can prove the $\mathbf{d}$-conjugation property of 
$Y_{H}^{R}$ which is equivalent to the 
$\mathbf{d}$-commutator formula for $Y_{H}^{R}$. We omit the proof.

\paragraph{3.} The identity property of $Y_{H}^{L}$ and the creation property
of $Y_{H}^{R}$: These properties follows directly from  (\ref{y-l-rigor}) 
and  (\ref{y-r-rigor}).

\paragraph{4.} The rationality and the pole-order condition: 
Let $w_{2}'\in W_{2}'$, $w_{1}\in W_{1}$, $v_{1}, \dots, v_{k+l}\in V$ and 
$\phi\in H^{N}$ be homogeneous. Since (\ref{y-h-l-r-3}) is of the form (\ref{rationality-3}), 
\begin{align}\label{rationality-9}
&z^{m}\prod_{1\le i<j\le k+l-1}(z_{i}-z_{j})^{p_{ij}}\prod_{i=1}^{k+l-1}z_{i}^{r_{i}}
\prod_{i=1}^{k+l-1}(z_{i}-z)^{s_{i}}\cdot \nn
&\cdot \langle w_{2}', ((Y_{H}^{L}(v_{1}, z_{1})\cdots Y_{H}^{L}(v_{k}, z_{k})Y_{H}^{R}(\phi, z)
Y_{V}(v_{k+1}, z_{k+1})\cdots Y_{V}(v_{k+l-1}, z_{k+l-1})
v_{k+l})(w_{1}))(z_{k+l})\rangle
\end{align}
must be convergent absolutely in the region $|z_{k+l}|>|z_{1}|>\cdots>|z_{k}|
>|z|>|z_{k+1}|>\cdots >|z_{k+l-1}|>0$ to
\begin{equation}\label{rationality-10}
g(z_{1}, \dots, z_{k+l}, z)(z+z_{k+l})^{-n}\prod_{i=1}^{k+l-1}(z_{i}+z_{k+l})^{-p_{i}}
 z_{k+l}^{-r_{k+l}}.
\end{equation}
So (\ref{rationality-9}) is the expansion of the rational function (\ref{rationality-10})
obtained by expanding the negative powers of $z+z_{k+l}$ and $z_{i}+z_{k+l}$ for $i=1, \dots, k+l-1$
as power series of $z$ and $z_{i}$, respectively. Thus (\ref{rationality-9})
must be a power series in $z$ and $z_{i}$ for $i=1, \dots, k+l-1$ and must be absolutely convergent 
to (\ref{rationality-10}) in the larger region $|z_{k+l}|>|z|, |z_{1}|, \dots, |z_{k+l-1}|>0$ than 
$|z_{k+l}|>|z_{1}|>\cdots>|z_{k}|
>|z|>|z_{k+1}|>\cdots >|z_{k+l-1}|>0$.

Since (\ref{rationality-9}) is a power series in $z$ and $z_{i}$ for $i=1, \dots, k+l-1$ 
for all  $w_{2}'\in W_{2}$ and $w_{1}\in W_{1}$ and all $z_{k+l}$ satisfying 
$|z_{k+l}|>|z|, |z_{1}|, \dots, |z_{k+l-1}|>0$ and 
$m$, $p_{ij}$ (for $1\le i<j\le k+l-1$), $r_{i}$ (for $i=1, \dots, k+l-1$), 
$s_{i}$ (for $i=1, \dots, k+l-1$) are independent of $w_{2}'\in W_{2}$ and $w_{1}\in W_{1}$, 
\begin{align}\label{rationality-11}
&z^{m}\prod_{1\le i<j\le k+l-1}(z_{i}-z_{j})^{p_{ij}}\prod_{i=1}^{k+l-1}z_{i}^{r_{i}}
\prod_{i=1}^{k+l-1}(z_{i}-z)^{s_{i}}\cdot \nn
&\cdot Y_{H}^{L}(v_{1}, z_{1})\cdots Y_{H}^{L}(v_{k}, z_{k})Y_{H}^{R}(\phi, z)
Y_{V}(v_{k+1}, z_{k+1})\cdots Y_{V}(v_{k+l-1}, z_{k+l-1})
v_{k+l}
\end{align}
is also a power series in $z$ and $z_{i}$ for $i=1, \dots, k+l-1$ with coefficients in $H^{N}$.
Let $\phi'\in (H^{N})'$. Since the grading of $(H^{N})'$ 
is bounded from below and the $\mathbf{d}$-commutator formula 
for $Y_{H}^{L}$ holds, 
\begin{align}\label{rationality-11.5}
\langle \phi', Y_{H}^{L}(v_{1}, z_{1})\cdots Y_{H}^{L}(v_{k}, z_{k})Y_{H}^{R}(\phi, z)
Y_{V}(v_{k+1}, z_{k+1})\cdots Y_{V}(v_{k+l-1}, z_{k+l-1})
v_{k+l}\rangle
\end{align}
has only finitely many positive power terms in $z_{1}$. Thus 
\begin{align}\label{rationality-12}
&z^{m}\prod_{1\le i<j\le k+l-1}(z_{i}-z_{j})^{p_{ij}}\prod_{i=1}^{k+l-1}z_{i}^{r_{i}}
\prod_{i=1}^{k+l-1}(z_{i}-z)^{s_{i}}\cdot \nn
&\cdot \langle \phi',  Y_{H}^{L}(v_{1}, z_{1})\cdots Y_{H}^{L}(v_{k}, z_{k})Y_{H}^{R}(\phi, z)
Y_{V}(v_{k+1}, z_{k+1})\cdots Y_{V}(v_{k+l-1}, z_{k+l-1})
v_{k+l}\rangle
\end{align}
also has only  finitely many positive power terms in $z_{1}$. So (\ref{rationality-12}) is a polynomial 
in $z_{1}$. 

Take the coefficient of a fixed nonnegative power of $z_{1}$ in (\ref{rationality-12}). Since 
$$z^{m}\prod_{1\le i<j\le k+l-1}(z_{i}-z_{j})^{p_{ij}}\prod_{i=1}^{k+l-1}z_{i}^{r_{i}}
\prod_{i=1}^{k+l-1}(z_{i}-z)^{s_{i}}$$
is a polynomial in $z_{1}$, this coefficient involves only finitely  many
terms in the Laurent series $Y_{H}^{L}(v_{1}, z_{1})$. Then this coefficient must be 
of the form
\begin{align}\label{rationality-13}
\Biggl\langle \phi', & \left(\sum_{n=n_{1}}^{n_{2}}
f_{n}(z_{2}, \dots, z_{k+l-1}, z)(Y_{H}^{L})_{n}(v_{1})\right)
Y_{H}^{L}(v_{2}, z_{2})\cdots Y_{H}^{L}(v_{k}, z_{k})\cdot\nn
&\quad\quad \quad\quad\cdot Y_{H}^{R}(\phi, z)
Y_{V}(v_{k+1}, z_{k+1})\cdots Y_{V}(v_{k+l-1}, z_{k+l-1})
v_{k+l}\Biggr\rangle,
\end{align}
where $n_{1}\le n_{2}$ are integers and $f_{n}(z_{2}, \dots, z_{k+l-1}, z)$ for $n=n_{1}, \dots, n_{2}$
are polynomials in $z_{2}, \dots, z_{k+l-1}, z$. Since  the grading of $(H^{N})'$ 
is bounded from below, the $\mathbf{d}$-commutator formula 
for $Y_{H}^{L}$ holds and $f_{n}(z_{2}, \dots, z_{k+l-1}, z)$ for $n=n_{1}, \dots, n_{2}$
are polynomials in $z_{2}, \dots, z_{k+l-1}, z$, (\ref{rationality-13}) has 
only  finitely many positive power terms in $z_{2}$. Since (\ref{rationality-13})  
as a coefficient of a fixed nonnegative power of $z_{1}$ in the power series (\ref{rationality-12})
must be a power series in $z_{2}, \dots, z_{k+l-1}, z$, it must be a polynomial 
in $z_{2}$. 

Since the coefficient of every fixed nonnegative power of $z_{1}$ in (\ref{rationality-12})
is a polynomial in $z_{2}$ and since (\ref{rationality-12}) is a polynomial in $z_{1}$, 
(\ref{rationality-12}) is also a polynomial in $z_{2}$.

Repeating these steps, we see that (\ref{rationality-12}) is in fact a polynomial $h(z_{1}, \dots, z_{k+l-1}, z)$ in 
$z$ and $z_{i}$ for $i=1, \dots, k+l-1$. Note that the form of the series (\ref{rationality-11.5})
is the same as the form of the Laurent series expansion in the region 
$|z_{1}|>\cdots>|z_{k}|>|z|>z_{k+1}|>\cdots>|z_{k+l-1}|>0$ of a rational function 
in $z$ and $z_{i}$ for $i=1, \dots, k+l-1$ with the only possible poles $z, z_{i}=0$ (for $i=1, \dots, k+l-1$),
$z_{i}=z_{j}$ (for $1\le i<j\le k+l-1$) and $z=z_{i}$ (for $i=1, \dots, k+l-1$). 
Thus (\ref{rationality-11.5}) must be absolutely convergent in the region
$|z_{1}|>\cdots>|z_{k}|>|z|>z_{k+1}|>\cdots>|z_{k+l-1}|>0$ to 
$$\frac{h(z_{1}, \dots, z_{k+l-1}, z)}{\displaystyle z^{m}\prod_{1\le i<j\le k+l-1}(z_{i}-z_{j})^{p_{ij}}\prod_{i=1}^{k+l-1}z_{i}^{r_{i}}
\prod_{i=1}^{k+l-1}(z_{i}-z)^{s_{i}}}.$$
Thus the rationality is proved. 

The pole-order condition follows immediately from the facts that (\ref{rationality-11})
is an element of $H^{N}[[z_{1}, \dots, z_{k+l-1}, z]]$ and that
$m, p_{ij}, r_{i}, s_{i}\in \N$ depend only on the pairs $(v_{k+l}, \phi)$, 
$(v_{i}, v_{j})$, $(v_{i}, v_{k+l})$, $(v_{i}, \phi)$, respectively. 

\paragraph{5.} The associativity:  We prove only the associativity of $Y_{H}^{L}$. 
The associativity of $Y_{H}^{R}$ and the associativity of $Y_{H}^{L}$ and $Y_{H}^{R}$ can be proved 
similarly and are omitted.
Take $k=2$, $l=1$, $v_{3}=\one$ and $z=0$ in (\ref{rationality-1}). Then we have proved above that 
for $w_{2}'\in W_{2}'$, $v_{1}, v_{2}\in V$, $w_{1}\in W_{1}$, $\phi\in H^{N}$
$$\langle w_{2}',  Y_{W_{2}}(v_{1}, z_{1}+z_{3})
Y_{W_{2}}(v_{2}, z_{2}+z_{3})(\phi(w_{1}))(z_{3})\rangle$$
is absolutely convergent in the region $|z_{1}+z_{2}|>|z_{2}+z_{3}|>|z_{3}|>0$ to a rational function 
\begin{equation}\label{y-l-ass-1}
R(\langle w_{2}',  Y_{W_{2}}(v_{1}, z_{1}+z_{3})
Y_{W_{2}}(v_{2}, z_{2}+z_{3})(\phi(w_{1}))(z_{3})\rangle)
\end{equation}
in $z_{1}, z_{2}, z_{3}$ with the only possible poles $z_{1}+z_{3}=0$, $z_{2}+z_{3}=0$, $z_{1}, z_{2},
z_{3}=0$
and $z_{1}=z_{2}$. Moreover, we have proved above that 
$$\langle w_{2}',   ((Y_{H}^{L}(v_{1}, z_{1})Y_{H}^{L}(v_{2}, z_{2})\phi)(w_{1}))(z_{3})\rangle$$
is also absolutely convergent in the region $|z_{3}|>|z_{1}|>|z_{2}|>0$ to (\ref{y-l-ass-1}). In particular, 
we have 
\begin{align}\label{y-l-ass-2}
R&(\langle w_{2}',  Y_{W_{2}}(v_{1}, z_{1}+z_{3})
Y_{W_{2}}(v_{2}, z_{2}+z_{3})(\phi(w_{1}))(z_{3})\rangle)\nn
&=R(\langle w_{2}',   ((Y_{H}^{L}(v_{1}, z_{1})Y_{H}^{L}(v_{2}, z_{2})\phi)(w_{1}))(z_{3})\rangle).
\end{align}

On the other hand, by the associativity of $Y_{W_{2}}$, we have 
\begin{align}\label{y-l-ass-3}
R&(\langle w_{2}',  Y_{W_{2}}(v_{1}, z_{1}+z_{3})
Y_{W_{2}}(v_{2}, z_{2}+z_{3})(\phi(w_{1}))(z_{3})\rangle)\nn
&=R(\langle w_{2}',  Y_{W_{2}}(Y_{V}(v_{1}, z_{1}-z_{2})v_{2}, z_{2}+z_{3})(\phi(w_{1}))(z_{3})\rangle).
\end{align}
Using the homogeneous basis $\{e^{V}_{(q; \mu)}\}_{q\in \Z,\; \mu\in M_{q}}$ of $V$
and the subset $\{(e^{V}_{(q; \mu)})'\}_{q\in \Z,\; \mu\in M_{q}}$ of $V'$ as above and then 
using (\ref{y-l-rigor-r}), the right-hand side of (\ref{y-l-ass-3}) can be expanded as 
the series
\begin{align}\label{y-l-ass-4}
\sum_{q\in \Z}&\Biggl(\sum_{\mu\in M_{q}} 
R(\langle w_{2}',  Y_{W_{2}}(e^{V}_{(q; \mu)}, z_{2}+z_{3})(\phi(w_{1}))(z_{3})\rangle)
\langle (e^{V}_{(q; \mu)})', Y_{V}(v_{1}, z_{1}-z_{2})v_{2}\rangle\Biggr)\nn
&=\sum_{q\in \Z}\Biggl(\sum_{\mu\in M_{q}} 
R(\langle w_{2}',  ((Y_{H}^{L}(e^{V}_{(q; \mu)}, z_{2})\phi)(w_{1}))(z_{3})\rangle
\langle (e^{V}_{(q; \mu)})', Y_{V}(v_{1}, z_{1}-z_{2})v_{2}\rangle\Biggr).
\end{align}
in the region $|z_{2}+z_{3}|>|z_{1}-z_{2}|>0$. But the right-hand side of (\ref{y-l-ass-4})
is absolutely convergent in the region $|z_{2}|>|z_{1}-z_{2}|>0$ to the rational function 
\begin{equation}\label{y-l-ass-5}
R(\langle w_{2}',  ((Y_{H}^{L}(Y_{V}(v_{1}, z_{1}-z_{2})v_{2}, z_{2})\phi)(w_{1}))(z_{3})\rangle.
\end{equation}
Using the calculation from  (\ref{y-l-ass-2})--(\ref{y-l-ass-5}), we obtain
\begin{align}\label{y-l-ass-6}
R&(\langle w_{2}',   ((Y_{H}^{L}(v_{1}, z_{1})Y_{H}^{L}(v_{2}, z_{2})\phi)(w_{1}))(z_{3})\rangle)\nn
&=R(\langle w_{2}',  ((Y_{H}^{L}(Y_{V}(v_{1}, z_{1}-z_{2})v_{2}, z_{2})\phi)(w_{1}))(z_{3})\rangle).
\end{align}

Taking $k=2$, $l=1$, $v_{3}=\one$ and $z=0$ in (\ref{rationality-3}), we see that for
$\tilde{p}_{12}\ge p_{12}$, $\tilde{r}_{1}\ge r_{1}$ and $\tilde{r}_{2}\ge r_{2}$
\begin{equation}\label{y-l-ass-6.5}
(z_{1}-z_{2})^{\tilde{p}_{12}}z_{1}^{\tilde{r}_{1}}z_{2}^{\tilde{r}_{2}}
\langle w_{2}',   ((Y_{H}^{L}(v_{1}, z_{1})Y_{H}^{L}(v_{2}, z_{2})\phi)(w_{1}))(z_{3})\rangle
\end{equation}
is absolutely convergent in the region $|z_{3}|>|z_{1}|, |z_{2}|>0$ to 
the rational function
\begin{equation}\label{y-l-ass-7}
\frac{\tilde{g}(z_{1}, z_{2}, z_{3})}{z_{3}^{n+r_{3}}(z_{1}+z_{3})^{p_{1}}(z_{2}+z_{3})^{p_{2}}},
\end{equation}
where $\tilde{g}(z_{1}, z_{2}, z_{3})$ is a polynomial in $z_{1}, z_{2}, z_{3}$. 
From (\ref{y-l-ass-6}), we also see that
\begin{equation}\label{y-l-ass-8}
(z_{1}-z_{2})^{\tilde{p}_{12}}z_{1}^{\tilde{r}_{1}}z_{2}^{\tilde{r}_{2}}
\langle w_{2}',  ((Y_{H}^{L}(Y_{V}(v_{1}, z_{1}-z_{2})v_{2}, z_{2})\phi)(w_{1}))(z_{3})\rangle
\end{equation}
must be absolutely convergent in the region $|z_{3}|>|z_{2}|>|z_{1}-z_{2}|>0$ to 
(\ref{y-l-ass-7}). Since there is no negative power term in (\ref{y-l-ass-7}), (\ref{y-l-ass-8})
is in fact absolutely convergent in the larger region $|z_{3}|>|z_{1}=z_{2}+(z_{1}-z_{2})|, |z_{2}|>0$
than $|z_{3}|>|z_{2}|>|z_{1}-z_{2}|>0$. 

We have proved that (\ref{y-l-ass-6.5}) is the expansion of
(\ref{y-l-ass-7}) by expanding the negative powers of $z_{1}+z_{3}$ and $z_{2}+z_{3}$ 
as power series of $z_{1}$ and $z_{2}$, respectively. We have also proved that
(\ref{y-l-ass-8}) is the expansion of
(\ref{y-l-ass-7}) by expanding the negative powers of $z_{1}+z_{3}$ and $z_{2}+z_{3}$ 
as power series of $z_{1}$ and $z_{2}$, respectively, and then expanding positive powers of 
$z_{1}=z_{2}+(z_{1}-z_{2})$ using the binomial expansion as polynomials in $z_{2}$ and $z_{1}-z_{2}$. 
Thus as a power series in $z_{2}$ and $z_{1}-z_{2}$, 
(\ref{y-l-ass-8}) can be obtained by expanding the positive powers of $z_{1}=z_{2}+(z_{1}-z_{2})$ 
in the series (\ref{y-l-ass-6.5})
using the binomial expansion as polynomials in $z_{2}$ and $z_{1}-z_{2}$. 
By the composability, $p_{12}$, $r_{1}$ and $r_{2}$  are
independent of $w_{2}'$, $w_{1}$ and $z_{3}$. Hence 
\begin{equation}\label{y-l-ass-9}
(z_{1}-z_{2})^{p_{12}}z_{1}^{r_{1}}z_{2}^{r_{2}}
Y_{H}^{L}(v_{1}, z_{1})Y_{H}^{L}(v_{2}, z_{2})\phi
\end{equation}
is a power series in $z_{1}$ and $z_{2}$ and 
\begin{equation}\label{y-l-ass-10}
(z_{1}-z_{2})^{p_{12}}z_{1}^{r_{1}}z_{2}^{r_{2}}
Y_{H}^{L}(Y_{V}(v_{1}, z_{1}-z_{2})v_{2}, z_{2})\phi
\end{equation}
is a power series in $z_{2}$ and $z_{1}-z_{2}$. Moreover, (\ref{y-l-ass-10})
can also be obtained from (\ref{y-l-ass-9}) by expanding positive powers of 
$z_{1}=z_{2}+(z_{1}-z_{2})$ using the binomial expansion as polynomials in $z_{2}$ and $z_{1}-z_{2}$. 
Let $\phi'\in (H^{N})'$. Then the same is also true for 
\begin{equation}\label{y-l-ass-11}
(z_{1}-z_{2})^{p_{12}}z_{1}^{r_{1}}z_{2}^{r_{2}}
\langle \phi', Y_{H}^{L}(v_{1}, z_{1})Y_{H}^{L}(v_{2}, z_{2})\phi\rangle
\end{equation}
and 
\begin{equation}\label{y-l-ass-12}
(z_{1}-z_{2})^{p_{12}}z_{1}^{r_{1}}z_{2}^{r_{2}}
\langle \phi', Y_{H}^{L}(Y_{V}(v_{1}, z_{1}-z_{2})v_{2}, z_{2})\phi\rangle.
\end{equation}
But as a special case of the rationality proved above,
(\ref{y-l-ass-11}) is a polynomial $h(z_{1}, z_{2})$ in $z_{1}$ and $z_{2}$. 
So (\ref{y-l-ass-12}) must be a polynomial in $z_{2}$ and $z_{1}-z_{2}$ obtained 
from (\ref{y-l-ass-11}) by expanding the positive powers of $z_{1}=z_{2}+(z_{1}-z_{2})$
using the binomial expansion as polynomials in $z_{2}$ and $z_{1}-z_{2}$. Thus 
\begin{align*}
R&(\langle \phi', Y_{H}^{L}(Y_{V}(v_{1}, z_{1}-z_{2})v_{2}, z_{2})\phi\rangle)\nn
&=\frac{h(z_{1}, z_{2})}{(z_{1}-z_{2})^{p_{12}}z_{1}^{r_{1}}z_{2}^{r_{2}}}\nn
&=R(\langle \phi', Y_{H}^{L}(v_{1}, z_{1})Y_{H}^{L}(v_{2}, z_{2})\phi\rangle.
\end{align*}
This is the associativity of $Y_{H}^{L}$. 

\paragraph{6.} The $D$-derivative property and the $D$-commutator formula: We prove only the 
$D$-derivative property and 
$D$-commutator formula for $Y_{H}^{R}$. The proof of these properties of $Y_{H}^{L}$ 
are similar and are omitted.  Calculating the derivatives and using 
(\ref{y-r-rigor-r}) and the definition of $D_{H}$, we obtain
\begin{align}\label{D-der}
R&\left(\left\langle w_{2}',  \left(\left(\frac{d}{dz_{1}}Y_{H}^{R}(\phi, z_{1})v\right)(w_{1})\right)(z_{2})
\right\rangle\right)\nn
&=\frac{\partial}{\partial z_{1}}R(\langle w_{2}',  ((Y_{H}^{R}(\phi, z_{1})v)(w_{1}))(z_{2})\rangle)\nn
&=\frac{\partial}{\partial z_{1}}R(\langle w_{2}',  \phi(Y_{W_{1}}(v, z_{2})w_{1})(z_{1}+z_{2})\rangle)\nn
&=\left(\frac{\partial}{\partial z_{0}}R(\langle w_{2}',  \phi(Y_{W_{1}}(v, z_{2})w_{1})(z_{0})\rangle)\right)
\lbar_{z_{0}=z_{1}+z_{2}}\nn
&=R(\langle w_{2}',  ((D_{H}\phi)(Y_{W_{1}}(v, z_{2})w_{1}))(z_{1}+z_{2})\rangle)\nn
&=R(\langle w_{2}',  ((Y_{H}^{R}(D_{H}\phi, z_{1})v)(w_{1}))(z_{2})\rangle)
\end{align}
for $w_{2}'\in W_{2}'$ and $w_{1}\in W_{1}$. Thus we obtain the $D$-derivative property
$$\frac{d}{dz_{1}}Y_{H}^{R}(\phi, z_{1})v=Y_{H}^{R}(D_{H}\phi, z_{1})v.$$

Using the equality that the left-hand side of (\ref{D-der}) is equal to the fourth line in (\ref{D-der}),
we obtain
\begin{align*}
R&\left(\left\langle w_{2}',  \left(\left(\frac{d}{dz_{1}}Y_{H}^{R}(\phi, z_{1})v\right)(w_{1})\right)(z_{2})
\right\rangle\right)\nn
&=\left(\frac{\partial}{\partial z_{0}}R(\langle w_{2}',  \phi(Y_{W_{1}}(v, z_{2})w_{1})(z_{0})\rangle)\right)
\lbar_{z_{0}=z_{1}+z_{2}}\nn
&=\frac{\partial}{\partial z_{2}}R(\langle w_{2}',  \phi(Y_{W_{1}}(v, z_{2})w_{1})(z_{1}+z_{2})\rangle)
-R\left(\left\langle w_{2}',  \phi\left(\frac{\partial}{\partial z_{2}}Y_{W_{1}}(v, z_{2})w_{1}\right)
(z_{1}+z_{2})\right\rangle\right)\nn
&=\frac{\partial}{\partial z_{2}}R(\langle w_{2}',  ((Y_{H}^{R}(\phi, z_{1})v)(w_{1}))(z_{2})\rangle)
-R(\langle w_{2}',  \phi(Y_{W_{1}}(D_{V}v, z_{2})w_{1})
(z_{1}+z_{2})\rangle)\nn
&=R(\langle w_{2}',  ((D_{H}Y_{H}^{R}(\phi, z_{1})v)(w_{1}))(z_{2})\rangle)
-R(\langle w_{2}',  ((Y_{H}^{R}(\phi, z_{1})D_{V}v)(w_{1}))(z_{2})\rangle)
\end{align*}
for $w_{2}'\in W_{2}'$ and $w_{1}\in W_{1}$. Thus we obtain 
the $D$-commutator formula for $Y_{H}^{R}$:
$$\frac{d}{dz_{1}}Y_{H}^{R}(\phi, z_{1})v=D_{H}Y_{H}^{R}(\phi, z_{1})v
-Y_{H}^{R}(\phi, z_{1})D_{V}v.$$

All the axioms have been verified and thus the theorem is proved. 
\epfv

For $S\subset H^{N}$, let $H^{(N, S)}$ be the $V$-subbimodule 
of $H^{N}$ generated by $S$.

\renewcommand{\theequation}{\thesection.\arabic{equation}}
\renewcommand{\thethm}{\thesection.\arabic{thm}}
\setcounter{equation}{0}
\setcounter{thm}{0}
\section{A $1$-cocycle constructed from a left $V$-module and a left $V$-submodule}

In this section, given a  left $V$-module $W$ and a
left $V$-submodule $W_{2}$  
of $W$ and assuming a composability condition, 
we construct a $1$-cocycle in $\hat{C}^{1}_{\infty}(V, H^{(N, F(V))})$ where 
$N$ is a lower bound of the weights of the elements of $V$, $F(V)$ 
is the image of a suitable linear map $F$ from $V$ to $H^{N}$, 
$H^{N}$ is the 
$V$-bimodule constructed in the preceding section and $H^{(N, F(V))}$ is the $V$-subbimodule 
of $H^{N}$ generated by $F(V)$. In fact $F(V)$ is independent of $N$. Thus we shall denote 
$H^{(N, F(V))}$ simply by $H^{F(V)}$.

Let $W$ be a  left $V$-module and $W_{2}$ a $V$-submodule of $W$. 
Let $W_{1}$ be a graded subspace of $W$ such that as a graded vector space, we have
$$W=W_{1}\oplus W_{2}.$$
Then we can also embed $W_{1}'$ and $W_{2}'$ into $W'$ and we have
$$W'=W_{1}'\oplus W_{2}'.$$

Let $\pi_{W_{1}}: W\to W_{1}$ and $\pi_{W_{2}}: W\to W_{2}$ be the projections
given by this graded space decomposition of $W$. For simplicity, we shall use 
the same notations $\pi_{W_{1}}$ and $\pi_{W_{2}}$ to denote their 
natural extensions to operators on $\overline{W}_{1}$ and $\overline{W}_{2}$, respectively. 
By definition, we have $\pi_{W_{1}}+\pi_{W_{2}}=1_{W}$, 
$\pi_{W_{1}}\circ \pi_{W_{1}}=\pi_{W_{1}}$, 
$\pi_{W_{2}}\circ \pi_{W_{2}}=\pi_{W_{2}}$, $\pi_{W_{1}}\circ \pi_{W_{2}}=\pi_{W_{2}}\circ 
\pi_{W_{1}}=0$.

Since $W_{2}$ is a submodule of $W$, we have 
$\pi_{W_{2}}\circ Y_{W}\circ (1_{V}\otimes \pi_{W_{2}})=Y_{W_{2}}$,
$D_{W_{2}}=\pi_{W_{2}} \circ D_{W}\circ \pi_{W_{2}}$ and 
$\mathbf{d}_{W_{2}}=\pi_{W_{2}}\circ \mathbf{d}_{W}\circ \pi_{W_{2}}$. We also have 
$\pi_{W_{1}}\circ Y_{W}\circ (1_{V}\otimes \pi_{W_{2}})=0$. 

Let $Y_{W_{1}}=\pi_{W_{1}}\circ Y_{W}\circ (1_{V}\otimes \pi_{W_{1}})$, 
$D_{W_{1}}=\pi_{W_{1}} \circ D_{W}\circ \pi_{W_{1}}$ and 
$\mathbf{d}_{W_{1}}=\pi_{W_{1}}\circ \mathbf{d}_{W}\circ \pi_{W_{1}}$ which is equal to the 
operator giving the grading on $W_{1}$. As we have done above, we use the 
same notations $D_{W_{1}}$ and 
$\mathbf{d}_{W_{1}}$ to denote their natural extensions to $\overline{W_{1}}$. 
We also use the same convention for notations
for extensions of operators on $W$ and $W_{2}$.

\begin{prop}
The graded vector space $W_{1}$ equipped with the vertex operator map $Y_{W_{1}}$
and the operator $D_{W_{1}}$ is a left $V$-module. 
\end{prop}
\pf
The axioms for the grading and the identity property are obvious. 

For $v\in V$,
\begin{align}\label{w-1-l-mod-1}
\frac{d}{dz}Y_{W_{1}}(v, z)
&=\pi_{W_{1}}\left(\frac{d}{dz}Y_{W}(v, z)\right) \pi_{W_{1}}\nn
&=\pi_{W_{1}} Y_{W}(D_{V}v, z)) \pi_{W_{1}}\nn
&=Y_{W_{1}}(D_{V}v, z).
\end{align}
Also using (\ref{w-1-l-mod-1}), we have
\begin{align}\label{w-1-l-mod-2}
\frac{d}{dz}Y_{W_{1}}(v, z)&=
\pi_{W_{1}} Y_{W}(D_{V}v, z)\pi_{W_{1}}\nn
&=\pi_{W_{1}}D_{W}Y_{W}(v, z)\pi_{W_{1}}
-\pi_{W_{1}}Y_{W}(v, z)D_{W} \pi_{W_{1}}\nn
&=\pi_{W_{1}}D_{W}\pi_{W_{1}}Y_{W}(v, z)\pi_{W_{1}}
+\pi_{W_{1}}D_{W}\pi_{W_{2}}Y_{W}(v, z)\pi_{W_{1}}\nn
&\quad -\pi_{W_{1}}Y_{W}(v, z)\pi_{W_{1}}D_{W} \pi_{W_{1}}
-\pi_{W_{1}}Y_{W}(v, z)\pi_{W_{2}}D_{W} \pi_{W_{1}}\nn
&=D_{W_{1}}Y_{W_{1}}(v, z)
+\pi_{W_{1}}D_{W}\pi_{W_{2}}Y_{W}(v, z)\pi_{W_{1}}\nn
&\quad -Y_{W_{1}}(v, z)D_{W_{1}} 
-\pi_{W_{1}}Y_{W}(v, z)\pi_{W_{2}}D_{W} \pi_{W_{1}}.
\end{align}
Since $W_{2}$ is a $V$-submodule of $W$, 
$\pi_{W_{1}}D_{W}\pi_{W_{2}}=\pi_{W_{1}}D_{W_{2}}\pi_{W_{2}}=0$. Again 
since $W_{2}$ is a $V$-submodule of $W$, $\pi_{W_{1}}Y_{W}(v, z)\pi_{W_{2}}
=\pi_{W_{1}}Y_{W_{2}}(v, z)\pi_{W_{2}}=0$. So the right-hand side 
of (\ref{w-1-l-mod-2}) is equal to  
$$D_{W_{1}}Y_{W_{1}}(v, z)-Y_{W_{1}}(v, z)D_{W_{1}}.$$
Thus both the $D$-derivative property and the $D$-commutator formula hold. 

For $v, v_{1}, \dots, v_{k}\in V$, $w_{1}'\in W_{1}'$ and $w_{1}\in W_{1}$, using the properties 
of $\pi_{W_{1}}$, $\pi_{W_{2}}$, $Y_{W_{1}}$, $Y_{W_{2}}$ given above, we have
\begin{equation}\label{w-1-l-mod-3}
\langle w_{1}', Y_{W_{1}}(v_{1}, z_{1})\cdots Y_{W_{1}}(v_{k}, z_{k})
w_{1}\rangle
=\langle w_{1}', Y_{W}(v_{1}, z_{1})\cdots Y_{W}(v_{k}, z_{k})
w_{1}\rangle.
\end{equation}
Since the right-hand side of (\ref{w-1-l-mod-3}) is absolutely convergent in the region 
$|z_{1}|>\cdots >|z_{k}|>0$ to a rational function in $z_{1}, \dots, z_{k}$ with the 
only possible poles $z_{i}=0$ for $i=1, \dots, k$ and $z_{i}=z_{j}$ for $1\le i<j\le k$,
 so is the left-hand side. This proves the rationality. 

For $v, v_{1}, v_{2}\in V$, $w_{1}'\in W_{1}'$ and $w_{1}\in W_{1}$, using 
the properties of $\pi_{W_{1}}$, $\pi_{W_{2}}$, $Y_{W_{1}}$, $Y_{W_{2}}$ again and 
the associativity for $Y_{W}$,
we obtain
\begin{align}\label{w-1-l-mod-4}
\langle w_{1}'&, Y_{W_{1}}(v_{1}, z_{1})Y_{W_{1}}(v_{2}, z_{2})w_{1}\rangle\nn
&=\langle w_{1}', Y_{W}(v_{1}, z_{1})Y_{W}(v_{2}, z_{2})w_{1}\rangle\nn
&=\langle w_{1}', Y_{W}(Y_{V}(v_{1}, z_{1}-z_{2})v_{2}, z_{2})w_{1}\rangle\nn
&=\langle w_{1}', Y_{W_{1}}(Y_{V}(v_{1}, z_{1}-z_{2})v_{2}, z_{2})w_{1}\rangle,
\end{align}
proving the associativity for $Y_{W_{1}}$. 
\epfv

\begin{rema}
{\rm Note that although $W_{1}$ is a graded subspace of $W$, $(W_{1}, Y_{W_{1}})$
is not a submodule of $W$ since the vertex operator $Y_{W_{1}}$ is not 
the restriction of the vertex operator $Y_{W}$ to $V\otimes W_{1}$. }
\end{rema}

We now have two left $V$-modules $W_{1}$ and $W_{2}$. From Theorem \ref{bimodule},
for $N\in \Z$, we have a $V$-bimodule $H^{N}\subset \hom(W_{1}, \widehat{(W_{2})}_{z})$.

We need the following assumption (called {\it composability condition}) on the map 
$\pi_{W_{2}}\circ Y_{W}\circ (1_{V}\otimes \pi_{W_{1}})$:

\begin{assum}[Composability condition]\label{composability}
For $v, v_{1}, \dots, v_{k+l}\in V$, $w_{2}'\in W_{2}'$ and $w_{1}\in W_{1}$, 
\begin{equation}\label{prod}
\langle w_{2}', Y_{W_{2}}(v_{1}, z_{1})\cdots Y_{W_{2}}(v_{k}, z_{k})\pi_{W_{2}}Y_{W}(v, z)\pi_{W_{1}}
Y_{W_{1}}(v_{k+1}, z_{k+1})\cdots Y_{W_{1}}(v_{k+l-1}, z_{k+l})
w_{1}\rangle
\end{equation}
is absolutely convergent in the region $|z_{1}|>\cdots |z_{k}|>|z|>|z_{k+1}|>\cdots >|z_{k+l}|>0$
to a rational function in $z_{1}, \dots, z_{k+l}, z$ with the only possible poles $z_{i}=0$ for $i=1, \dots, k+l$, 
$z=0$, $z_{i}=z_{j}$ for $1\le i<j\le k+l$ and $z_{i}=z$ for $i=1, \dots, k+l$. Moreover, 
the orders of the poles $z_{i}=z$ for $i=1, \dots, k+l$  and 
$z_{i}=z_{j}$ for $i, j=1, \dots, k+l$, $i\ne j$ are bounded above by nonnegative integers
depending only on the pairs $(v_{i}, v)$ and $(v_{i}, v_{j})$, respectively, and there exists $N\in \Z$
such that when (\ref{prod}) is expanded as a Laurent series 
in the region $|z_{k+l}|>|z_{i}-z_{k+l}|>\cdots >|z_{k}-z_{k+l}|>|z-z_{k+l}|>|z_{k+1}|>\cdots
>|z_{k+l-1}-z_{k+l}|>0$  as a Laurent series in $z_{i}-z_{k+l}$ for $i=1, \dots, k+l-1$ and $z-z_{k+l}$
with Laurent polynomials in $z_{k+l}$ as coefficients, the total degree of each monomial 
in $z_{i}-z_{k+l}$ for $i=1, \dots, k+l-1$ and $z-z_{k+l}$ in the expansion 
is larger than or equal to $N-\sum_{i=1}^{k+l}\wt v_{i}+\wt v$. 
\end{assum}

We now assume that $\pi_{W_{2}}\circ Y_{W}\circ (1_{V}\otimes \pi_{W_{1}})$ satisfies the composability condition. 
For $v\in V$, let $F(v)\in \hom(W_{1}, \widehat{(W_{2})}_{z})$ be given by
$$((F(v))(w_{1}))(z)=\pi_{W_{2}}Y_{W}(v, z)\pi_{W_{1}}w_{1}$$
for $w_{1}\in W_{1}$. Thus we obtain a linear map $F: V\to \hom(W_{1}, \widehat{(W_{2})}_{z})$.

\begin{prop}\label{F(v)}
For $N\in \Z$ such that $\wt v\ge N$ for any homogeneous $v\in V$,
the image of $F$ is in fact in $H^{N}$ and is thus a map from $V$ to $H^{N}$.
Moreover, $F$ preserves the gradings.
\end{prop}
\pf
Let $v\in V$ be homogeneous. 
For $a\in \C^{\times}$ and $w_{1}\in W_{1}$, 
\begin{align*}
a^{\mathbf{d}_{W_{2}}}((F(v))(w_{1}))(z)&=a^{\mathbf{d}_{W_{2}}}\pi_{W_{2}}Y_{W}(v, z)\pi_{W_{1}}w_{1}\nn
&=\pi_{W_{2}}Y_{W}(a^{\mathbf{d}_{V}}v, az)\pi_{W_{1}}a^{\mathbf{d}_{W_{1}}}w_{1}\nn
&=a^{\swt v}((F(v))(a^{\mathbf{d}_{W_{1}}}w_{1}))(az),
\end{align*}
proving the $\mathbf{d}$-conjugation property of $F(v)$ and $\wt F(v)=\wt v$.

For $k, l\in \N$ and $v_{1}, \dots, v_{k+l}\in V$, $w_{1}\in W_{1}$ and $w_{2}'\in W_{2}$,
by the definition of $F(v)$, we have
\begin{align}\label{F(v)-1}
&\langle w_{2}', Y_{W_{2}}(v_{1}, z_{1})\cdots Y_{W_{2}}(v_{k}, z_{k})((F(v))(
Y_{W_{1}}(v_{k+1}, z_{k+1})\cdots Y_{W_{1}}(v_{k+l-1}, z_{k+l-1})
w_{1}))(z)\rangle\nn
&=\langle w_{2}', Y_{W_{2}}(v_{1}, z_{1})\cdots Y_{W_{2}}(v_{k}, z_{k})\pi_{W_{2}}Y_{W}(v, z)\pi_{W_{1}}
Y_{W_{1}}(v_{k+1}, z_{k+1})\cdots Y_{W_{1}}(v_{k+l-1}, z_{k+l-1}).
w_{1}\rangle
\end{align}
Then by Assumption \ref{composability}, the composability for $F(v)$ holds. 

Also by Assumption \ref{composability}, the sum of the orders of the possible poles 
$z_{i}=z$ for $i=1, \dots, k+l$  and 
$z_{i}=z_{j}$ for $i, j=1, \dots, k+l$, $i\ne j$ of the rational function that
the right-hand side of (\ref{F(v)-1}) converges to is 
less than or equal to $\sum_{i=1}^{k+l}\wt v_{i}+\wt v-N$. 
By (\ref{F(v)-1}) and the fact $\wt F(v)=\wt v$, we see that the $N$-weight-degree condition for 
$F(v)$ holds. 
\epfv

Proposition \ref{F(v)} says in particular that $F(V)$ is independent of such lower bound $N$
of the weights of $V$. Thus we shall denote $H^{(N, F(V))}$ simply by $H^{F(V)}$.
Now we construct a $1$-cochain $\Psi\in \hat{C}_{\infty}^{1}(V, H^{F(V)})$. Since 
$$\hat{C}_{\infty}^{1}(V, H^{F(V)})\subset \hom(V, \widetilde{(H^{F(V)})}_{z}),$$ 
to avoid confusion with 
the variable in $H^{F(V)}\subset \hom(W_{1}, \widehat{(W_{2})}_{z})$, we use different notations 
to denote these variables. For example, we might use $z_{1}$ and $z_{2}$ to denote the 
variables in $\hom(W_{1}, \widehat{(W_{2})}_{z_{1}})$ and $\hom(V, \widetilde{(H^{F(V)})}_{z_{2}})$,
respectively.
For $v\in V$, $w_{1}\in W_{1}$ and $w_{2}'\in W_{2}'$, 
\begin{align*}
\langle w_{2}', ((e^{z_{2}D_{H}}F(v))(w_{1}))(z_{1})\rangle
&=\langle w_{2}', ((F(v))(w_{1}))(z_{1}+z_{2})\rangle\nn
&=\langle w_{2}', \pi_{W_{2}}Y_{W}(v, z_{1}+z_{2})w_{1}\rangle
\end{align*}
in the region $|z_{1}|>|z_{2}|$.  Let $E(e^{z_{2}D_{H}}F(v))\in \widetilde{(H^{F(V)})}_{z_{2}}$ be
defined by 
\begin{align*}
\langle w_{2}', ((E(e^{z_{2}D_{H}}F(v)))(w_{1}))(z_{1})\rangle
&=\langle w_{2}', \pi_{W_{2}}Y_{W}(v, z_{1}+z_{2})w_{1}\rangle
\end{align*}
$v\in V$, $w_{1}\in W_{1}$ and $w_{2}'\in W_{2}'$ in the region $z_{1}+z_{2}\ne 0$. 

We define 
$$(\Psi(v))(z_{2})=E(e^{z_{2}D_{H}}F(v)).$$
More explicitly, for $v\in V$, $w_{1}\in W_{1}$ and $w_{2}'\in W_{2}'$, 
\begin{align}\label{Psi-expl-def}
\langle w_{2}', (((\Psi(v))(z_{2}))(w_{1}))(z_{1})\rangle
&=\langle w_{2}', \pi_{W_{2}}Y_{W}(v, z_{1}+z_{2})w_{1}\rangle.
\end{align}
in the region $z_{1}+z_{2}\ne 0$.
In the region $|z_{1}|>|z_{2}|$, the series $e^{z_{2}D_{H}}F(v)$ is convergent absolutely to 
$\Psi(v)$. We shall also use $e^{z_{2}D_{H}}F(v)$ to denote $(\Psi(v))(z_{2})$ in the region $|z_{1}|>|z_{2}|$.
By definition, 
$\Psi(v)\in \widetilde{(H^{F(V)})}_{z_{2}}$ and thus $\Psi\in \hom(V, \widetilde{(H^{F(V)})}_{z_{2}})$.

\begin{thm}\label{Psi}
$\Psi\in \ker \hat{\delta}_{\infty}^{1}\subset \hat{C}_{\infty}^{1}(V, H^{F(V)})$.
\end{thm}
\pf
From Theorem \ref{1-coh-der}, we need only show that $\Psi(\cdot)(0)=F$ is a derivation from 
$V$ to $H^{F(V)}$. Using (\ref{y-l-rigor}) and (\ref{y-r-rigor}) in the region 
$|z+z_{1}|>|z_{1}|>|z|>0$, we have
\begin{align}\label{Psi-0}
&\langle w_{2}', ((F(Y_{V}(u, z)v))(w_{1}))(z_{1})\rangle\nn
&\quad =\langle w_{2}', \pi_{W_{2}}Y_{W}(Y_{V}(u, z)v, z_{1})\pi_{W_{1}}w_{1}\rangle\nn
&\quad =\langle w_{2}', \pi_{W_{2}}Y_{W}(u, z+z_{1})Y_{W}(v, z_{1})\pi_{W_{1}}w_{1}\rangle\nn
&\quad =\langle w_{2}', \pi_{W_{2}}Y_{W}(u, z+z_{1})\pi_{W_{1}}
Y_{W}(v, z_{1})\pi_{W_{1}}w_{1}\rangle\nn
&\quad \quad+\langle w_{2}', \pi_{W_{2}}Y_{W}(u, z+z_{1})\pi_{W_{2}}
Y_{W}(v, z_{1})\pi_{W_{1}}w_{1}\rangle\nn
&\quad =\langle w_{2}', ((F(u))(Y_{W_{1}}(v, z_{1})w_{1}))(z+z_{1})\rangle\nn
&\quad \quad+\langle w_{2}', Y_{W_{2}}(u, z+z_{1})((F(v))(w_{1}))(z_{1})\rangle\nn
&\quad =\langle w_{2}', ((Y_{H}^{R}(F(u), z)v)(w_{1}))(z_{1})\rangle\nn
&\quad \quad+\langle w_{2}', ((Y_{H}^{L}(u, z)F(v))(w_{1}))(z_{1})\rangle
\end{align}
for $u, v\in V$, $w_{1}\in W_{1}$ and $w_{2}'\in W_{2}'$. By Proposition
\ref{part-y-h-l-r}, the left-hand side and the two terms in the right-hand side 
of (\ref{Psi-0}) are absolutely convergent in the same region $|z_{1}|>|z|>0$.
Hence the left-hand side and the right-hand side of (\ref{Psi-0}\
are equal in the larger region $|z_{1}|>|z|>0$. Thus we obtain 
$$F(Y_{V}(u, z)v)=Y_{H}^{R}(F(u), z)v+Y_{H}^{L}(u, z)F(v)$$
for $u, v\in V$, proving that $F$ is indeed a derivation. 
\epfv

\renewcommand{\theequation}{\thesection.\arabic{equation}}
\renewcommand{\thethm}{\thesection.\arabic{thm}}
\setcounter{equation}{0}
\setcounter{thm}{0}
\section{The main theorem}

In this section, we formulate and prove our main result  on complete reducibility 
of modules for a meromorphic open-strong vertex algebra $V$.

Let $W$ be a  left $V$-module. Assume that $W$ is not irreducible. 
Then there exists a proper nonzero left $V$-submodule $W_{2}$ of $W$. 
We say that the pair $(W, W_{2})$ {\it satisfies the composability condition} if there exists
a graded subspace $W_{1}$ of $W$ such that $W=W_{1}\oplus W_{2}$ as a graded 
vector space such that $\pi_{W_{2}}\circ Y_{W}\circ (1_{V}\otimes \pi_{W_{1}})$ 
satisfies Assumption \ref{composability}. If for every 
proper nonzero left $V$-submodule $W_{2}$ of $W$, 
the pair $(W, W_{2})$ satisfies the composability condition, we say that {\it $W$ satisfies the composability
condition}. 

\begin{prop}\label{reducibility-composability}
Let $W$ be a completely reducible left $V$-module. Then $W$ satisfies the composability condition. 
\end{prop}
\pf
Let $W_{2}$ be a left $V$-submodule of $W$. Since $W$ is completely reducible, 
there is a left $V$-submodule $W_{1}$ of $W$ such that $W$ as a left $V$-module
is the direct sum of the left $V$-modules $W_{1}$ and $W_{2}$. 
Then $\pi_{W_{1}}$ and $\pi_{W_{2}}$ are module maps. Thus 
$\pi_{W_{2}}\circ Y_{W}\circ (1_{V}\otimes \pi_{W_{1}})$ satisfies Assumption \ref{composability}.
\epfv

Now let $W$ be a  left $V$-module which is not irreducible and $W_{2}$
a proper nonzero left $V$-submodule $W_{2}$ of $W$. Assume that 
the pair $(W, W_{2})$ satisfies the composability condition. Then there exists
a graded subspace $W_{1}$ of $W$ such that as a graded vector space,
$W$ is the direct sum of $W_{1}$ and $W_{2}$ and $\pi_{W_{2}}\circ Y_{W}\circ (1_{V}\otimes \pi_{W_{1}})$
satisfies Assumption \ref{composability}. By Theorem \ref{bimodule}, Proposition \ref{F(v)} and 
Theorem \ref{Psi}, there exist a left $V$-module structure on $W_{1}$, 
a $V$-bimodule $H^{N}\subset \hom(W_{1}, \widehat{(W_{2})}_{z_{1}})$ for a lower bound $N$ of $V$,
a grading preserving linear map $F: V\to H^{N}$ and 
$\Psi\in \ker \hat{\delta}_{\infty}^{1}\subset \hat{C}_{\infty}^{1}(V, H^{F(V)})$, where 
$H^{F(V)}$ is the $V$-subbimodule of $H^{N}$ generated by $F(V)$.

We see from Proposition \ref{wt-1-der} that besides inner derivations, 
there are also zero-mode derivations obtained from suitable elements of a bimodule.
For example, when $V$ is a vertex algebra and $W=V$, then every weight $1$ element of $V$ 
gives a derivation from $V$ to $V$. In particular, the first cohomology $\hat{H}_{\infty}^{1}(V, V)$
is not $0$, even if all $V$-modules are completely reducible. Moreover,  we have the following
result:

\begin{prop}
Let $V$ be a grading-restricted vertex algebra generated by its homogeneous subspace 
$V_{(1)}$ of weight $1$ elements. 
Assume that the finite-dimensional Lie algebra $V_{(1)}$ with the Lie bracket given by 
$(u, v)\mapsto (Y_{V})_{0}(u)v$ for $u, v\in V_{(1)}$ is semisimple. Then $\hat{H}_{\infty}^{1}(V, V)
=\hat{Z}_{\infty}^{1}(V, V)$ and $\hat{Z}_{\infty}^{1}(V, V)$ is linearly isomorphic to 
$V_{(1)}$. In particular, $\hat{H}_{\infty}^{1}(V, V)$ is linearly isomorphic to 
$V_{(1)}$. 
\end{prop}
\pf
Since $V$ is a vertex algebra, an inner derivation from $V$ to $V$ is $0$. To prove this result, 
we need only prove that every derivation from $V$ to $V$ is a zero-mode derivation and every 
zero-mode derivation can be identified linearly with an element of $V_{(1)}$. 

Let $f$ be a derivation from $V$ to $V$. For $u, v\in V_{(1)}$, from $f(Y_{V}(u, x)v)=Y_{V}(f(u), x)v+
Y_{V}(u, x)f(v)$, we obtain $f((Y_{V})_{0}(u)v)=(Y_{V})_{0}(f(u))v+(Y_{V})_{0}(u)f(v)$. 
So the restriction of $f$ to $V_{(1)}$ is a derivation of the Lie algebra $V_{(1)}$. 
Since this Lie algebra is semisimple, $f$ must be an inner derivation of this Lie algebra. 
Thus there exists a unique $w\in V_{(1)}$ such that $f(v)=(Y_{V})_{0}(w)v$ for $v\in V_{(1)}$. 
Since $w\in V_{(1)}$, $g_{w}: V\to V$ defined by $g_{w}(v)=(Y_{V})_{0}(w)v$ for $v\in V$
is a zero-mode derivation from $V$ to $V$. Since $V$ is generated by $V_{(1)}$ and both 
$f$ and $g_{w}$ are derivations from $V$ to $V$, $f(v)=g_{w}(v)$ for $v\in V_{(1)}$ implies 
$f=g_{w}$. 
\epfv

From this result, we see that even for an affine Lie algebra vertex operator algebra $V$
associated to a  finite-dimensional simple Lie algebra $\mathfrak{g}$
such that all weak $V$-modules are completely reducible, $\hat{H}_{\infty}^{1}(V, V)
=\hat{Z}_{\infty}^{1}(V, V)$ which is in turn isomorphic to $\mathfrak{g} \ne 0$.
However, we have the 
following result:

\begin{thm}\label{mosva-red}
Let $W$, $W_{1}$, $W_{2}$ and $H^{F(V)}$ be as above. 
If  $\hat{H}_{\infty}^{1}(V, H^{F(V)})=\hat{Z}_{\infty}^{1}(V, H^{F(V)})$ (see 
the end of Section 2), then 
there exists another left $V$-submodule $\widetilde{W}_{1}$ of $W$ such that $W$ is the
direct sum of the left $V$-submodules $\widetilde{W}_{1}$ and $W_{2}$.
\end{thm}
\pf
Since $\hat{H}^{1}_{\infty}(V, H^{F(V)})=\hat{Z}_{\infty}^{1}(V, H^{F(V)})$, 
$\Psi$ constructed in the preceding section must be the sum of a coboundary and 
a $0$-cochain obtained from a zero-mode derivation. That is, 
there exist $\Phi_{1}\in \hat{C}^{0}_{\infty}(V, H)$ and $\Phi_{2}\in 
H^{F(V)}_{[1]}$ satisfying 
\begin{equation}\label{main--1}
e^{zD_{H}}Y_{H}^{L}(v, -z)\Phi_{2}-Y_{H}^{R}(\Phi_{2}, z)v\in H^{F(V)}[z, z^{-1}].
\end{equation}
for $v\in V$
such that 
\begin{equation}\label{main-0}
(\Psi(v))(z_{2})=(\hat{\delta}_{\infty}^{0} \Phi_{1})(v)(z_{2})+e^{z_{2}D_{H}}
\res_{z}Y_{H}^{R}(\Phi_{2}, z)v.
\end{equation}
Note that 
the $0$-cochain $\Phi_{1}$ is an element
of $H^{F(V)}_{[0]}$ such that  $D_{H}\Phi_{1}=0$ and, in particular, 
$((\hat{\delta}_{\infty}^{0} \Phi_{1})(v))(z)$ is an $\overline{H^{F(V)}}$-valued holomorphic
function on $\C$. By the definition of $((\hat{\delta}_{\infty}^{0} \Phi_{1})(v))(z_{2})$, we obtain 
from (\ref{main-0})
\begin{equation}\label{main-1}
(\Psi(v))(z_{2})=Y_{H}^{L}(v, z_{2})\Phi_{1}-e^{z_{2}D_{H}}Y_{H}^{R}(\Phi_{1}, -z_{2})v
+e^{z_{2}D_{H}}
\res_{z}Y_{H}^{R}(\Phi_{2}, z)v.
\end{equation}
Applying both sides of (\ref{main-1}) to $w_{1}\in W_{1}$, evaluating at $z_{1}$, pairing with 
$w'\in W'$, using (\ref{Psi-expl-def}) and properties of $D_{H}$ and then taking the limit 
$z_{2}\to 0$,
we obtain in the region $|z_{1}|>|z_{2}|>0$, 
\begin{align}\label{main-2}
\langle w'&, \pi_{W_{2}}Y_{W}(v, z_{1})w_{1}\rangle\nn
&=\langle w',  ((\lim_{z_{2}\to 0}Y_{H}^{L}(v, z_{2})\Phi_{1})(w_{1}))(z_{1})\rangle
-\langle w',  ((\lim_{z_{2}\to 0}Y_{H}^{R}(\Phi_{1}, -z_{2})v)(w_{1}))(z_{1})\rangle\nn
&\quad +\res_{z} \langle w', 
((Y_{H}^{R}(\Phi_{2}, z)v)(w_{1}))(z_{1})\rangle.
\end{align}

Note that since $D_{H}\Phi_{1}=0$, $(\Phi_{1}(w_{1}))(z_{1})$
is in fact independent of $z_{1}$. In particular, $(\Phi_{1}(w_{1}))(0)$ exists and 
$(\Phi_{1}(w_{1}))(z)=(\Phi_{1}(w_{1}))(0)$ for $z\in \C$.
We now define a linear map $\eta: W_{1}\to W$ by
\begin{equation}\label{main-2.5}
\eta(w_{1})=w_{1}-(\Phi_{1}(w_{1}))(0)+\res_{z}(\Phi_{2}(w_{1}))(z)
\end{equation}
for homogeneous $w_{1}\in W_{1}$. 
Note that the weight of $(\Phi_{1}(w_{1}))(0)$ is equal to the weight of $w_{1}$. It is clear that 
the weight of $\res_{z}(\Phi_{2}(w_{1}))(z)$ is also equal to the weight of $w_{1}$. 
Thus the map $\eta$ preserves the gradings. 
Let $\widetilde{W}_{1}=\eta(W_{1})$.
We show that $\widetilde{W}_{1}$ is in fact a left $V$-submodule of $W$. 

Applying $Y_{W}(v, z_{1})$ to the right-hand side of (\ref{main-2.5}) and pairing the result with 
$w'\in W'$,
using (\ref{main-2})
and then using (\ref{y-l-rigor-r}) and (\ref{y-r-rigor-r}), we obtain 
\begin{eqnarray}\label{main-3}
\lefteqn{\langle w', Y_{W}(v, z_{1})(w_{1}-(\Phi_{1}(w_{1}))(0)
+\res_{z}(\Phi_{2}(w_{1}))(z))\rangle}\nn
&&=\langle w', (\pi_{W_{1}}Y_{W}(v, z_{1})w_{1}+\pi_{W_{2}}Y_{W}(v, z_{1})w_{1}\nn
&&\quad\quad \quad\quad  \quad 
-Y_{W_{2}}(v, z_{1})(\Phi_{1}(w_{1}))(0)+\res_{z}
Y_{W_{2}}(v, z_{1})(\Phi_{2}(w_{1}))(z))\rangle\nn
&&=\langle w', (Y_{W_{1}}(v, z_{1})w_{1}+
((\lim_{z_{2}\to 0}Y_{H}^{L}(v, z_{2})\Phi_{1})(w_{1}))(z_{1})\nn
&&\quad\quad \quad \quad \quad  -((\lim_{z_{2}\to 0}Y_{H}^{R}(\Phi_{1}, -z_{2})v)(w_{1}))(z_{1})
+\res_{z} 
((Y_{H}^{R}(\Phi_{2}, z)v)(w_{1}))(z_{1})\nn
&&\quad\quad \quad \quad \quad -Y_{W_{2}}(v, z_{1})(\Phi_{1}(w_{1}))(0)
+\res_{z}
Y_{W_{2}}(v, z_{1})(\Phi_{2}(w_{1}))(z))\rangle\nn
&&=\langle w', (Y_{W_{1}}(v, z_{1})w_{1}
-(\Phi_{1}(Y_{W_{1}}(v, z_{1})w_{1}))(0)\nn
&&\quad\quad \quad \quad \quad +\res_{z} 
((Y_{H}^{R}(\Phi_{2}, z)v)(w_{1}))(z_{1}) +\res_{z}
Y_{W_{2}}(v, z_{1})(\Phi_{2}(w_{1}))(z))\rangle.
\end{eqnarray}

We now calculate the last two terms in (\ref{main-3}). From (\ref{y-r-rigor-r}),
we have 
\begin{equation}\label{main-3.1}
R(\langle w', ((Y_{H}^{R}(\Phi_{2}, z_{2})v)(w_{1}))(z_{1}) \rangle)
=R(\langle w', (\Phi_{2}(Y_{W_{1}}(v, z_{1})w_{1}))(z_{1}+z_{2}) \rangle).
\end{equation}
Using (\ref{main--1}) and (\ref{y-l-rigor-r}), we see that there exists a rational function 
$g(z_{1}, z_{2})$ with the only possible pole $z_{1}=0$ and $z_{2}=0$ such that  
\begin{align}\label{main-3.2}
R&(\langle w', ((Y_{H}^{R}(\Phi_{2}, z_{2})v)(w_{1}))(z_{1}) \rangle)\nn
&=R(\langle w', ((e^{z_{2}D_{H}}Y_{H}^{L}(v, -z_{2})\Phi_{2})(w_{1}))(z_{1}) \rangle)+g(z_{1}, z_{2})\nn
&=R(\langle w', ((Y_{H}^{L}(v, -z_{2})\Phi_{2})(w_{1}))(z_{1}+z_{2}) \rangle)+g(z_{1}, z_{2})\nn
&=R(\langle w', Y_{W_{2}}(v, z_{1})(\Phi_{2}(w_{1}))(z_{1}+z_{2}) \rangle)+g(z_{1}, z_{2}).
\end{align}
Let $z_{2}=z-z_{1}$ in (\ref{main-3.1}) and (\ref{main-3.2}). Then we obtain 
\begin{align}\label{main-3.3}
R&(\langle w', ((Y_{H}^{R}(\Phi_{2}, z-z_{1})v)(w_{1}))(z_{1}) \rangle)\nn
&=R(\langle w', (\Phi_{2}(Y_{W_{1}}(v, z_{1})w_{1}))(z) \rangle)\nn
&=R(\langle w', Y_{W_{2}}(v, z_{1})(\Phi_{2}(w_{1}))(z) \rangle)+g(z_{1}, z-z_{1}).
\end{align}
Let $C_{\infty}$, $C_{z_{1}}$ and $C_{0}$ be circles on the $z$-complex plane centered at $0$, $z_{1}$ and $0$,
respectively, with radii larger than $|z_{1}|$, less than $|z_{1}|$ and less than $|z_{1}|$, respectively.
Since (\ref{main-3.3}) is a rational function in $z$ and $z_{1}$ with the only possible poles 
$z=0$, $z_{1}=0$ and $z=z_{1}$, by the convergence regions of the series 
$\langle w', (\Phi_{2}(Y_{W_{1}}(v, z_{1})w_{1}))(z) \rangle$, 
$\langle w', ((Y_{H}^{R}(\Phi_{2}, z-z_{1})v)(w_{1}))(z_{1}) \rangle$,
$\langle w', Y_{W_{1}}(v, z_{1})(\Phi_{2}(w_{1}))(z) \rangle$, 
Cauchy's integral theorem and the fact that $g(z_{1}, z-z_{1})$ is analytic on the disk with $C_{0}$ as 
its boundary as a function in $z$, we have
\begin{align*}
&\res_{z}\langle w', (\Phi_{2}(Y_{W_{1}}(v, z_{1})w_{1}))(z) \rangle\nn
&\quad =\oint_{C_{\infty}}\langle w', (\Phi_{2}(Y_{W_{1}}(v, z_{1})w_{1}))(z) \rangle dz\nn
&\quad =\oint_{C_{\infty}}R(\langle w', (\Phi_{2}(Y_{W_{1}}(v, z_{1})w_{1}))(z) \rangle)dz
\quad\quad\quad\quad\quad\quad\quad\quad\quad\quad\quad\quad\quad\quad\quad\quad
\end{align*}
\begin{align}\label{main-3.4}
&\quad=\oint_{C_{z_{1}}}R(\langle w', ((Y_{H}^{R}(\Phi_{2}, z-z_{1})v)(w_{1}))(z_{1}) \rangle)dz\nn
&\quad  \quad +\oint_{C_{0}}(R(\langle w', Y_{W_{2}}(v, z_{1})(\Phi_{2}(w_{1}))(z) \rangle)+g(z_{1}, z-z_{1}))dz\nn
&\quad=\oint_{C_{z_{1}}}\langle w', ((Y_{H}^{R}(\Phi_{2}, z-z_{1})v)(w_{1}))(z_{1}) \rangle dz
+\oint_{C_{0}}\langle w', Y_{W_{2}}(v, z_{1})(\Phi_{2}(w_{1}))(z) \rangle dz\nn
&\quad=\res_{z}\langle w', ((Y_{H}^{R}(\Phi_{2}, z)v)(w_{1}))(z_{1}) \rangle
 +\res_{z}\langle w', Y_{W_{2}}(v, z_{1})(\Phi_{2}(w_{1}))(z) \rangle.
\end{align}
From (\ref{main-3.4}) and (\ref{main-3}), we obtain
\begin{align*}
\langle w', &\; Y_{W}(v, z_{1})(w_{1}-(\Phi_{1}(w_{1}))(0)
+\res_{z}(\Phi_{2}(w_{1}))(z))\rangle\nn
&=\langle w', (Y_{W_{1}}(v, z_{1})w_{1}
-(\Phi_{1}(Y_{W_{1}}(v, z_{1})w_{1}))(0)+\res_{z} (\Phi_{2}(Y_{W_{1}}(v, z_{1})w_{1}))(z) )\rangle,
\end{align*}
or equivalently,
\begin{equation}\label{main-4}
Y_{W}(v, z_{1})\eta(w_{1})=\eta(Y_{W_{1}}(v, z_{1})w_{1}).
\end{equation}
The formula (\ref{main-4}) means in particular that 
the space $\widetilde{W}_{1}$ is invariant under the action of 
$Y_{W}$. Thus $\widetilde{W}_{1}$ is a submodule of $W$. 
Moreover, the sum of $\widetilde{W}_{1}$ and $W_{2}$ is clearly $W$ and 
the intersection of $\widetilde{W}_{1}$ and $W_{2}$ is clearly $0$. So 
$W$ is equal to the direct sum of $\widetilde{W}_{1}$ and $W_{2}$, proving the theorem.
\epfv

We shall need the following result:

\begin{prop}\label{submod-composability}
Let $W$ be a left $V$-module satisfying the composability condition. 
Then every left $V$-submodule of $W$ also satisfies the composability condition.
\end{prop}
\pf
Let $W_{0}$ be a left $V$-submodule of $W$. Then any proper nonzero 
left $V$-submodule $W_{2}$ of $W_{0}$ is also a proper nonzero  
left $V$-submodule of $W$. Then there is a graded subspace $W_{3}$ of $W$
such that $W=W_{3}\oplus W_{2}$ as a graded vector space and
 $\pi^{W}_{W_{2}}\circ Y_{W}\circ (1_{V}\otimes \pi^{W}_{W_{3}})$ 
satisfies Assumption \ref{composability} with 
$W_{1}$ in Assumption \ref{composability} replaced by $W_{3}$, where 
$\pi^{W}_{W_{2}}$ and $\pi^{W}_{W_{3}}$ are projections from $W$ to 
$W_{2}$ and $W_{3}$, respectively. Let $W_{1}=W_{3}\cap W_{0}$. 
Then $W_{0}=W_{1}\oplus W_{2}$ as a graded vector space. 
Let $\pi^{W_{0}}_{W_{1}}$ and $\pi^{W_{0}}_{W_{2}}$ be the projections from 
$W_{0}$ to $W_{1}$ and $W_{2}$, respectively. Then $\pi^{W_{0}}_{W_{1}}
=\left.\pi^{W}_{W_{3}}\right|_{W_{0}}$ and $\pi^{W_{0}}_{W_{2}}=\left.\pi^{W}_{W_{2}}\right|_{W_{0}}$.
So we have
$$\pi^{W_{0}}_{W_{2}}\circ Y_{W_{0}}\circ (1_{V}\otimes \pi^{W_{0}}_{W_{1}})
=\left.\pi^{W}_{W_{2}}\circ Y_{W}\circ (1_{V}\otimes \pi^{W}_{W_{3}})\right|_{V\otimes W_{0}}.$$
Since $\pi^{W}_{W_{2}}\circ Y_{W}\circ (1_{V}\otimes \pi^{W}_{W_{3}})$ 
satisfies Assumption \ref{composability}, 
$\left.\pi^{W}_{W_{2}}\circ Y_{W}\circ (1_{V}\otimes \pi^{W}_{W_{3}})\right|_{V\otimes W_{0}}$
satisfies Assumption \ref{composability} with $W$ in Assumption \ref{composability} replaced by $W_{0}$. Thus 
$\pi^{W_{0}}_{W_{2}}\circ Y_{W_{0}}\circ (1_{V}\otimes \pi^{W_{0}}_{W_{1}})$ 
satisfies Assumption \ref{composability} with $W$ replaced by $W_{0}$.
\epfv

{


From Theorem \ref{mosva-red} and Proposition \ref{submod-composability}, 
we obtain immediately the following main result of 
this paper:

\begin{thm}\label{main}
Let $V$ be a meromorphic open-string vertex algebra. If 
$\hat{H}_{\infty}^{1}(V, M)=\hat{Z}_{\infty}^{1}(V, M)$ for every $\Z$-graded $V$-bimodule $M$,  
then every left $V$-module satisfying the composability condition 
is completely reducible. 
Assume in addition that the following condition also holds:
For every left $V$-module $W$ satisfying the composability condition and
every nonzero proper left $V$-submodule $W_{2}$ of $W$,
there exists a graded
subspace $W_{1}$ of $W$ such that $W=W_{1}\oplus W_{2}$ as a graded vector space,
$\pi_{W_{2}}\circ Y_{W}\circ (1_{V}\otimes \pi_{W_{1}})$
satisfies Assumption \ref{composability} and the submodule $H^{F(V)}$ of $H^{N}$ is grading restricted
for the grading preserving linear map $F: V\to H^{N}$ given by 
Proposition \ref{F(v)}.
Then the conclusion still holds if  $\hat{H}_{\infty}^{1}(V, M)=\hat{Z}_{\infty}^{1}(V, M)$ only for every
grading-restricted $\Z$-graded $V$-bimodule $M$. 
\end{thm}
\pf
Let $W$ be a left $V$-module satisfying the composability condition
and let $W_{0}$ be the direct sum of all irreducible left $V$-submodule of $W$. We want to show
$W_{0}=W$. If $W_{0}\ne W$, there exists $w\in W$ but $w\not\in W_{0}$. 
By Zorn's lemma, we can find a 
left $V$-submodule $W_{2}$ of $W$ such that $W_{2}$ is maximal among 
left $V$-submodules containing $W_{0}$ but not $w$. Clearly $W_{2}\ne W$. 
Since $W$ satisfies the composability condition,
there is a graded 
subspace $W_{1}$ of $W$ such that   $W=W_{1}\oplus W_{2}$ as a vector space and 
$\pi_{W_{2}}\circ Y_{W}\circ (1_{V}\otimes \pi_{W_{1}})$  satisfies Assumption \ref{composability}. 
By Theorem \ref{mosva-red}, there is a nonzero left $V$-submodule 
$\widetilde{W}_{1}$ of $W$ such that 
$W=\widetilde{W}_{1}\oplus W_{2}$ as a left $V$-module. 
Since $W_{0}\subset W_{2}$,
$\widetilde{W}_{1}$ cannot be irreducible. Let $W_{4}$ be a nonzero proper 
left $V$-submodule of $\widetilde{W}_{1}$. 
By Proposition \ref{submod-composability}, $\widetilde{W}_{1}$ also satisfies 
the composability condition. Therefore
there is a graded 
subspace $W_{3}$ of $W$ such that   $\widetilde{W}_{1}=W_{3}\oplus W_{4}$ as a vector space and 
$\pi_{W_{4}}\circ Y_{\widetilde{W}_{1}}\circ (1_{V}\otimes \pi_{W_{3}})$  satisfies Assumption \ref{composability}. By Theorem \ref{mosva-red}, there is a nonzero proper left $V$-submodule 
$\widetilde{W}_{3}$ of $\widetilde{W}_{1}$ such that 
$\widetilde{W}_{1}=\widetilde{W}_{3}\oplus W_{4}$ as a left $V$-module. 
Thus we have $W=\widetilde{W}_{3}\oplus W_{4}\oplus W_{2}$. 

Since $W_{2}$ is maximal among 
left $V$-submodules containing $W_{0}$ but not $w$ and $\widetilde{W}_{3}\oplus W_{2}$
contains $W_{2}$,  $\widetilde{W}_{3}\oplus W_{2}$ must contain $w$. 
Similarly,   $W_{4}\oplus W_{2}$ must contain $w$. Thus 
$w\in (\widetilde{W}_{3}\oplus W_{2})\cap (W_{4}\oplus W_{2})=W_{2}$.
Contradiction. So $W_{0}=W$.

In the case that the additional condition also holds, since $H^{F(V)}$ is grading restricted, 
$\hat{H}_{\infty}^{1}(V, M)=\hat{Z}_{\infty}^{1}(V, M)$ for every
grading-restricted $\Z$-graded $V$-bimodule $M$ implies in particular 
$\hat{H}_{\infty}^{1}(V, H^{F(V)})=\hat{Z}_{\infty}^{1}(V, H^{F(V)})$.
Then the conclusion of Theorem \ref{mosva-red} still holds. Thus $W$ is completely reducible. 
\epfv

\begin{rema}\label{main-der}
{\rm By Theorems \ref{1-coh-der} and \ref{main}, we can replace the condition 
$\hat{H}_{\infty}^{1}(V, M)=\hat{Z}_{\infty}^{1}(V, M)$
by the statement that every derivation from $V$ to $M$ is the sum of an inner derivation and 
a zero-mode derivation.}
\end{rema}

We also have the following conjecture:

\begin{conj}
Let $V$ be a meromorphic open-string vertex algebra. 
If every left $V$-module is completely reducible,
then $\hat{H}_{\infty}^{1}(V, M)=\hat{Z}_{\infty}^{1}(V, M)$ for every $V$-bimodule $M$.
\end{conj}

Note that because of Proposition \ref{reducibility-composability}, we do not need to 
require that left $V$-modules satisfy the composability condition in this conjecture.

\section{Application to the reductivity of vertex algebras}

In this section, we apply our results in the preceding section to lower-bounded and grading-restricted 
vertex algebras to obtain results on the complete reducibility of $V$-modules. 
See Remark \ref{va-defn} for our definitions of lower-bounded and grading-restricted 
vertex algebra.

We prove that when $V$ is a grading-restricted  vertex algebra containing a vertex subalgebra 
$V_{0}$ satisfying certain natural conditions (see below), then every grading-restricted 
left $V$-module satisfies the composability condition and thus in this case,
the results in the preceding sections holds for all grading-restricted left $V$-modules. 
We also prove in this section that in this case, $F(V)$ (see the preceding section) is grading restricted 
 and thus in this case, the condition that $\hat{H}_{\infty}^{1}(V, M)=\hat{Z}_{\infty}^{1}(V, M)$ for 
every $\Z$-graded $V$-bimodule $M$ in Theorem \ref{main} can be
weakened to require 
that  $\hat{H}_{\infty}^{1}(V, M)=\hat{Z}_{\infty}^{1}(V, M)$ only for every $\Z$-graded 
$V$-bimodule $M$ generated by a grading-restricted subspace. 

Before we give our applications, we need to clarify the terminology on modules. 
Note that in Section 2, left modules, right modules and bimodules are graded.
When $V$ is a vertex operator algebra (in particular, there are Virasoro operators acting on 
$V$), a lower-bounded generalized $V$-module $W$ (that is, a weak $V$ module $W$ with a 
grading given by generalized eigenspaces of $L(0)$ with a lower bound 
on the weights of the elements of $W$) is a left $V$-module defined 
in Section 2 when $V$ is viewed as a meromorphic open-string vertex algebra. 
A grading-restricted generalized $V$-module is a grading-restricted left $V$-module 
defined in Section 2. But to avoid confusion, we shall use the terminology of this paper to 
call these generalized $V$-modules left $V$-modules or grading-restricted left $V$-modules,
though for a vertex algebra, the notions of left module and right module are equivalent. 
The reader should not confuse left $V$-modules in this section with $V$-modules in the literature 
on vertex (operator) algebras. We shall also use the same terminology for modules for subalgebras of $V$. 
We also note that a $V$-bimodule when $V$ is viewed as a meromorphic open-string vertex algebra
is not the same as a left $V$-module. See the discussion below.

First, Theorem \ref{main} in the case that $V$ is lower-bounded vertex algebra becomes:

\begin{thm}\label{va}
Let $V$ be a  lower-bounded vertex algebra (in particular, a grading-restricted vertex algebra or 
a vertex operator algebra). If $\hat{H}_{\infty}^{1}(V, M)=\hat{Z}_{\infty}^{1}(V, M)$
(or every derivation from $V$ to $M$ is the sum of an inner derivation and 
a zero-mode derivation) for   every $\Z$-graded $V$-bimodule $M$,
then every left $V$-module 
satisfying the composability condition is completely reducible. 
\end{thm}

\begin{rema}
{\rm 
Note that when $V$ is a grading-restricted  vertex algebra (in particular, a 
vertex operator algebra) and $M$ is a left $V$-module, we have the chain 
complexes $C^{n}_{\infty}(V, M)$ and the cohomologies $H_{\infty}^{n}(V, M)$ introduced in 
\cite{Hcoh}.
In the case $n=1$, if we view $M$ as a $V$-bimodule with the right $V$-module structure
obtained from the left $V$-module structure using the skew-symmetry formula, 
then $C^{1}_{\infty}(V, M)=\hat{C}^{1}_{\infty}(V, M)$ and $H_{\infty}^{1}(V, M)
=\hat{H}_{\infty}^{1}(V, M)$.
But in Theorem \ref{va}, the assumption is that
$\hat{H}_{\infty}^{1}(V, M)=\hat{Z}_{\infty}^{1}(V, M)$  
for all $V$-bimodules $M$, including those $V$-bimodules for which 
the right actions are not obtained from the left actions. 
For a $V$-bimodule $M$ when $V$ is viewed as
a meromorphic open-string vertex algebra, $M$ is in fact a $\C$-graded space $M$ equipped with two 
commuting $V$-module structures $Y_{M}^{1}$ and $Y_{M}^{2}$ 
(commuting in the sense of the  
commutativity)  when $V$ is viewed as a 
grading-restricted vertex algebra such that the operator $D_{M}^{1}$ and $D_{M}^{2}$ 
for the two $V$-module structures are equal (denoted simply by $D_{M}$). 
Then an inner derivation from $V$ to $M$ is a linear map 
$f_{w}$ from $V$ to $W$ defined by 
$$f_{w}(v)=\lim_{z\to 0}(Y_{M}^{1}(v, z)w-Y_{M}^{2}(v, z)w)
=(Y_{M}^{1})_{-1}(v)w-(Y_{M}^{2})_{-1}(v)w$$
for $v\in V$, where $w\in W$ is of weight $0$ and satisfies $D_{M}w=0$. On the other hand,
a zero-mode derivation is a linear map $g_{w}$ from  $V$ to $W$ defined by 
$$g_{w}(v)=\res_{z}e^{zD_{M}}Y_{M}^{2}(v, -z)w)=\sum_{n\in \N}\frac{(-1)^{-n-1}}{n!}
D_{M}^{n}(Y_{M}^{2})_{n}(v)w$$
for $v\in V$, where $w\in W$ is of weight $1$ and satisfies 
$Y_{M}^{1}(v, z)w-Y_{M}^{2}(v, z)w\in M[z, z^{-1}]$.
The condition $\hat{H}_{\infty}^{1}(V, M)=\hat{Z}_{\infty}^{1}(V, M)$ is 
now equivalent to the statement that 
every derivation $f$ from $V$ to $M$ is of the sum of derivations of the forms $f_{w}$
and $g_{w}$ as above (certainly with different $w$'s). }
\end{rema}

The remark above provides a method to determine whether  
$\hat{H}_{\infty}^{1}(V, M)=\hat{Z}_{\infty}^{1}(V, M)$ when $V$ is 
a grading-restricted vertex algebra (in particular, a vertex operator 
algebra). But  this method is practical only if we know all bimodules, not just those 
modules for which the right module structure is obtained from the left module
structure using the skew-symmetry. These will have to be studied 
in future papers. 

In Theorem \ref{va}, only left $V$-modules satisfying the composability condition
are completely reducible.  We now show that in the case that $V$ has a nice vertex subalgebra,
this condition holds for every grading-restricted left $V$-module. 
Certainly the existence of such a vertex subalgebra is a very strong condition.
But that this condition follows from the existence of such a vertex subalgebra
is an evidence that this condition is a natural and important condition. More importantly, 
our proof below uses the theory of intertwining operators and the tensor category theory 
for module categories for grading-restricted vertex algebras satisfying suitable conditions.
This use reveals the deep connection between the cohomology theory and the theory of 
intertwining operators or the tensor category theory.

In the remaining part of this section, $V$ is a grading-restricted vertex algebra and 
$V_{0}$ is a vertex subalgebra of $V$. Note that the grading and the operator $D_{V_{0}}$ 
for the vertex subalgebra $V_{0}$ of $V$ are induced from those  for $V$. 
In particular, when $V$ is a vertex operator algebra such that $\mathbf{d}_{V}=L_{V}(0)$ and 
$D_{V}=L_{V}(-1)$ and $V_{0}$ is a vertex subalgebra with a different conformal vector,
$\mathbf{d}_{V_{0}}=L_{V}(0)|_{V_{0}}$ and $D_{V_{0}}=L_{V}(-1)|_{V_{0}}$ are in general 
different from  $L_{V_{0}}(0)$ and $L_{V_{0}}(-1)$, respectively. 

Now $V$ is a left $V_{0}$-module and any left $V$-module is also a left $V_{0}$-module.
Here we emphasize that the left $V_{0}$-module structure on a left $V$-module 
$W$ has the the same grading and the operator $D_{W}$ as the 
grading and the operator $D_{W}$ for the $V$-module structure. 
Moreover, we have:

\begin{prop}\label{int-op}
Assume that every grading-restricted left $V_{0}$-module (or generalized 
$V_{0}$-module satisfying the grading-restriction conditions 
in the terminology in \cite{HLZ}) is completely reducible. Let $W$ be a grading-restricted 
left $V$-module and $W_{2}$ a left $V$-submodule of $W$. 
Then there exists a left $V_{0}$-submodule $W_{1}$ of the left $V_{0}$-module $W$ such that 
$W=W_{1}\oplus W_{2}$ as a left $V_{0}$-module and the map 
$\pi_{W_{2}}\circ Y_{W}\circ (1_{V}\otimes \pi_{W_{1}})$ is an intertwining operator 
of type $\binom{W_{2}}{VW_{1}}$ for the left $V_{0}$-modules $V$, $W_{1}$ and $W_{2}$. 
\end{prop}
\pf
The existence of $W_{1}$ follows from the complete reducibility of the grading-restricted left 
$V_{0}$-module. 
When $V$ and $W$ are viewed as left $V_{0}$-modules, $Y_{W}$ is in fact an intertwining 
operator of type ${W\choose VW}$. Since $W_{1}$, $W_{2}$ and $W$ are all left $V_{0}$-modules
and $W=W_{1}\oplus W_{2}$ as a left $V_{0}$-module, the projections $\pi_{W_{1}}$ and $\pi_{W_{2}}$
are $V_{0}$-module maps. Thus the compositions of these module maps with $Y_{W}$
are still intertwining operators among left $V_{0}$-modules. In particular, 
$\pi_{W_{2}}\circ Y_{W}\circ (1_{V}\otimes \pi_{W_{1}})$ is an intertwining operator 
of type $\binom{W_{2}}{VW_{1}}$ for the $V_{0}$-modules $V$, $W_{1}$ and $W_{2}$. 
\epfv

\begin{thm}\label{va-composability}
Assume the following conditions holds for intertwining operators among grading-restricted left $V_{0}$-modules:

\begin{enumerate}

\item  For any $n\in \Z_{+}$, 
products of $n$ intertwining operators among grading-restricted left $V_{0}$ modules 
evaluated at $z_{1}, \dots, z_{n}$
are absolutely convergent in the region $|z_{1}|>\cdots >|z_{n}|>0$ and
can be analytically extended to  (possibly multivalued) analytic functions in $z_{1}, \dots, z_{n}$
with the only possible 
singularities (branch points or poles)
$z_{i}=0$ for $i=1, \dots, n$ and $z_{i}=z_{j}$ for 
$i, j=1, \dots, n$, $i\ne j$. 

\item The associativity of intertwining operators
among grading-restricted left $V_{0}$ modules holds. 

\end{enumerate}
Let $W$ be a grading-restricted left $V$-module and $W_{2}$ a left $V$-submodule of $W$. 
Then for any left $V_{0}$-submodule $W_{1}$ of $W$ such that $W=W_{1}\oplus W_{2}$ as 
a left $V_{0}$-module, Assumption \ref{composability} holds for the map 
$\pi_{W_{2}}\circ Y_{W}\circ (1_{V}\otimes \pi_{W_{1}})$. 
\end{thm}
\pf
Since $Y_{W_{1}}$ and $Y_{W_{2}}$ are vertex operator maps for left $V$-modules, 
they are also intertwining operators among left $V_{0}$-modules. By Proposition \ref{int-op},
$\pi_{W_{2}}\circ Y_{W}\circ (1_{V}\otimes \pi_{W_{1}})$ is also an intertwining operator 
among left $V_{0}$-modules. 

Then by Condition 1,
(\ref{prod}) for $v, v_{1}, \dots, v_{k+l}\in V$, $w_{2}'\in W_{2}'$ and $w_{1}\in W_{1}$
are absolutely convergent in the region $|z_{1}|>\cdots>|z_{k}|>|z|>|z_{k+1}|>\cdots>|z_{k+l}|>0$ 
and can be analytically extended to a (possibly multivalued) analytic function with the only possible 
singularities (branch points or poles) $z_{i}=0$ for $i=1, \dots, k+l$ and $z_{i}=z_{j}$ for 
$i, j=1, \dots, k+l$, $i\ne j$. The proof of Theorem 2.2 in \cite{HVir} in fact 
shows that the commutativity of intertwining operators among grading-restricted left $V_{0}$-modules also holds. 
Since the intertwining operators $Y_{W_{1}}$,  $Y_{W_{2}}$ and 
$\pi_{W_{2}}\circ Y_{W}\circ (1_{V}\otimes \pi_{W_{1}})$ involve only integral powers of 
the variables, using the associativity and commutativity of intertwining operators among 
grading-restricted left $V_{0}$-modules,  we see that
 the singularities $z_{i}=0$ for $i=1, \dots, k+l$ and $z_{i}=z_{j}$ for 
$i, j=1, \dots, k+l$, $i\ne j$ must be poles. Thus (\ref{prod}) is absolutely convergent in the region 
$|z_{1}|>\cdots>|z_{k}|>|z|>|z_{k+1}|>\cdots>|z_{k+l}|>0$  to a rational function with the only possible poles
$z_{i}=0$ for $i=1, \dots, k+l$ and $z_{i}=z_{j}$ for 
$i, j=1, \dots, k+l$, $i\ne j$. So the first part of Assumption \ref{composability}
for $\pi_{W_{2}}\circ Y_{W}\circ (1_{V}\otimes \pi_{W_{1}})$ 
holds. The second part of Assumption \ref{composability} holds for 
$\pi_{W_{2}}\circ Y_{W}\circ (1_{V}\otimes \pi_{W_{1}})$ also because of the 
 associativity and commutativity of intertwining operators among 
grading-restricted left $V_{0}$-modules.
\epfv

We also have the following result:

\begin{thm}\label{F(v)-g-r}
Let $W$ be a grading-restricted left $V$-module,  $W_{2}$ a left $V$-submodule of $W$ and
$W_{1}$ a left $V_{0}$-submodule  of $W$ such that $W=W_{1}\oplus W_{2}$ as 
a left $V_{0}$-module. Assume that for some $z_{0}\in \C^{\times}$, the $Q(z_{0})$-tensor product 
$W_{2}'\boxtimes_{Q(z_{0})} W_{1}$ of left $V_{0}$-modules $W_{2}'$ and $W_{1}$ 
 is still a  grading-restricted 
generalized left $V_{0}$-module.
Then for the map $F: V\to H^{N}$ defined in 
Proposition \ref{F(v)},  the graded vector space $F(V)$ is grading restricted and, in particular,
the $V$-bimodule $H^{F(V)}$ is generated by the grading-restricted subspace $F(V)$.
\end{thm}
\pf
We consider the category $\mathcal{C}$ of grading-restricted 
left $V_{0}$-modules. By Proposition 5.69 in \cite{HLZ1.5}, since 
$W_{2}'\boxtimes_{Q(z_{0})} W_{1}$ is still in $\mathcal{C}$, $W_{2}'\hboxtr_{Q(z_{0})} W_{1}$ is also
in $\mathcal{C}$. In particular, $W_{2}'\hboxtr_{Q(z_{0})} W_{1}$ is grading restricted. 

For  $p\in \Z$, by Remark 5.62 and Proposition 5.63 in \cite{HLZ1.5},
$$(I_{\pi_{W_{2}}\circ Y_{W}\circ (1_{V}\otimes \pi_{W_{1}}), p}^{Q(z_{0})})'(V)
\subset W_{2}'\hboxtr_{Q(z_{0})} W_{1}.$$ 
Since $W_{2}'\hboxtr_{Q(z_{0})} W_{1}$ is grading restricted,
$(I_{\pi_{W_{2}}\circ Y_{W}\circ (1_{V}\otimes \pi_{W_{1}}), p}^{Q(z_{0})})'(V)$ is also grading restricted. 
On the other hand,  for $v\in V$, $w_{2}'\in W_{2}'$ and $w_{1}\in W_{1}$, 
\begin{align*}
&((I_{\pi_{W_{2}}\circ Y_{W}\circ (1_{V}\otimes \pi_{W_{1}}), p}^{Q(z_{0})})'(v))(w_{2}'\otimes w_{1})\nn
&\quad\quad =\langle w_{2}', \pi_{W_{2}}Y_{W}(v, z_{0})\pi_{W_{1}}w_{1}\rangle\nn
&\quad\quad=\langle w_{2}', ((F(v))(w_{1}))(z_{0})\rangle
\end{align*}
By the $D$-derivative property for $Y_{W}$, 
$((F(v))(w_{1}))(z)=((F(e^{(z-z_{0})D_{V}}v))(w_{1}))(z_{0})$,
that is, $F(v)$ is determined completely by its evaluation at $z_{0}$.
So we obtain a linear map 
from $(I_{\pi_{W_{2}}\circ Y_{W}\circ (1_{V}\otimes \pi_{W_{1}}), p}^{Q(z_{0})})'(V)$
to $H^{F(V)}$ given by 
$$((I_{\pi_{W_{2}}\circ Y_{W}\circ (1_{V}\otimes \pi_{W_{1}}), p}^{Q(z_{0})})'(v))
\mapsto F(v).$$ 
It is clear from the definition that this map is a linear isomorphism and 
preserves weights. Since $(I_{\pi_{W_{2}}\circ Y_{W}\circ (1_{V}\otimes \pi_{W_{1}}), p}^{Q(z_{0})})'(V)$ 
is grading restricted, the space $F(V)$ is also grading restricted. 
\epfv

\begin{rema}
{\rm Since $W_{1}$ is a $V_{0}$-submodule of $W$, it is 
invariant under the action of $D_{W}$. In particular, $D_{W_{1}}$ is equal to the 
restriction of $D_{W}$ on $W_{1}$. From this fact, the definition of $F(v)$ for $v\in V$ and the 
$D$-derivative property of the vertex operator $Y_{W}(v, z)$, we  see that 
$F(v)$ satisfies the 
the $D$-derivative property
$$\frac{d}{dz}((F(v))(w_{1}))(z)=D_{W}((F(v))(w_{1}))(z)
-((F(v))(D_{W_{1}}w_{1}))(z).$$
Thus elements of $F(V)$ satisfy the $D$-derivative property.}
\end{rema}

From Theorems \ref{va}, \ref{va-composability} and \ref{F(v)-g-r}, we obtain immediately our main theorem 
in this section:

\begin{thm}
Let $V$ be a grading-restricted vertex algebra and 
$V_{0}$ a vertex subalgebra of $V$. Assuming that the two conditions in Theorem 
 \ref{va-composability} hold. Also assume that every  grading-restricted left $V_{0}$-module  is 
completely reducible and  the $Q(z)$-tensor product of any two grading-restricted left $V_{0}$-modules
is still a grading-restricted left $V_{0}$-module. Then $\hat{H}_{\infty}^{1}(V, M)=\hat{Z}_{\infty}^{1}(V, M)$ 
for every $\Z$-graded $V$-bimodule $M$ generated by a grading-restricted subspace 
implies the complete reducibility of every grading-restricted
left $V$-module. 
\end{thm}

The two conditions in Theorem \ref{va-composability} and the condition in Theorem \ref{F(v)-g-r}
indeed hold when $V_{0}$ and $V_{0}$-modules
satisfy some algebraic conditions. We discuss these conditions in the following remark:

\begin{rema}
{\rm A left $V_{0}$-module $W$ is $C_{1}$-cofinite if $\dim W/C_{1}(W)<\infty$ 
where $C_{1}(W)$ is the subspace 
of $W$ spanned by elements of the form $(Y_{W})_{-1}(v)w$ for $v\in V_{0}$ and $w\in W$.
Note that in \cite{Hdiff-eqn}, a generalized module for a vertex (operator) algebra means a 
weak module graded by the eigenvalues of $L(0)$ but the grading does not have to be lower bounded 
and the homogeneous subspaces do not have to be finite dimensional.  Here we 
call such a generalized module an $L(0)$-semisimple left weak module. 
By Theorems 3.5 and 3.7 in \cite{Hdiff-eqn},  the assumptions
in Theorems \ref{va-composability} and \ref{F(v)-g-r} hold when the following three conditions
on (generalized) $V_{0}$-modules hold:
\begin{enumerate}

\item Every $L(0)$-semisimple left weak $V_{0}$-module is a direct sum of irreducible
grading-restricted left $V_{0}$-modules.

\item There are only finitely many inequivalent irreducible grading-restricted left $V_{0}$-modules
and they are all $\R$-graded.

\item Every irreducible grading-restricted left $V_{0}$-module satisfies the $C_{1}$-cofiniteness condition.

\end{enumerate}
The vertex algebra $V_{0}$ is $C_{2}$-cofinite if $\dim V_{0}/C_{2}(V_{0})<\infty$ 
where $C_{2}(V_{0})$ is the subspace 
of $V_{0}$ spanned by elements of the form $(Y_{V_{0}})_{-2}(u)v$ for $u, v\in V_{0}$.
By Theorem 3.9 in \cite{Hdiff-eqn}, the three conditions above hold when the following three conditions holds:
\begin{enumerate}

\item For $n < 0$, $(V_{0})_{(n)} = 0$ and $(V_{0})_{(0)} = \C\one$.

\item Every grading-restricted left $V_{0}$-module is completely reducible.

\item $V_{0}$ is $C_{2}$-cofinite.

\end{enumerate}
By Proposition 12.5 in \cite{HL},  Corollary 9.30 in \cite{HLZ2}, Theorem 11.8 in \cite{HLZ3} 
and Theorem 3.1 in \cite{Hz-correction}, the 
assumptions in Theorems \ref{va-composability} and \ref{F(v)-g-r} also hold when the following three conditions
hold:
\begin{enumerate}

\item Every grading-restricted left $V_{0}$-module is completely reducible.

\item  For any left $V_{0}$-modules $W_{1}$ and $W_{2}$ and any $z\in \C^{\times}$, if 
the weak $V$-module $W_{\lambda}$ 
generated by a generalized eigenvector  $\lambda\in 
(W_{1}\otimes W_{2})^{*}$ for 
$L_{P(z)}(0)$ satisfying the $P(z)$-compatibility condition is lower bounded, then $W_{\lambda}$
is a grading-restricted left $V_{0}$-module.

\item Every irreducible grading-restricted left $V_{0}$-module satisfies the $C_{1}$-cofiniteness condition.
\end{enumerate}
Thus the conclusions of Theorems \ref{va-composability} and \ref{F(v)-g-r} 
hold when one of these three sets of three conditions holds.}
\end{rema}

\noindent {\small \sc Max Planck Institute for Mathematics, Vivatsgasse 7
53111 Bonn, Germany}
\vspace{.5em}

\noindent {\it and}
\vspace{.5em}

\noindent {\small \sc Department of Mathematics, Rutgers University,
110 Frelinghuysen Rd., Piscataway, NJ 08854-8019, USA  (permanent address)}
\vspace{.5em}

\noindent {\em E-mail addresses}: yzhuang@math.rutgers.edu
\vspace{1em}

\noindent {\small \sc Department of Mathematics, Rutgers University,
110 Frelinghuysen Rd., Piscataway, NJ 08854-8019, USA}
\vspace{.5em}

\noindent {\em E-mail addresses}: volkov\_cl@hotmail.com

\end{document}